\documentclass[12pt]{article}
\usepackage{amsmath, amsfonts, amscd, latexsym, amsthm, amssymb, graphicx}
\newtheorem{theor}{\hspace{1cm}{\sc Theorem}}[section]
\newtheorem{utver}[theor]{\hspace{1cm}{\sc Proposition}}
\newtheorem{sledst}[theor]{\hspace{1cm}{\sc Corollary}}
\newtheorem{lemma}[theor]{\hspace{1cm}{\sc Lemma}}
\newtheorem{conj}[theor]{\hspace{1cm}{\sc Conjecture}}
\newtheorem*{utver*}{\hspace{1cm}{\sc Proposition}}
\theoremstyle{definition}

\newtheorem{defin}[theor]{\hspace{1cm}{\sc Definition}}
\newtheorem{exa}[theor]{\hspace{1cm}{\sc Example}}
\newtheorem*{rem}{\hspace{1cm}{\sc Remark}}

\newcommand{\Vol}{\mathop{\rm Vol}\nolimits}

\newcommand{\codim}{\mathop{\rm codim}\nolimits}

\newcommand{\id}{\mathop{\rm id}\nolimits}

\newcommand{\spec}{\mathop{\rm spec}\nolimits}

\newcommand{\conv}{\mathop{\rm conv}\nolimits}

\newcommand{\BP}{\mathop{\rm BP}\nolimits}
\newcommand{\MV}{\mathop{\rm MV}\nolimits}
\newcommand{\MF}{\mathop{\rm MP}\nolimits}
\newcommand{\HP}{\mathop{\rm HP}\nolimits}
\newcommand{\MS}{\mathop{\rm MS}\nolimits}

\def\R{\mathbb R}
\def\N{\mathbb N}
\def\Z{\mathbb Z}
\def\Q{\mathbb Q}
\def\C{\mathbb C}
\def\CC{({\mathbb C}\setminus 0)}
\def\T{\mathbb T}

\def\CP{\mathbb C\mathbb P}

\begin{document}

\begin{center}
{\Large \textbf{Newton polyhedra of discriminants of
projections.}}

A. Esterov\footnote{Partially supported by RFBR-JSPS-06-01-91063,
RFBR-07-01-00593 and INTAS-05-7805 grants. Email address:
\texttt{esterov@mccme.ru}}
\end{center}


\begin{quote} For a system of polynomial equations, whose coefficients
depend on parameters, the Newton polyhedron of its discriminant
is computed in terms of the Newton polyhedra of the coefficients.
This leads to an explicit formula (involving Euler obstructions
of toric varieties) in the unmixed case, suggests certain open
questions in general, and generalizes similar known results
(\cite{gkz}, \cite{sturmf}, \cite{mcd}, \cite{gp}, \cite{ekh}).
\end{quote}

{\large \textbf{Introduction.}}

Let $F_0,\ldots,F_l$ be Laurent polynomials on the complex torus
$\CC^k$, whose coefficients are Laurent polynomials on the
parameter space $\CC^n$. Consider the set $\Sigma\subset\CC^n$ of
all values of the parameter, such that the corresponding system
of polynomial equations $F_0=\ldots=F_l=0$ defines a singular set
in $\CC^k$. In most cases (see below for details), the closure of
$\Sigma$ is a hypersurface, and its defining equation is called
the \textit{discriminant} of $F_0=\ldots=F_l=0$.

In this paper, we compute the Newton polyhedron of the
discriminant in terms of Newton polyhedra of the coefficients of
the polynomials $F_0,\ldots,F_l$. The answer is known in many
special cases, and we give a number of references as examples of
various approaches to this problem: the universal special case
for $l=0$ and $l=k$ was studied in \cite{gkz}, \cite{sturmf} and
\cite{dfs} (\textit{universal case} means that $\CC^n$
parameterizes all collections of polynomials $F_0,\ldots,F_l$,
whose monomials are contained in a given finite set of
monomials), the general case for $l=k-1$, for $l=k-1=0$ and for
$l=k$ was studied in \cite{mcd}, \cite{gp} and \cite{ekh}.

To formulate the answer in general, we need the following
notation: we denote the Minkowski sum $\{a+b\; |\; a\in A,\; b\in
B\}$ of polyhedra $A$ and $B$ by $A+B$, and denote the mixed fiber
polyhedron of the polyhedra $\Delta_0,\ldots,\Delta_l$ in
$\R^n\oplus\R^k$ by the monomial $\Delta_0\cdot\ldots\cdot
\Delta_l$ (it is a certain polyhedron in $\R^k$, see \cite{mcm},
\cite{ekh},  Definition \ref{defrelmixfib}, or Appendix). To a set
$A\subset\Z^k$ and a face $B$ of its convex hull we associate its
Euler obstruction $e^{B,A}\in\Z$, whose combinatorial definition
is given in Subection \ref{seulobstr}, and whose geometrical
meaning is $(-1)^{k-\dim B}$ times the Euler obstruction of the
$A$-toric variety at its orbit, corresponding to $B$ (see also
\cite{tak} and a remark at the end of Subsection
\ref{snewtdiscrim}).

Considering $F_i$ as a polynomial on $\CC^k$ with polynomial
coefficients, denote the set of its monomials by $A_i\subset\R^k$;
considering the same $F_i$ as a polynomial on $\CC^n\times\CC^k$
with complex coefficients, denote its Newton polyhedron by
$\Delta_i\subset\R^{n}\oplus\R^{k}$. Denote the preimage of a set
$A'$ under the natural projection $\Delta_i\to \mbox{(convex hull
of } A_i)$ by $\Delta_i(A')$. For simplicity, we assume here that
$A_0=\ldots=A_l=A$, and pairwise differences of its points
generate $\Z^k$.

\vspace{3mm}

\textsc{Theorem.} \textit{If $F_0,\ldots,F_l$ are generic
polynomials 
with Newton polyhedra $\Delta_0,\ldots,\Delta_l$, then the Newton
polyhedron of the discriminant equals}

$$\mathcal{N}\; =\; \sum_{A'\subset A}
\qquad e^{A',A}\quad\cdot \sum_{a_0>0,\ldots,a_l>0\atop
a_0+\ldots+a_l=\dim A'+1}
\Delta_0(A')^{a_0}\cdot\ldots\cdot\Delta_l(A')^{a_l},$$
\textit{where $A'$ runs over all faces of dimension $l$ and
greater in the convex hull of $A$.}

\textit{More precisely, if the polyhedron $\mathcal{N}$ consists
of one point, then the discriminant set $\Sigma$ has no
codimension 1 components; otherwise, the closure of $\Sigma$ is a
hypersurface, and the Newton polyhedron of its equation equals
$\mathcal{N}$.}

\vspace{3mm}

The word ``generic" means that the set of collections
$F_0,\ldots,F_l$, satisfying the statement, is dense in the space
of collections with the given Newton polyhedra
$\Delta_0,\ldots,\Delta_l$. Note that coefficients $e^{A',A}$ may
be negative, and the formula above involves subtraction of
polyhedra. The difference of polyhedra $P$ and $Q$ is by
definition the solution of the equation $Q+X=P$, which is always
unique if exists; see the end of Subsection \ref{losungex} for
details and related computability questions.

The result that we actually prove is somewhat more general than
the theorem stated above in the following sense. Firstly, together
with the Newton polyhedron, we describe the leading coefficients
of the discriminant (i.e. the coefficients of monomials in the
boundary of the Newton polyhedron) in terms of leading
coefficients of the polynomials $F_0,\ldots,F_l$ (see Theorem
\ref{elimth} for $l=k$, Proposition \ref{thleadcoefdiscr} for
$l=0$, and Theorem \ref{cayley} that reduces the general case to
$l=0$). Secondly, we do not assume that $A_0=\ldots=A_l$.
Thirdly, we prefer a slightly more general context throughout the
paper (see \cite{gp} for motivation): instead of polynomials on
$\CC^n$, we consider analytic functions on an arbitrary affine
toric variety. Nevertheless, all our results and proofs can be
translated back into the global setting, word by word,
substituting germs of analytic functions on affine toric
varieties with Laurent polynomials on complex tori, unbounded
polyhedra with bounded ones, and the local version of the
elimination theorem \ref{elimth} with its global version
\cite{ekh}. For instance, the theorem stated above is exactly the
global version of Theorem \ref{losung} with $A_0=\ldots=A_l$. We
illustrate it with an example in Subsection \ref{losungex}.

Denote the Newton polyhedron of the discriminant of
$F_0=\ldots=F_l=0$ by $\mathcal{N}(\Delta_0,\ldots,\Delta_l)$ for
generic equations $F_0,\ldots,F_l$ with given Newton polyhedra
$\Delta_0,\ldots,\Delta_l$. Theorem \ref{losung} leads to a
certain relation between Newton polyhedra of discriminants
(higher additivity, see Subsection \ref{snewtdiscrci3} for
details), provided that the convex hulls of $A_0,\ldots,A_l$ are
large enough and analogous (i.e. have the same number of faces as
the convex hull of their sum $A_0+\ldots+A_l$):
$$\mathcal{N}(\Delta_0+\Delta_1,\Delta_2,\ldots,\Delta_l) =
\sum_{\mu=1}^{\infty}
\mathcal{N}(\underbrace{\Delta_0,\ldots,\Delta_0}_{\mu},\underbrace{
\Delta_1,\ldots,\Delta_1}_{\mu-1},\Delta_2,\ldots,\Delta_l)+$$
$$+\mathcal{N}(\underbrace{\Delta_0,\ldots,\Delta_0}_{\mu-1},\underbrace{
\Delta_1,\ldots,
\Delta_1}_{\mu},\Delta_2,\ldots,\Delta_l)+2\mathcal{N}(\underbrace{\Delta_0,\ldots,\Delta_0}_{\mu},\underbrace{
\Delta_1,\ldots, \Delta_1}_{\mu},\Delta_2,\ldots,\Delta_l)$$ (all
but finitely many polyhedra in this sum are equal to $\{0\}$; in
particular, for $l=k+1$, this is the conventional additivity).
Unexpectedly, the assumption that $A_0,\ldots,A_l$ are analogous
can be significantly relaxed in some cases (see Subsection
\ref{snewtdiscrci3} for examples), and it would be interesting to
know, to what extent it can be relaxed in general.

\vspace{3mm}

The paper is organized as follows. In Section 1, we recall
necessary facts and notation, related to convex geometry and
Newton polyhedra. In section 2, we study the universal case of
our problem, which generalizes results of \cite{gkz} and
\cite{sturmf}. In Section 3, we study the special case $l=k$ of
our problem (Theorem \ref{elimth}), which is a local version of
elimination theory in the context of Newton polyhedra \cite{ekh},
and is based on a certain local version of D.~Bernstein's formula
(Theorem \ref{relbernst}). In Section 4, we apply elimination
theory to study the special case $l=0$ (Theorem \ref{mainth}). In
Section 5, we reduce the general case $0\leqslant l\leqslant k$
to the case $l=0$ by means of a classical technique, known as
Cayley trick, or Lagrange multipliers (Theorem \ref{cayley}).

In particular, if the Newton polyhedra of $F_0,\ldots,F_l$ are not
too ``thin", then the discriminant set $\Sigma$ is a hypersurface
(see Propositions \ref{gm11} for $l=0$ or \ref{gm21} for arbitrary
$l$), and its Newton polyhedron can be computed by Theorems
\ref{mainth} and \ref{cayley}. By ``thin" we mean Newton
polyhedra, such that the collection $A_0,\ldots,A_l$ (see above)
is dual defect. One simple test for non-dual-defectiveness is
provided by Propositions \ref{pthin} (for $l=0$) and
\ref{equidim}.

We also study the same problem as formulated in the beginning,
with another definition of the discriminant: we can define the
discriminant set as the minimal set $S$ in $\CC^n$, such that the
restriction of the projection $\CC^n\times\CC^k\to\CC^n$ to the
complete intersection $\{F_0=\ldots=F_l=0\}$ is a fiber bundle
outside of $S$. This version of the problem is outlined in
Subsection \ref{sfibration}. For example, if $A$ is the set of
integer points of a Delzant polyhedron in the assumptions of the
theorem stated above, then the Newton polyhedron of the equation
of $S$ equals
$$\sum_{a_0>0,\ldots,a_l>0\atop a_0+\ldots+a_l=k+1}
\Delta_0^{a_0}\cdot\ldots\cdot\Delta_l^{a_l}$$ (see Corollary
\ref{losung2}). The counterpart of dual-defectiveness for this
problem seems to behave much simpler: see Proposition
\ref{conj1sim} and Conjecture \ref{conj1}.

All answers are formulated in terms of mixed fiber polyhedra.
This notion is introduced in Section 1.2 and, in more detail, in
Appendix (existence, uniqueness and monotonicity of mixed fiber
bodies is proved and a formula for the support function of a
mixed fiber body is given).

I am grateful to Askold Khovanskii, Yutaka Matsui, Kiyoshi
Takeuchi, Pedro Gonz\'alez-P\'erez and Semion Tregub for fruitful
discussions. \vspace{-12mm} \small\tableofcontents\normalsize

\section{Mixed volumes, mixed fiber polyhedra, and Euler obstructions.}

In this section, we recall relevant facts and notation from convex
geometry: mixed fiber polyhedra (\cite{mcm}, \cite{ekh}), Euler
obstructions of polyhedra (\cite{tak}), relative mixed volume
(\cite{me05}, \cite{me06}, \cite{mearx}), and the corresponding
relative version of the Kouchnirenko-Bernstein formula
(\cite{bernst}, \cite{kh77}). Subsections \ref{skbk} and
\ref{smultci} contain generalizations of the relative
Kouchnirenko-Bernstein formula that provide a simple proof for
the Matsui-Takeuchi formula for Euler obstructions and for the
Gelfand-Kapranov-Zelevinsky decomposition formula; we do not need
these generalizations for other purposes. Material from other
subsections is used in the proof of our main result.

\subsection{Relative mixed volume.}\label{smixed}

\textsc{Classical mixed volume.} Recall the notion of the mixed
volume of bounded polyhedra. The set $\mathcal{M}$ of all bounded
polyhedra in $\R^m$ is a semigroup with respect to Minkowski
summation $P+Q=\{p+q\; |\; p\in P,\; q\in Q\}$.
\begin{defin} \label{defmixvol}
\textit{The mixed volume} of polyhedra is the symmetric
multilinear function
$\MV:\underbrace{\mathcal{M}\times\ldots\times
\mathcal{M}}_m\to\R$, such that $\MV(P,\ldots,P)$ equals the
volume of $P$ for every $P\in\mathcal{M}$.
\end{defin}
\begin{lemma}[\cite{kh77}]\label{zeromv} ${\MV(\Delta_1,\ldots,\Delta_m)=0}$ if and only if
$\dim\Delta_{i_1}+\ldots+\Delta_{i_J}<J$ for some
$i_1<\ldots<i_J$.
\end{lemma}
This fact is mentioned as obvious in \cite{kh77}, but we prefer
to give a proof for the sake of completeness.

\textsc{Proof.} $(\Leftarrow)$ follows by an explicit computation
if $\Delta_{i_1}=\ldots=\Delta_{i_J}=C$ is a cube, and all other
$\Delta_i$ are equal to an $m$-dimensional cube with a face $C$;
the general case can be reduced to this one by monotonicity of
the mixed volume. In the other direction, consider points
$a_i\in\Delta_i$ and $b_i\in\Delta_i$ such that the vectors
$a_1-b_1,\ldots,a_m-b_m$ are in general position in the sense
that the dimension of the space generated by
$a_{i_1}-b_{i_1},\ldots,a_{i_J}-b_{i_J}$ is the maximal possible
one for every subset $\{i_1,\ldots,i_J\}\subset\{1,\ldots,m\}$.
By monotonicity of the mixed volume, the mixed volume of the
segments, connecting $a_i$ and $b_i$, equals zero, which means
that the vectors $a_1-b_1,\ldots,a_m-b_m$ are linearly dependent.
In particular, there exists a minimal subset
$\{i_1,\ldots,i_J\}\subset\{1,\ldots,m\}$ such that the vectors
$a_{i_1}-b_{i_1},\ldots,a_{i_J}-b_{i_J}$ are linearly dependent.
They generate a proper subspace $L\subset\R^m$, and every $J-1$
of them form a basis of $L$. If there exists
$b'_{i_j}\in\Delta_{i_j}$ such that $a_{i_j}-b'_{i_j}\notin L$,
then the vectors $a_{i_1}-b_{i_1},\ldots,a_{i_J}-b_{i_J}$ with
$a_{i_j}-b'_{i_j}$ instead of $a_{i_j}-b_{i_j}$ generate a
subspace $L'\varsupsetneq L$, which contradicts the condition of
general position. Thus, $\Delta_{i_1},\ldots,\Delta_{i_J}$ are
contained in a $(J-1)$-dimensional subspace $L$, up to a parallel
translation, and $\dim\Delta_{i_1}+\ldots+\Delta_{i_J}<J$.
$\quad\Box$

\vspace{0.3cm} \textsc{Relative mixed volume.} We need the
following relative version of the mixed volume. For a convex
polyhedral $m$-dimensional cone $\tau\subset(\R^m)^*$, denote its
dual cone $\{x\in\R^m\; |\; \gamma(x)>0 \mbox{ for }
\gamma\in\tau\}$ by $\tau^{\vee}$, and let
$\mathcal{M}_{\tau^{\vee}}$ be the semigroup of all (unbounded)
polyhedra of the form
\begin{center} $\tau^{\vee}+$ a bounded polyhedron. \end{center} \noindent
Consider the set
$\mathcal{P}_{\tau^{\vee}}\subset\mathcal{M}_{\tau^{\vee}}\times\mathcal{M}_{\tau^{\vee}}$
of all ordered pairs of polyhedra $(P,Q)$, such that the
symmetric difference $P\vartriangle Q$ is bounded.
$\mathcal{P}_{\tau^{\vee}}$ is a semigroup with respect to
Minkowski summation of pairs $(P,Q)+(C,D)=(P+C,Q+D)$.

\begin{defin}[\cite{me05}, \cite{me06}, \cite{mearx}] \label{defrelmixvol}
\textit{The volume} $V(P,Q)$ of a pair of polyhedra
$(P,Q)\in\mathcal{P}_{\Gamma}$ is defined to be the difference
$\Vol(P\setminus Q)-\Vol(Q\setminus P)$. \textit{The mixed
volume} of pairs of polyhedra is defined to be the symmetric
multilinear function
$\MV:\underbrace{\mathcal{P}_{\tau^{\vee}}\times\ldots\times
\mathcal{P}_{\tau^{\vee}}}_m\to\R$, such that
$\MV\bigl((P,Q),\ldots,(P,Q)\bigr)=V(P,Q)$ for every pair
$(P,Q)\in \mathcal{P}_{\tau^{\vee}}$.
\end{defin}
Existence and uniqueness are proved in \cite{me06} and
\cite{mearx}.
\begin{exa} If
$\tau^{\vee}=\{0\}$, then $\mathcal{P}_{\tau^{\vee}}$ is the set
of pairs of convex bounded polyhedra, and the mixed volume of
pairs $\MV\Bigl((P_1,Q_1),\ldots,(P_m,Q_m)\Bigr)$ equals the
difference of classical mixed volumes of the collections
$P_1,\ldots,P_m$ and $Q_1,\ldots,Q_m$.
\end{exa}
In general, the mixed volume of pairs can be expressed in terms of
the classical mixed volume as follows.
\begin{lemma}[\cite{me06}, \cite{mearx}]\label{mvmv} Let $\widetilde P_i\subset P_i$ and $\widetilde Q_i\subset Q_i$
be bounded polyhedra in $\R^m$, such that $\widetilde
P_i\setminus\widetilde Q_i=P_i\setminus Q_i$ and $\widetilde
Q_i\setminus\widetilde P_i=Q_i\setminus P_i$ for $i=1,\ldots,m$.
Then
$$\MV\Bigl( (P_1,Q_1),\ldots,(P_m,Q_m)\Bigr)=\MV(\widetilde P_1,\ldots,\widetilde P_m)-
\MV(\widetilde Q_1,\ldots,\widetilde Q_m).$$\end{lemma} Note
that, for any pair $(P_i,Q_i)$, we can always find the requested
bounded polyhedra $\widetilde P_i$ and $\widetilde Q_i$.

The cone $\tau^{\vee}$ plays the role of the unit in the semigroup
$\mathcal{P}_{\tau^{\vee}}$:
\begin{lemma}[\cite{me06}, \cite{mearx}]\label{mvcone}
$$\mbox{1) } \MV\Bigl( (\tau^{\vee},Q_1),(P_2,Q_2),\ldots,(P_m,Q_m)\Bigr)=\MV\Bigl(
(\tau^{\vee},Q_1),(Q_2,Q_2),\ldots,(Q_m,Q_m)\Bigr),$$

i.e. the left hand side does not depend on the choice of
$P_2,\ldots,P_m$;
\newline 2) $\MV\Bigl( (\tau^{\vee},\tau^{\vee}),(P_2,Q_2),\ldots,(P_m,Q_m)\Bigr) =
0.$
\end{lemma}

\vspace{0.3cm} \textsc{Mixed volume of a prism.} Let
$e_1,\ldots,e_l$ be the standard basis of $\R^l$, and $e_0$ be
$0\in\R^l$. For bounded polyhedra $P_0,\ldots,P_l$ in $\R^m$ and
a subset $I\subset\{0,\ldots,l\}$, denote the convex hull of the
union of the polyhedra $P_i\times\{e_i\}\subset\R^m\oplus\R^l,\;
i\in I,$ by $P_I$. In what follows, it will be convenient to
denote the mixed volume of bounded polyhedra $Q_1,\ldots,Q_m$ in
$\R^m$ by the monomial $Q_1\cdot\ldots\cdot Q_m$.
\begin{lemma}\label{chimv}
$$\sum_{I\subset\{0,\ldots,l\}}(-1)^{l+1-|I|}(m+|I|-1)!\,\Vol(P_I)=
\sum_{a_0>0,\ldots,a_l>0\atop a_0+\ldots+a_l=m}
m!\,P_0^{a_0}\cdot\ldots\cdot P_l^{a_l}.$$
\end{lemma}
\textsc{Proof.} Pick generic polynomials $g_0,\ldots,g_l$ whose
Newton polyhedra are $P_0,\ldots,P_l$. Compute the Euler
characteristic of the hypersurface
$\lambda_0g_0+\ldots+\lambda_lg_l=0$ in
$\CP^l_{\lambda_0:\ldots:\lambda_l}\times\CC^m$ in the following
two ways.\newline 1) The subdivision of the toric variety
$\CP^l_{\lambda_0:\ldots:\lambda_l}\times\CC^m$ into complex tori
$T_I=\{\lambda_i=0\, \mbox{ for } i\notin I\}$ induces the
subdivision of the hypersurface into pieces, whose Euler
characteristics equal $(-1)^{l+1-|I|}(m+|I|-1)!\,\Vol(P_I)$ by the
Kouchnirenko-Khovanskii formula \cite{kh77}. By additivity, the
Euler characteristic of the hypersurface
$\lambda_0g_0+\ldots+\lambda_lg_l=0$ equals the left hand side of
the desired equality.
\newline 2) Considering the projection of the hypersurface
$\lambda_0g_0+\ldots+\lambda_lg_l=0$ to $\CC^m$, we note that the
fiber of this projection over a point $y\in\CC^m$ equals
$\CP^{l}$ or $\CP^{l-1}$ depending on whether
$y\in\{g_0=\ldots=g_l=0\}$ or not. Thus, integrating the Euler
characteristic over fibers of this projection, we conclude that
the Euler characteristic of the hypersurface
$\lambda_0g_0+\ldots+\lambda_lg_l=0$ equals the Euler
characteristic of the complete intersection
$\{g_0=\ldots=g_l=0\}$. Computing the latter one by the
Kouchnirenko-Bernstein-Khovanskii formula \cite{kh77}, we get the
right hand side of the desired equality. $\quad\Box$

\subsection{Mixed fiber polyhedra.}\label{smixedfib}

\textsc{Minkowski integral $\int \Delta$} (\cite{bs}). Denote the
projections of the direct sum $\R^n\oplus\R^k$ onto the summands
by $p$ and $q$ respectively. Consider an integer polyhedron
$\Delta\subset\R^n\oplus\R^k$, whose projection $q(\Delta)$ is
bounded, and denote the affine span of $q(\Delta)$ by $S$. Choose
the volume form ${\rm d}x$ on $S$ such that the volume of the
image of $S$ under the projection $\R^k\to\R^k/\Z^k$ equals
$(1+\dim S)!$, and consider the set $I\subset\R^n\oplus\R^k$ of
points of the form $\int_{q(\Delta)}\, s(x)\, {\rm d}x$, where $s$
runs over all continuous sections of the projection $q:\Delta\to
q(\Delta)$.

\begin{defin}\label{integr} The \textit{Minkowski integral}, or the \textit{fiber polytope} $\int \Delta$ is the closure of $p(I)$.
\end{defin}

\textsc{Mixed Minkowski integral} (\cite{mcm}). Let
$\tau^{\vee}\subset\R^n\subset\R^n\oplus\R^k$ be a convex
polyhedral cone that does not contain a line. \textit{The mixed
Minkowski integral, or the mixed fiber polyhedron,} is the
symmetric multilinear polyhedral-valued function
$\MF:\underbrace{\mathcal{M}_{\tau^{\vee}}\times\ldots\times
\mathcal{M}_{\tau^{\vee}}}_{k+1}\to\mathcal{M}_{\tau^{\vee}},\,$
such that, for every polyhedron $\Delta\in
\mathcal{M}_{\tau^{\vee}}$, $$\MF(\Delta,\ldots,\Delta)=\left\{
\begin{array}{ll}\int \Delta \; &\mbox{ if } \dim q(\Delta)=n \\
\tau^{\vee}&\mbox{ otherwise }\end{array}\right. .$$ See
\cite{mcm}, \cite{ekh}, or Appendix for existence, uniqueness and
properties of the mixed fiber polyhedron.

\begin{exa}\label{examf} If $\Delta$ is the product of $P\subset\R^n$ and
$Q\subset\R^k$, then $\int \Delta$ equals $\Bigl((n+1)!\Vol
Q\Bigr)\cdot P$. If $n=1$ and $\tau^{\vee}=\{0\}$, then
$\MF(\Delta_0,\ldots,\Delta_n)$ is a segment of length
$(n+1)!\MV(\Delta_0,\ldots,\Delta_n)$.
\end{exa}


We need the following formula for the support function of the
mixed fiber polyhedron. Let $\Delta_0,\ldots,\Delta_k$ be
polyhedra in $\mathcal{M}_{\tau^{\vee}}$, denote the set of
positive real numbers by $\R_+\subset\R$, and the product
$\R_+\times q(\Delta_i)\subset\R\oplus\R^k$ by $B_i$. For every
linear function $l:\R^n\to\R$, denote the image of $\Delta_i$
under the map $(l,\id):\R^n\oplus\R^k\to\R\oplus\R^k$ by
$l\Delta_{i}$.
\begin{utver}[see Appendix]\label{supportmixfib} The minimal value of a linear function $l:\R^n\to\R$ on
$\MF(\Delta_0,\ldots,\Delta_k)$ equals
$$(k+1)!\MV\Bigl((l\Delta_0,\; B_0),\; \ldots,\; (l\Delta_k,\; B_k)\Bigr)$$ if
$l\in\tau$ and equals $-\infty$ otherwise.
\end{utver}
The proof is given in Appendix under an inessential assumption
that the polyhedra are bounded.

We also need a little more flexible version of the notation
above. For an $l$-dimensional vector space $L\subset\R^k$,
consider the semigroup $\mathcal{M}_{\tau^{\vee}}(L)$ of all
polyhedra of the form $$Q +
\tau^{\vee}\times\{x\}\subset\R^n\oplus\R^k,$$ where $x$ is a
point in $\R^k$ and $Q$ is a bounded polyhedron in $\R^n\oplus L$.
\begin{defin} \label{defrelmixfib}
\textit{The mixed Minkowski integral, or the mixed fiber
polyhedron,} is the symmetric multilinear polyhedral-valued
function
$$\MF:\underbrace{\mathcal{M}_{\tau^{\vee}}(L)\times\ldots\times
\mathcal{M}_{\tau^{\vee}}(L)}_{l+1}\to\mathcal{M}_{\tau^{\vee}}(0),\,$$
such that, for every polyhedron $\Delta\in
\mathcal{M}_{\tau^{\vee}}(L)$, $$\MF(\Delta,\ldots,\Delta)=\left\{
\begin{array}{ll}\int \Delta \; &\mbox{ if } \dim q(\Delta)=l \\
\tau^{\vee}&\mbox{ otherwise }\end{array}\right. .$$
\end{defin}
Note that the value $\MF(\Delta_0,\ldots,\Delta_l)$ does not
depend on the choice of the cone $\tau^{\vee}$ and the space $L$
in the definition above. More precisely, the cone $\tau^{\vee}$ is
uniquely defined by the arguments $\Delta_0,\ldots,\Delta_l$. The
space $L$ is also uniquely defined by $\Delta_0,\ldots,\Delta_l$,
provided that the projection $q(\Delta_0+\ldots+\Delta_l)$ is
$l$-dimensional; otherwise,
$\MF(\Delta_0,\ldots,\Delta_l)=\tau^{\vee}$ independently of the
choice of $L$.

\vspace{0.3cm} \textsc{Minkowski integral of a prism.} In what
follows, it is convenient to denote the mixed fiber polyhedron
$\MF(\Delta_0,\ldots,\Delta_l)$ by the monomial
$\Delta_0\cdot\ldots\cdot\Delta_l$, as well as we do for the
mixed volume (this agrees with Example \ref{examf}). Let
$e_1,\ldots,e_l$ be the standard basis of $\R^l$, and $e_0$ be
$0\in\R^l$. For polyhedra $P_0,\ldots,P_l$ in
$\mathcal{M}_{\tau^{\vee}}$ and a subset
$I\subset\{0,\ldots,l\}$, denote the convex hull of the union of
the polyhedra $P_i\times\{e_i\}\subset\R^n\oplus\R^k\oplus\R^l,\;
i\in I,$ by $P_I$.
\begin{lemma}\label{chimf} For polyhedra $P_0,\ldots,P_l$ in
$\mathcal{M}_{\tau^{\vee}}$,
$$1)\quad P_{\{0,\ldots,l\}}^{k+l+1}=
\sum_{a_0\geqslant 0,\ldots,a_l\geqslant 0\atop
a_0+\ldots+a_l=k+1} P_0^{a_0}\cdot\ldots\cdot P_l^{a_l},$$
$$2)\quad \sum_{I\subset\{0,\ldots,l\}}(-1)^{l+1-|I|}P_I^{k+|I|}=
\sum_{a_0>0,\ldots,a_l>0\atop a_0+\ldots+a_l=k+1}
P_0^{a_0}\cdot\ldots\cdot P_l^{a_l}.$$
\end{lemma}
\textsc{Proof.} These two formulas are equivalent by the
inclusion-exclusion formula, and we prove the second one.
Substituting mixed fiber polyhedra with mixed volumes of pairs by
Proposition \ref{supportmixfib}, and then with classical mixed
volumes by Lemma \ref{mvmv}, it is enough to prove the same
formula for mixed volumes of bounded polyhedra, which is the
statement of Lemma \ref{chimv}. $\quad\Box$

\subsection{Kouchnirenko-Bernstein formula.}\label{srelbernst}

The relative version of the classical mixed volume participates
in a certain relative version of the classical
Kouchnirenko-Bernstein formula. To formulate it, we need some
notation related to toric varieties, Newton polyhedra and
intersection numbers.

\vspace{0.3cm} \textsc{Toric varieties.} For a rational fan
$\Sigma$ in $(\R^m)^*$, the corresponding toric variety is
denoted by $\T^{\Sigma}$. For every codimension 1 orbit $T$ of
$\T^{\Sigma}$, the primitive generator of the corresponding
1-dimensional cone of $\Sigma$ is denoted by $\gamma(T)$. We
assume that the union of cones of $\Sigma$ is a closed convex cone
$\tau$, and denote its dual by $\tau^{\vee}\subset\R^m$.

If $I$ is a very ample line bundle on $\T^{\Sigma}$, and a
meromorphic section $s$ of the bundle $I$ has no zeros and no
poles in the maximal torus of $\T^{\Sigma}$, then there exists a
unique polyhedron $\Delta\in\mathcal{M}_{\tau^{\vee}}$,
such that the multiplicity of every codimension 1 orbit $T$ of
$\T^{\Sigma}$ in the divisor of poles and zeroes of the section
$s$ equals the maximal value of the linear function
$-\gamma(T):\R^m\to\R$ on the polyhedron ${\Delta}$.
Since the pair $(I,s)$ is uniquely determined by this polyhedron
$\Delta$, we denote the line bundle $I$ by $I_{\Delta}$ and the
section $s$ by $s_{\Delta}$.

\vspace{0.3cm} \textsc{Newton polyhedra.} The union of all
precompact orbits of the toric variety $\T^{\Sigma}$ (i.e. the
orbits, corresponding to the cones of $\Sigma$ in the interior of
$\tau$) is denoted by $\T^{\Sigma}_{comp}$ and is called the
\textit{compact part of} $\T^{\Sigma}$ (it is indeed a compact
set).

If $f$ is an arbitrary germ of a holomorphic section of
$I_{\Delta}$ near the compact set $\T^{\Sigma}_{comp}$, then the
function $f/s_{\Delta}$ can be represented as a power series
$\sum_{a\in\Delta} c_ax^a$ for $x$ in the maximal torus $\CC^m$
of the toric variety $\T^{\Sigma}$. The convex hull of the set
$\{a\; |\; c_a\ne 0\}+\tau^{\vee}$ is an integer polyhedron in
$\mathcal{M}_{\tau^{\vee}}$. It is called the \textit{Newton
polyhedron} of $f$ and is denoted by $\Delta_f$. For any bounded
$\Gamma\subset\R^m$, the polynomial $\sum_{a\in\Gamma} c_ax^a$ on
$\CC^m$ is denoted by $f^{\Gamma}$. If $a$ is contained in a
bounded face of the Newton polyhedron $\Delta_f$, then the
coefficient $c_a$ is called a \textit{leading coefficient} of
$f$. Every section has finitely many leading coefficients.

\vspace{0.3cm} \textsc{Intersection numbers.} Let
$f_1,\ldots,f_k$ be continuous sections of complex line bundles
$I_1,\ldots,I_k$ on a $k$-dimensional complex algebraic variety
$V$, such that the set $\{f_1=\ldots=f_k=0\}$ is compact.
Consider the Chern class $c_i\in H^2\bigl(V,\, \{f_i\ne 0\};\;
\Z\bigr)$ of the bundle $I_i$, localized near the zero locus of
its section $f_i$. Then the \textit{intersection number} of the
divisors of the sections $f_1,\ldots,f_k$ is defined as
$c_1\smile\ldots\smile c_k\in H^{2k}\bigl(V,\, \cup_i\{f_i\ne
0\};\; \Z\bigr)=\Z$ and is denoted by $m(f_1\cdot\ldots\cdot
f_k\cdot{V})$.

In other words, if we consider smooth (non-holomorphic)
perturbations $\tilde f_1,\ldots,\tilde f_k$ of the sections
$f_1,\ldots,f_k$, such that the system $\tilde f_1=\ldots=\tilde
f_k=0$ has finitely many regular solutions near the set
$\{f_1=\ldots=f_k=0\}$, then each of the solutions can be assigned
a weight $\pm 1$, depending on orientation of the base $d\tilde
f_1,\ldots,d\tilde f_k$ at this point. The intersection number
$m(f_1\cdot\ldots\cdot f_k\cdot{V})$ by definition equals the sum
of these weights.

\vspace{0.3cm} \textsc{Relative Kouchnirenko-Bernstein formula.}
Let $\Delta_1,\ldots,\Delta_m$ be integer polyhedra in
$\mathcal{M}_{\tau^{\vee}}$, and let
$\widetilde\Delta_1,\ldots,\widetilde\Delta_m$ be the Newton
polyhedra of germs of sections $f_1,\ldots,f_m$ of the line
bundles $I_{\Delta_1},\ldots,I_{\Delta_m}$ on the toric variety
$\T^{\Sigma}$. We compute the intersection number of the divisors
of the sections $f_1,\ldots,f_m$ in terms of the polyhedra
$\Delta_1,\ldots,\Delta_m$ and
$\widetilde\Delta_1,\ldots,\widetilde\Delta_m$, provided that the
leading coefficients of $f_1,\ldots,f_m$ are in general position.
\begin{defin} \label{defcompat}
For every face $\Gamma$ of the sum of polyhedra
$\Delta_1,\ldots,\Delta_m$ in $\R^p$, the maximal collection of
faces $\Gamma_1\subset\Delta_1,\ldots,\Gamma_m\subset\Delta_m$,
such that $\Gamma_1+\ldots+\Gamma_m=\Gamma$, is said to be
\textit{compatible}.
\end{defin}
For bounded faces, the word ``maximal" can be omitted in this
definition.
\begin{defin} \label{genposbernst}
The leading coefficients of the sections $f_1,\ldots,f_m$ are
said to be \textit{in general position}, if, for every collection
of bounded compatible faces
$\widetilde\Gamma_1,\ldots,\widetilde\Gamma_m$ of the polyhedra
$\widetilde\Delta_1,\ldots,\widetilde\Delta_m$,
the system of polynomial equations
$f_{1}^{\widetilde\Gamma_{1}}=\ldots=f_{m}^{\widetilde\Gamma_{m}}=0$
has no roots in the maximal torus $\CC^m$.
\end{defin}

\begin{theor}[Relative Kouchnirenko-Bernstein formula, \cite{me05}, \cite{me06}, \cite{mearx}] \label{relbernst}
Let $\Delta_1,\ldots,\Delta_m$ be integer polyhedra in
$\mathcal{M}_{\tau^{\vee}}$, and let
$\widetilde\Delta_1,\ldots,\widetilde\Delta_m$ be the Newton
polyhedra of sections $f_1,\ldots,f_m$ of the line bundles
$I_{\Delta_1},\ldots,I_{\Delta_m}$, such that the difference
$\Delta_i\setminus\widetilde\Delta_i$ is bounded for every $i$.
Then
\newline
1) The intersection number $m(f_1\cdot\ldots\cdot
f_m\cdot\T^{\tau})$ is greater or equal than the mixed volume
$m!\MV\Bigl((\Delta_1,\widetilde\Delta_1),\ldots,(\Delta_m,\widetilde\Delta_m)\Bigr)$.
\newline
2) This inequality turns into an equality if and only if leading
coefficients of the sections $f_1,\ldots,f_m$ are in general
position in the sense of Definition \ref{genposbernst}.
\end{theor}

\subsection{Kouchnirenko-Bernstein-Kho\-van\-skii
formula.}\label{skbk} In the assumptions of Theorem
\ref{relbernst}, suppose that the first line bundle
$\mathcal{I}_{\Delta_1}$ is trivial, i.e. $\Delta_1=\tau^{\vee}$.
The relative version of the Kouchnirenko-Bernstein-Khovanskii
formula computes the Euler characteristic of the Milnor fiber of
the function $f_1$ on the complete intersection
$f_2=\ldots=f_k=0$ for $k\leqslant m$, in terms of the Newton
polyhedra of the sections $f_1,\ldots,f_k$. At the end of this
subsection, we also explain how to drop the assumption on
triviality of $\mathcal{I}_{\Delta_1}$.

To define the Milnor fiber of $s_1$, it is convenient to fix a
family of neighborhoods for the compact part of the toric variety
$\T^{\Sigma}$. For instance, choose an integer point $a_i$ on
every infinite edge of $\Delta_1$, and let $B_{\varepsilon}$ be
the set of all $x\in\CC^m$ such that $\sum_i |x^{a_i}| \leqslant
\varepsilon$. Then its closure in the toric variety $\T^{\Sigma}$
is a neighborhood of the compact part $\T^{\Sigma}_{comp}$.
\begin{defin} The Milnor fiber of the function $f_1$ on the
complete intersection $f_2=\ldots=f_k=0$ is the manifold
$\{f_1-\delta=f_2=\ldots=f_k=0\}\cap B_{\varepsilon}$ for
$|\delta|\ll\varepsilon\ll 1$.
\end{defin}
\begin{defin} \label{genpos2} The leading coefficients of $f_1,\ldots,f_k$, are
said to be in general position, if, for every collection of
bounded compatible faces
$\widetilde\Gamma_1,\ldots,\widetilde\Gamma_k$ of the polyhedra
$\widetilde\Delta_1,\ldots,\widetilde\Delta_k$,
the systems of polynomial equations
$f_{1}^{\widetilde\Gamma_{1}}=\ldots=f_{k}^{\widetilde\Gamma_{k}}=0$
and
$f_{2}^{\widetilde\Gamma_{2}}=\ldots=f_{k}^{\widetilde\Gamma_{k}}=0$
define regular varieties in the maximal torus $\CC^m$.
\end{defin}

We denote the mixed volume of pairs of polyhedra
$(P_1,Q_1),\ldots,(P_m,Q_m)$ in $\R^m$ by the monomial
$(P_1,Q_1)\cdot\ldots\cdot(P_m,Q_m)$.

\begin{theor} \label{relbkh} In the above assumptions,
the Euler characteristic of the Milnor fiber of $f_1$ on the
complete intersection $\{f_2=\ldots=f_k=0\}$ equals
$$(-1)^{m-k}m!\sum_{a_1>0,\ldots,a_k>0\atop a_1+\ldots+a_k=m}
(\Delta_1,\widetilde\Delta_1)^{a_1}\cdot\ldots\cdot
(\Delta_k,\widetilde\Delta_k)^{a_k},$$ provided that the leading
coefficients of $f_1,\ldots,f_k$ are in general position in the
sense of Definition \ref{genpos2}.
\end{theor}

This is proved in \cite{oka} for a regular affine toric variety
(based on the idea of \cite{varch}), and in \cite{tak2} for an
arbitrary affine toric variety (in a more up to date language).
Both arguments can be easily applied to an arbitrary (not
necessary affine) toric variety, and also provide a formula for
the $\zeta$-function of monodromy of the function $f_1$. However,
since we restrict our consideration to the Milnor number in this
paper, we prefer to give a much simpler proof by reduction to the
global Kouchnirenko-Bernstein-Khovanskii formula.

\textsc{Proof.} If the leading coefficients are in general
position, then topology of the Milnor fiber only depends on the
Newton polyhedra of $f_1,\ldots,f_k$, and we can assume without
loss of generality that $f_i=s_{\Delta_i}\cdot\tilde f_i$, where
$\tilde f_1,\ldots,\tilde f_k$ are Laurent polynomials on $\CC^m$
and satisfy the condition of general position of \cite{kh77}.
Denote the Newton polyhedra of the polynomials $\tilde
f_1,\ldots,\tilde f_k$ and $\tilde f_1-\delta$ with $\delta\ne 0$
by $D_1,\ldots,D_k$ and $\widetilde D_1$.

By the global Kouchnirenko-Bernstein-Khovanskii formula
\cite{kh77}, the Euler characteristics of $\{\tilde
f_1=\ldots=\tilde f_k=0\}$ and $\{\tilde f_1-\delta=\tilde
f_2=\ldots=\tilde f_k=0\}$ equal
$(-1)^{m-k}m!\sum\limits_{a_1>0,\ldots,a_k>0\atop
a_1+\ldots+a_k=m} D_1^{a_1}\cdot\ldots\cdot D_k^{a_k}$ and
$(-1)^{m-k}m!\sum\limits_{a_1>0,\ldots,a_k>0\atop
a_1+\ldots+a_k=m} \widetilde D_1^{a_1}\cdot
D_2^{a_2}\cdot\ldots\cdot D_k^{a_k}$ respectively.

Since the boundary of $B_{\varepsilon}$ subdivides the set
$\{\tilde f_1-\delta=\tilde f_2=\ldots=\tilde f_k=0\}$ into two
parts, homeomorphic to the set $\{\tilde f_1=\ldots=\tilde
f_k=0\}$ and the Milnor fiber of $f_1$ on $\{f_2=\ldots=f_n=0\}$,
the Euler characteristic of the latter equals
$$(-1)^{m-k}m!\sum_{a_1>0,\ldots,a_k>0\atop a_1+\ldots+a_k=m}
\widetilde D_1^{a_1}\cdot D_2^{a_2}\cdot\ldots\cdot D_k^{a_k}\,
-\, (-1)^{m-k}m!\sum_{a_1>0,\ldots,a_k>0\atop a_1+\ldots+a_k=m}
D_1^{a_1}\cdot\ldots\cdot D_k^{a_k}$$ by additivity of Euler
characteristic. This difference is equal to
$$(-1)^{m-k}m!\sum_{a_1>0,\ldots,a_k>0\atop a_1+\ldots+a_k=m}
(\Delta_1,\widetilde\Delta_1)^{a_1}\cdot\ldots\cdot
(\Delta_k,\widetilde\Delta_k)^{a_k}$$ by Lemmas \ref{mvmv} and
\ref{mvcone}. $\Box$

Summing up the Euler characteristics of the Milnor fibers of
$f_1|_{\{f_2=\ldots=f_k=0\}}$ over all non-compact toric
subvarieties in $\T^{\Sigma}$, we have the following formula for
the Euler characteristic of the closure of Milnor fiber.

The fact that $\mathcal{I}_{\Delta_i}$ is a line bundle on the
toric variety $\T^{\Sigma}$ implies that, for every cone
$\sigma\in\Sigma$, all interior points of $\sigma$, being
considered as linear functions on the polyhedron $\Delta_i$,
attain their minimum on the same face of $\Delta_i$. Denote this
face by $A_i$, the codimension of the cone $\sigma$ by $q$, and
pick a vector $a_i\in A_i$. Then the shifted pairs
$\mathcal{A}_i=\Bigl(A_i-a_i,(A_i\cap\widetilde\Delta_i)-a_i\Bigr)$
are contained in the same rational $q$-dimensional subspace of
$\R^m$, and their $q$-dimensional mixed volumes make sense. Denote
the number $(-1)^{q-k}q!\sum\limits_{a_1>0,\ldots,a_k>0\atop
a_1+\ldots+a_k=q} \mathcal{A}_1^{a_1}\cdot\ldots\cdot
\mathcal{A}_k^{a_k}$ by $\chi_{\sigma}$.

\begin{sledst} \label{relbkhcomp} In the above assumptions, the Euler
characteristic of the closure of the Milnor fiber of $f_1$
on the complete intersection $\{f_2=\ldots=f_n=0\}$ equals
$\sum_{\sigma\in\Sigma}\chi_{\sigma}$, provided that the leading
coefficients of $s_1,\ldots,s_k$ are in general position.
\end{sledst}

We also formulate a more general version of this theorem, with no
assumptions on the triviality of the first line bundle (we do not
need this generalization in what follows; since the proof is
similar to that of Corollary \ref{relbkhcomp}, we omit it). Let
$f_1,\ldots,f_k$ be germs of holomorphic sections of arbitrary
line bundles
$\mathcal{I}_{\Delta_1},\ldots,\mathcal{I}_{\Delta_k}$ near the
compact part of a toric variety $\T^{\Sigma}$, and pick
holomorphic sections $t_1,\ldots,t_k$ of these bundles in the
closure of the set $B_{\varepsilon}$ for a small $\varepsilon$.
Varieties $B_{\varepsilon}\cap\{f_1-t_1=\ldots=f_k-t_k=0\}$ are
diffeomorphic to each other for almost all collections of small
sections $(t_1,\ldots,t_k)$. Such variety is called the Milnor
fiber of the complete intersection $\{f_1=\ldots=f_k=0\}$.

\begin{theor} In the above assumptions, the Euler characteristic of the Milnor fiber of
the complete intersection $\{f_1=\ldots=f_k=0\}$ equals
$$(-1)^{m-k}m!\sum_{a_1>0,\ldots,a_k>0\atop a_1+\ldots+a_k=m}
(\Delta_1,\widetilde\Delta_1)^{a_1}\cdot\ldots\cdot
(\Delta_k,\widetilde\Delta_k)^{a_k},$$ provided that the leading
coefficients of $f_1,\ldots,f_k$ are in general position.
\end{theor}

\subsection{Euler obstructions of polyhedra.} \label{seulobstr}

\textsc{Euler obstructions of varieties.} Let
$\bigsqcup_{\alpha\in\Lambda} U_{\alpha}$ be a Whitney
stratification of a complex algebraic variety $U$. Pick a point
$x_0$ in a stratum $U_{\alpha'}$ and consider a germ of an
analytic function $f:(U,U_{\alpha'})\to(\C,0)$ at this point. The
Euler characteristic of the set $\{x\in U_{\alpha}\;\; |\;\;
f(x)=\delta,\, |x-x_0|\leqslant\varepsilon\}$ takes the same
value for almost all germs $f$ and all positive numbers
$\delta\ll\varepsilon\ll 1$. This value does not depend on the
choice of $x_0\in U_{\alpha'}$, and we denote the negative of this
value by $\mu^{\alpha',\, \alpha}$ (note that it equals 0 unless
$U_{\alpha}$ is adjacent to $U_{\alpha'}$). We also set
$\mu^{\alpha,\, \alpha}=1$ for every $\alpha\in\Lambda$ and denote
the $|\Lambda|\times|\Lambda|$ matrix with entries
$\mu^{\alpha',\, \alpha}$ by $M$.

\begin{defin}
The $(\alpha',\,\alpha)$-entry of the inverse matrix $M^{-1}$ is
denoted by $\epsilon^{\alpha',\, \alpha}$ and is called the
\textit{Euler obstruction} of the closure of the stratum
$U_{\alpha}$ at a point of the stratum $U_{\alpha'}$.
\end{defin}
\begin{rem}
Adjacency of strata induces a partial order structure on the set
$\Lambda$. If we consider the collections $\mu^{\alpha',\,
\alpha}$ and $\epsilon^{\alpha',\, \alpha}$ as functions on the
poset $\Lambda\times\Lambda$, these functions are called
\textit{M\"obius inverse} (see e.g. \cite{hall} for details).
\end{rem}
\begin{exa} Let $(U,0)$ be an $m$-dimensional isolated toric singularity (it consists of
two strata $U\setminus\{0\}$ and $\{0\}$). Denote the
$m$-dimensional cone of its fan by $\tau$, the convex hull of
$\tau^{\vee}\cap\Z^n\setminus\{0\}$ by $H$, and the volume of
$\tau^{\vee}\setminus H$ by $V$. Then, by Corollary
\ref{relbkhcomp}, the above matrix $M$ equals $\begin{pmatrix}
  1 & m!V-2 \\
  0 & 1
\end{pmatrix}$ for even $m$ and $\begin{pmatrix}
  1 & -m!V \\
  0 & 1
\end{pmatrix}$ for odd $m$. Thus, the Euler
obstruction of $U$ equals $2-m!V$ for even $m$ and $m!V$ for odd
$m$.
\end{exa}

\textsc{Combinatorial coefficients $c^{A',A}$ and $e^{A',A}$.}
\begin{defin} A subset $A'\subset A\subset\Z^k$ is called a \textit{face}
of $A$, if it can be represented as the intersection of $A$ with
a face of the convex hull of $A$.\end{defin} For a face $A'$ of a
finite set $A\subset\R^k$, let $M'$ and $M\subset\R^k$ be the
vector spaces, parallel to the affine spans of $A'$ and $A$
respectively. We denote the projection $\R^k\to\R^k/M'$ by $s$,
and choose the volume form $\eta$ on $M/M'$ such that the volume
of $M/(M'+\Z^k)$ equals $(\dim M/M')!$. Then the $\eta$-volume of
the difference of the convex hulls of $s(A)$ and $s(A\setminus
A')$ is denoted by $c^{A',\, A}\in\Z$. Set $c^{A,\, A}$ to 1 and
$c^{A',\, A}$ to 0 if $A'$ is not a face of $A$.
\begin{exa}
The above difference is shown by hatching below for the following
two examples:
\newline 1) $A_1=\{(0,0,0),\; (0,0,1),\; (2,1,0),\; (1,2,1)\}$
and $A'_1$ is its vertical edge. \newline 2) $A_2$ is the set of
integer lattice points in the convex hull of $2\cdot A_1$, and
$A'_2$ is its vertical edge.
\end{exa}
\begin{center}
\noindent\includegraphics[width=9cm]{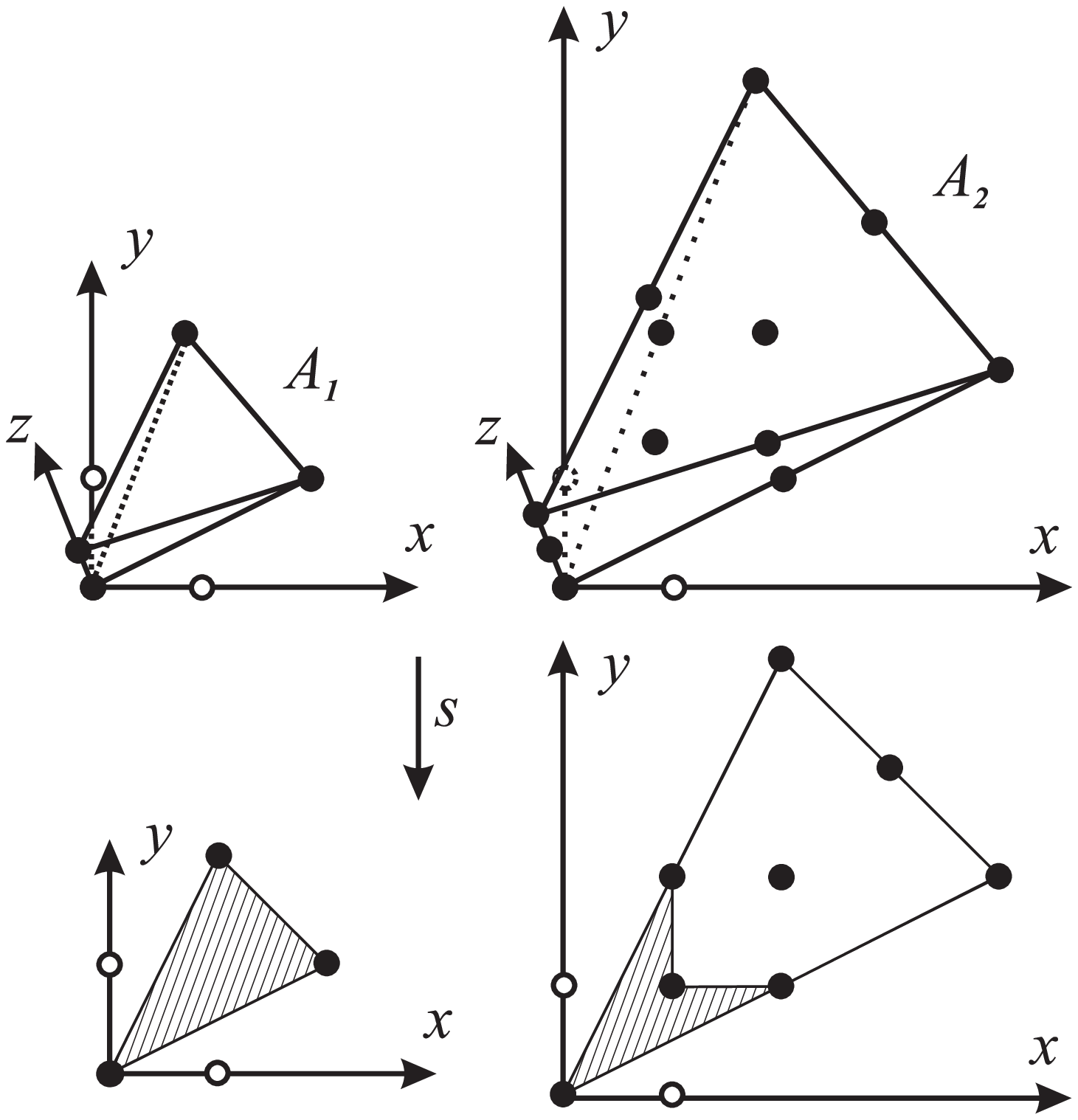} \end{center}

Consider the square matrix $C$ with entries $c^{A'',A'}$, where
$A''$ and $A'$ run over the set of all faces of $A$, and define
$e^{A'',\, A'}$ as the $(A'',\, A')$-entry of the inverse of $C$.
Note that $C$ is upper triangular with $1$'s on the diagonal, if
we order faces of $A$ by their dimension; in particular, the
determinant of $C$ equals 1, and its inverse is integer (although
positivity of the entries of $C$ does not imply positivity of the
entries of its inverse).
\begin{defin} The number $e^{A',\, A}$ is called the \textit{Euler
obstruction of the set} $A$ at its face $A'$.
\end{defin}
For instance, restricting our attention to the four faces of the
set $A_2$ adjacent to the vertical edge $A'_2$ in the previous
example (including $A'_2$ itself), we obtain the Euler obstruction
$e^{A_2,A'_2}$ as the top right element in the matrix $C^{-1}$:

$$C=\begin{pmatrix}
   1 &  1 &  1 &  2 \\
     &  1 &  0 &  1 \\
     &    &  1 &  1 \\
     &    &    &  1 \\
\end{pmatrix},\quad C^{-1}=\begin{pmatrix}
   1 &  -1 &  -1 &  0 \\
     &  1 &  0 &  -1 \\
     &    &  1 &  -1 \\
     &    &    &  1 \\
\end{pmatrix}.$$

For an integer polyhedron $P\subset\R^m$ and its face $Q$, we
call the numbers $c^{P\cap\Z^m,\, Q\cap\Z^m}$ and $e^{P\cap\Z^m,\,
Q\cap\Z^m}$ \textit{the Milnor number and the Euler obstruction
of the polyhedron} $P$ \textit{at its face} $Q$, and denote them
by $c^{Q,\, P}$ and $e^{Q,\, P}$.
\begin{exa} \label{exachain} For an integer polygon
$P$ at its vertex $Q$, denoting the convex hull of $\Z^2\cap
P\setminus \{Q\}$ by $P_Q$, we have $$c^{Q,\, P}=2\cdot\mbox{area
of } P\setminus P_Q,\quad e^{Q,\, P}=2-2\cdot\mbox{area of }
P\setminus P_Q.$$ In particular, if $c^{Q,\, P}=e^{Q,\, P}=1$,
then the adjacent edges of the vertex $Q$ can be brought to the
lines
$$y=0 \mbox{ and } x=0$$ by a suitable affine change of
coordinates in $\R^2$ that preserves the integer lattice. If
$c^{Q,\, P}=2$ (and $e^{Q,\, P}=0$), then the adjacent edges of
the vertex $Q$ can be brought to the lines of the form
$$y=0 \mbox{ and } y=\left(1+\frac{1}{n}\right)x, \mbox{ where } n\in\N.$$

More generally, if $c^{Q,\, P}=s$, then  the adjacent edges of the
vertex $Q$ can be brought to the lines of the form
$$y=0 \mbox{ and } y=\left(p_0+\frac{1}{p_1+\cdots\frac{1}{p_{2k}}}\right)x, \mbox{ where } p_{i}\in\N \mbox{ and } \sum_i p_{2i}=s$$
(the number in the brackets is the continued fraction of the
sequence $p_0,\ldots,p_{2k}$). It would be interesting to extend
this to classification of $r$-dimensional rational convex
polyhedral cones $P$ such that $c^{0,\, P}=s$ or $e^{0,\, P}=s$
for small $r$ and $|s|$; this problem is related to problems of
classification of multidimensional continued fractions and sails
(see e.g. \cite{karp}).
\end{exa}

\vspace{0.3cm} \textsc{Geometric meaning of  $c^{A',A}$ and
$e^{A',A}$.} For a set $A=\{a_1,\ldots,a_N\}\subset\Z^k$, such
that the differences $a_i-a_j$ generate the lattice $\Z^k$, the
closure of the image of the torus $\CC^k$ under the inclusion
$j:\CC^k\to\CP^A,\; j(t)=(t^{a_1},\ldots,t^{a_N})$, is a toric
variety, whose orbits are in one to one correspondence with faces
of $A$. Its subdivision into the orbits is a Whitney
stratification $\bigsqcup U_{A'}$, where $A'$ runs over all faces
of $A$. This stratification gives rise to the numbers
$\mu^{A'',\, A'}$ and $\epsilon^{A'',\, A'}$ for every pair of
faces $A'$ and $A''$, as defined in the beginning of this
subsection.

Since, by Theorem \ref{relbkh}, we have $\mu^{A''\,
,A'}=(-1)^{\dim A''-\dim A'}c^{A'',\, A'}$, then
$\epsilon^{A'',\, A'}=(-1)^{\dim
A''-\dim A'}e^{A'',\, A'}$, 
which proves

\begin{utver}[\cite{tak}]\label{thmt} The Euler obstruction of the set $A\subset\Z^k$ at its
face $A'$ equals $(-1)^{\dim A'-\dim A}$ times the Euler
obstruction of the toric variety, corresponding to $A$, at a point
of its orbit, corresponding to $A'$.
\end{utver}

\subsection{Multiplicities of non-degenerate complete intersections.}\label{smultci}

Here we formulate a corollary of the relative
Kouchnerenko-Bernstein formula that leads to a simple proof the
Gelfand-Kapranov-Zelevinsky decomposition formula (we do not need
it for other purposes).

\vspace{0.3cm} \textsc{Varieties with multiplicities.} Let
$f_1,\ldots,f_l$ be germs of holomorphic sections of complex line
bundles $I_1,\ldots,I_l$ on a germ of a $k$-dimensional complex
algebraic variety $(V,x)$, such that the set
$\{f_1=\ldots=f_l=0\}$ is smooth and $(k-l)$-dimensional. Then we
can choose germs of holomorphic functions $f_{l+1},\ldots,f_k$ on
$(V,x)$, such that the differentials of their restrictions to
$\{f_1=\ldots=f_l=0\}$ are linearly independent. The local
topological degree of the map $(f_1,\ldots,f_k):(V,x)\to(\C^k,0)$
does not depend on the choice of the germs $f_{l+1},\ldots,f_k$
and is called \textit{the multiplicity of the (local) complete
intersection} $f_1=\ldots=f_l=0$ at its point $x$.

If $S_i,\; i=1,\ldots,I,$ are the irreducible components of a
complete intersection $\{f_1=\ldots=f_l=0\}$, and $a_i$ is the
multiplicity of this complete intersection at a smooth point of
$S_i$, then we denote the cycle $\sum_{i=0}^I a_iS_i$ (i.e. a
formal sum of irreducible varieties) by $[f_1=\ldots=f_l=0]$.

If $S=\sum a_iS_i$ is a cycle in $V$ (every $S_i$ is an
irreducible variety in $V$, and every $a_i$ is a positive
number), and $f:V\to W$ is a proper map, then we define the
\textit{image} $f_*(S)$ as follows. For every component $S_i$,
denote the topological degree of the map $f:S_i\to f(S_i)$ by
$d_i$, provided that $\dim f(S_i)=\dim S_i$ (otherwise, set
$d_i=0$ by definition). Then the \textit{image} $f_*(S)$ is
defined to be the sum $\sum a_id_if(S_i)$.

\vspace{0.3cm} \textsc{Multiplicities of nondegenerate complete
intersections.} Let $\Delta\subset\R^m$ be an $m$-dimensional
integer polyhedron, let $\Sigma$ be its dual fan, and let
$\widetilde\Delta$ be the Newton polyhedron of holomorphic
sections $f_1,\ldots,f_l$ of the line bundle $I_{\Delta}$ on the
toric variety $\T^{\Sigma}$. Denote the orbit of $\T^{\Sigma}$,
corresponding to a codimension $l$ face $\Gamma$ of $\Delta$, by
$\T^{\Gamma}$. Then the multiplicity of the complete intersection
$f_1=\ldots=f_l=0$ at a generic point of the orbit $\T^{\Sigma}$
can be computed as follows.

Let $p:\R^m\to\R^l$ be a projection such that $p(\Z^m)=\Z^l$ and
$p(\Gamma)$ is a point.
\begin{lemma} \label{orbmult}
If $p(\widetilde\Delta)$ touches all faces of $p(\Delta)$ except
the vertex $p(\Gamma)$, and the leading coefficients of
$f_1,\ldots,f_l$ are in general position, then the closure of the
orbit $\T^{\Gamma}$ is a component of multiplicity
$l!\Vol\Bigl(p(\Delta)\setminus p(\widetilde\Delta)\Bigr)$ in the
complete intersection $f_1=\ldots=f_l=0$.
\end{lemma}
The assumption that the Newton polyhedra of the sections
$f_1,\ldots,f_l$ coincide is obviously redundant; it is
introduced to simplify the notation.

\textsc{Proof.} Let $\T_{\Gamma}\subset\CC^m$ be the
$l$-dimensional subtorus that acts trivially on the orbit
$\T^{\Gamma}$. For every point $x\subset\T^{\Gamma}$, consider
the $\T_{\Gamma}$-invariant $l$-dimensional closed toric
subvariety $H_x\subset\T^{\Sigma}$ that intersects the orbit
$\T^{\Gamma}$ at the point $x$.

We should prove that, for a generic point $x$, the intersection
number $\mathcal{J}$ of the variety $H_x$ and the complete
intersection $f_1=\ldots=f_l=0$ at $x$ makes snese and equals
$l!\Vol\Bigl(p(\Delta)\setminus p(\widetilde\Delta)\Bigr)$. To do
so, we denote the restrictions of the line bundle $I_{\Delta}$
and its sections $f_1,\ldots,f_l$ to $H_x$ by $I'$ and
$f'_1,\ldots,f'_l$ respectively and note that $I'=I_{p(\Delta)}$
and the Newton polyhedra of $f'_1,\ldots,f'_l$ are equal to
$p(\widetilde\Delta)$. On one hand, the intersection number of
the divisors of $f'_1,\ldots,f'_l$ is equal to the desired
intersection number $\mathcal{J}$, while on the other hand it
makes sense and equals $l!\Vol\Bigl(p(\Delta)\setminus
p(\widetilde\Delta)\Bigr)$ by Theorem \ref{relbernst} for the
sections $f'_1,\ldots,f'_l$. $\quad\Box$

\section{$A$-discriminants.}

In this section, we discuss the universal case of our problem.
For a collection of finite sets $A_0,\ldots,A_k$ in $\Z^k$, we
recall the definition of the $(A_0,\ldots,A_k)$-resultant
(Subsection \ref{sadeterm}) and the $A_0$-discriminant
(Subsection \ref{sadiscrim}). More generally, for every
$l\leqslant k$, we introduce the so called
$(A_0,\ldots,A_l)$-discriminant (Subsection \ref{sgkzdiscrimci}),
and express it in terms of $A$-discriminants by means of the
Cayley trick (Subsection \ref{scayley}).

The collection $A_0,\ldots,A_l$ is called dual defective, if the
$(A_0,\ldots,A_l)$-dis\-cri\-mi\-nant set is not a hypersurface,
and we give a number of examples of sufficient conditions for
non-dual defectiveness of a collection (Subsection \ref{sdualdef}
for $l=0$ and Proposition \ref{equidim} for arbitrary $l$). We
also consider an alternative version of the
$(A_0,\ldots,A_l)$-discriminant set (the bifurcation set, see
Definition \ref{defbifuniv}), which is presumably always a
hypersurface (see Proposition \ref{conj1sim} and the subsequent
conjecture). The technical proof of Propositions \ref{equidim} and
\ref{conj1sim} and the alternative prove of the
Gelfand-Kapranov-Zelevinsky decomposition formula are postponed
till the end of this section.

\subsection{Resultants and $A$-determinants (\cite{sturmf} and \cite{gkz}).} \label{sadeterm}
\textsc{Resultant.} For a finite set $A\subset\Z^k$, denote the
set of all Laurent polynomials of the form $\sum_{a\in A} c_ax^a$
on the complex torus $\CC^k$ by $\C[A]$. Consider finite sets
$A_0,\ldots,A_l$ in $\Z^k$, such that the dimension of the convex
hull of $A_0+\ldots+A_l$ is not greater than $l$. Let $\Sigma$ be
the set of all collections
$(\varphi_0,\ldots,\varphi_l)\in\C[A_0]\oplus\ldots\oplus\C[A_l]$
such that $\varphi_0(y)=\ldots=\varphi_l(y)=0$ for some
$y\in\CC^k$.

If the closure of $\Sigma$ is a hypersurface in
$\C[A_0]\oplus\ldots\oplus\C[A_l]$, then

1) it is defined by the equation $G=0$ for a certain irreducible
polynomial $G$ of positive degree.

2) For a generic collection
$(\varphi_0,\ldots,\varphi_l)\in\Sigma$, the set $\{y\; |\;
\varphi_0(y)=\ldots=\varphi_l(y)=0\}$ can be represented as
$J\cdot T$, where $$T=\{z\, |\, z^a=z^b \mbox{ for all } a \mbox{
and } b \mbox{ in the afine span of } A_0+\ldots+A_l\}$$ is a
subtorus in $\CC^k$, and $J\subset\CC^k/T$ is a certain finite
set, whose cardinality does not depend on the choice of
$(\varphi_0,\ldots,\varphi_l)$; we denote this cardinality by
$d(A_0,\ldots,A_l)$.

\begin{defin} \label{defgkzres} If the closure of $\Sigma$ is a hypersurface,
then the polynomial $G^{d(A_0,\ldots,A_l)}$ is called the
\textit{$(A_0,\ldots,A_l)$-resultant} and is denoted by
$R_{A_0,\ldots,A_l}$, otherwise we set $R_{A_0,\ldots,A_l}=1$ by
definition.
\end{defin}

\begin{rem} The resultant
is uniquely defined up to multiplication by a non-zero constant,
and all equalities involving resultants should be understood
correspondingly.

This definition differs from \cite{sturmf} and \cite{gkz} by the
exponent $d(A_0,\ldots,A_l)$. This exponent slightly simplifies
computations and can be easily expressed in terms of $A$ as
follows.
\end{rem}

\vspace{0.3cm} \textsc{Multiplicity of the resultant.} Let
$A_0,\ldots,A_k$ be finite sets in $\Z^k$. The resultant
$R_{A_0,\ldots,A_k}$ is by definition a certain power
$d(A_0,\ldots,A_k)$ of an irreducible polynomial. Here we recall
an explicit formula for the number $d(A_0,\ldots,A_k)$ and a
criterion for triviality of the resultant $R_{A_0,\ldots,A_k}$ in
terms of the sets $A_0,\ldots,A_k$.
\begin{defin} The dimension of the convex hull of a finite set $A\subset\R^k$
is called the dimension of $A$ and is denoted by $\dim A$.
\end{defin}
\begin{defin} For every non-empty subset $J\subset\{0,\ldots,k\}$, the
difference\newline $\dim\sum_{j\in J}A_j - |J|$ is called the
\textit{codimension of the collection} $A_j,\, j\in J$, and is
denoted by $\codim J$.
\end{defin}
\begin{utver}[\cite{sturmf}]\label{essent} $\,$\newline 1) There exists $J$ with $\codim J<-1$ if and
only if $R_{A_0,\ldots,A_k}=1$. \newline 2) If $\codim
J\geqslant-1$ for every $J$, then there exists the minimal set
$J_0\subset\{0,\ldots,k\}$, such that $\codim J_0=-1$.
\end{utver}
Under the assumption of Proposition \ref{essent}(2), let
$L_{\Z}\subset\Z^{k+1}$ be the lattice generated by the set
$(\sum_{j\in J_0}A_j)\times\{1\}\subset\Z^k\oplus\Z$, and let
$L\subset\R^{k+1}$ be the linear span of $L_{\Z}$, denote the
number $\bigl|(L\cap\Z^{k+1})/L_{\Z}\bigr|$ by $d_1$.

Denote the projection $\R^{k+1}\to\R^{k+1}/L$ by $s$, and choose
the volume form $\eta$ on $\R^{k+1}/L$ such that the volume of
$\R^{k+1}/(L+\Z^{k+1})$ equals $(k-|J_0|+1)!$. Denote the
$\eta$-mixed volume of the convex hulls of the sets
$s(A_j\times\{1\}),\; j\in\{0,\ldots,k\}\setminus J_0$, by $d_2$.
\begin{utver}[\cite{detandnewt}, \cite{mearx}]\label{lmultres}
In the notation above, $d(A_0,\ldots,A_k)=d_1\cdot d_2$.
\end{utver}

\vspace{0.3cm} \textsc{$A$-determinant.} Let $l$ be the dimension
of the affine span $L$ of a finite set $A\subset\Z^k$, let
$t_1,\ldots,t_k$ be the standard coordinates on $\R^k$, and let
$y_1,\ldots,y_k$ be the standard coordinates on $\CC^k$. Choose
numbers $i_1,\ldots,i_l$ such that the functions
$t_{i_1},\ldots,t_{i_l}$ form a system of coordinates on $L$ (for
example, if $\dim A=k$, then $\{i_1,\ldots,i_l\}=\{1,\ldots,k\}$).

\begin{defin}\label{defgkzdeterm}
The $A$-\textit{determinant} is the polynomial $E_A$ on $\C[A]$,
defined by the equality
$$E_A(\varphi)=R_{A,\ldots,A}\Bigl(\varphi,\;y_{i_1}\frac{\partial\varphi}{\partial
y_{i_1}},\;\ldots,\;y_{i_l}\frac{\partial\varphi}{\partial
y_{i_l}}\Bigr).$$
\end{defin}
\begin{rem} The $A$-determinant
is uniquely defined up to multiplication by a non-zero constant,
and does not depend on the choice of the collection
$i_1,\ldots,i_l$, because changing this collection results in
multiplication of the vector
$\Bigl(\varphi,\;y_{i_1}\frac{\partial\varphi}{\partial
y_{i_1}},\;\ldots,\;y_{i_l}\frac{\partial\varphi}{\partial
y_{i_l}}\Bigr)$ by a square non-degenerate matrix (by definition,
the resultant $R_{A,\ldots,A}$ is invariant with respect to
multiplication of its argument by a non-degenerate matrix).
\end{rem}

\subsection{$A$-discriminants (\cite{gkz}).} \label{sadiscrim}

For a finite set $A\in\Z^k$, let $\Sigma_A\subset\C[A]$ be the
set of all polynomials $\varphi\in\C[A]$ such that both $\varphi$
and its differential ${\rm d}\varphi$ vanish at some point
$y\in\CC^k$.

If the closure of $\Sigma_A$ is a hypersurface, then

1) it is defined by the equation $G=0$ for a certain irreducible
polynomial $G$ of positive degree.

2) For a generic $\varphi\in\Sigma_A$, its singular locus $\{y\;
|\; \varphi(y)={\rm d}\varphi(y)=0\}$ has the form $J\cdot T$,
where $T=\{z\, |\, z^a=z^b \mbox{ for all } a \mbox{ and } b
\mbox{ in the affine span of } A\}$ is a subtorus in $\CC^k$, and
$J\subset\CC^k/T$ is a finite set, whose cardinality $|J|$ does
not depend on $\varphi$.

\begin{defin} \label{defgkzdiscr} If the closure of $\Sigma_A$ is a hypersurface, then
the polynomial $G^{|J|}$ is called the \textit{$A$-discriminant}
and is denoted by $D_{A}$; otherwise we set $D_A=1$.
\end{defin}

\begin{rem} The discriminant
is uniquely defined up to multiplication by a non-zero constant,
and all equalities involving discriminants should be understood
correspondingly.

This definition differs from \cite{gkz} by the exponent $|J|$.
This exponent slightly simplifies computations in what follows,
and can be easily expressed in terms of $A$: let
$L_{\Z}\subset\Z^{k+1}$ be the lattice generated by the set
$A\times\{1\}\subset\Z^k\oplus\Z$, and let $L\subset\R^{k+1}$ be
its linear span, then $|J|=\bigl|(L\cap\Z^{k+1})/L_{\Z}\bigr|$.
\end{rem}
\begin{lemma} \label{lirred} For every $A\subset\Z^m$, the discriminant $D_A$ is
a power of an irreducible polynomial.
\end{lemma}
\textsc{Proof.} The set of all pairs
$(\varphi,y)\in\C[A]\times\CC^m$, such that $\varphi(y)={\rm
d}\varphi(y)=0$, is the total space of a vector bundle
$\mathcal{V}$ with the base space $\CC^m$ and the projection
$\C[A]\times\CC^m\to\CC^m$. Thus it is irreducible, thus its
image $\Sigma_A$ under the projection to $\C[A]$ is also
irreducible. $\quad\Box$

Recall that the dimension of the convex hull of $A\subset\Z^k$ is
denoted by $\dim A$.
\begin{lemma} \label{ldimsing} For a generic $\varphi\in\C[A]$,
the codimension of the set of all points $y\in\CC^k$, such that
$\varphi(y)={\rm d}\varphi(y)=0$, equals $\dim
A-\codim\Sigma_A+1$.
\end{lemma}
\textsc{Proof.} The fiber of the vector bundle $\mathcal{V}$,
introduced in the previous proof, has codimension $1+\dim A$ in
$\C[A]$. $\quad\Box$

$A$-discriminants and $A$-determinants are related as follows.
For a polynomial $\varphi(y)=\sum_{a\in A}\varphi_ay^a$ and a
face $A'$ of the set $A$, denote the polynomial $\sum_{a\in
A'}\varphi_ay^a$ by $\varphi^{A'}$. Recall that coefficients
$c^{A',A}$ are introduced in Subsection \ref{seulobstr}.
\begin{utver}[\cite{gkz}]\label{gkzdecomp} $$E_A(\varphi)= \prod_{A'}
\Bigl(D_{A'}(\varphi^{A'})\Bigr)^{c^{A',\, A}},$$ where $A'$ runs
over all faces of $A$, including $A'=A$.
\end{utver}
This also follows from Theorem \ref{relbernst}, see Subsection
\ref{sproofdecomp}.

\subsection{Dual defectiveness.}\label{sdualdef}
\begin{defin} \label{defnondeg} A finite set $A\subset\Z^k$ is said to be \textit{dual defect}, if $D_A=1$.
\end{defin}
We recall a few simple facts about dual defect sets. There is
also a well known way to decide combinatorially if the set is dual
defect or not (Corollary \ref{ldegree}), which will follow from
our results on Newton polyhedra of discriminants. One more
prospective combinatorial criterion for dual defectiveness is
given by Conjecture \ref{conjdualdef2} below. Note that these
facts do not provide classification of dual defect sets, which is
a much more complicated problem, and is solved only for Delzant
polytopes by now, see \cite{roc}. An obvious but useful
reformulation of this definition is as follows.
\begin{lemma} \label{dualdef2} A finite set $A\subset\Z^k$, whose convex hull is
$k$-dimensional, is dual defect if and only if a generic
polynomial in $\C[A]$ has a singular point.
\end{lemma}
\textsc{Proof.} If a generic polynomial $\varphi\in\C[A]$ has a
singular point, then a generic line of the form $\{\varphi-c \,
|\, c\in\C$ intersects the set $\Sigma_A$, thus
$\codim\Sigma_A=1$.

If $\codim\Sigma_A=1$ and the convex hull of $A$ is
$k$-dimensional, then a generic polynomial in $\Sigma_A$ has an
isolated singular point, thus all nearby polynomials in
$\Sigma_A$ has a singular point as well. $\quad\Box$

\begin{utver}[Monotonicity; see also \cite{cds}] \label{lmonot} If a subset $A'$ of a finite set $A\subset\Z^k$ is not dual defect and
is not contained in an affine hyperplane, then $A$ is not dual
defect.\end{utver} \textsc{Proof.} We prove this by induction on
$|A|$. Since dual defectiveness is by definition invariant with
respect to parallel tranlations, the inductive step can be reduced
to the following fact: if $A'$ is not dual defect and is not
contained in an affine hyperplane, then $A'\cup\{0\}$ is not dual
defect. The latter implicatio follows by Lemma \ref{dualdef2}
$\quad\Box$

Dual defect sets are ``thin" in many senses, for example:
\begin{utver} \label{pthin} If $A$ is not contained in a
union of two parallel hyperplanes, then it is not dual defect.
\end{utver}
This can be deduced from \cite{cat}; we also give a simple
self-contained proof, which consists of two simple properties of
iterated circuits.

\begin{defin}
A set $B\subset\R^m$ of cardinality $m+2$ is called a
\textit{circuit}, if none of its cardinality $m+1$ subsets is
contained in an affine hyperplane.

A set $B\subset\R^m$ is called an \textit{iterated circuit}, if
it is not contained in an affine hyperplane, and if, after a
suitable parallel translation, it can be represented as a disjoint
union $\{0\}\sqcup B_1\sqcup\ldots\sqcup B_p$, such that the
following condition is satisfied: denote the linear span of
$\{0\}\sqcup B_1\sqcup\ldots\sqcup B_i$ by $L_i$ for
$i=0,\ldots,p$, then the projection $L_{i+1}\to L_{i+1}/L_i$ maps
the union $\{0\}\sqcup B_{i+1}$ injectively onto a circuit in
$L_{i+1}/L_i$.

The minimal possible number $p$ in this representation is called
the \textit{depth} of the iterated circuit $B$.
\end{defin}
\begin{exa} Some iterated circuits do not contain circuits. The
simplest example is the set of integer points in the unit ball
centered at the origin in $\R^m$ (this set is an iterated
circuit: the sets $B_i$ above are the pairs of its opposite
points). If the cardinality of $B\subset\R^m$ is $2m$ or greater,
then it is not an iterated circuit.
\end{exa}
\begin{lemma}\label{ic1} An iterated circuit $B\subset\Z^m$ is not dual defect.
\end{lemma}
\textsc{Proof.} Proving this by induction on the depth of $B$,
the inductive step can be reduced to the following statement by a
suitable $\Q$-affine change of coordinates in $\R^m$: if
$B\subset\Z^k$ is not dual defect and $B'\sqcup\{0\}\subset\Z^l$
is a circuit, then $B\sqcup B'\subset\Z^k\oplus\Z^l$ is not dual
defect.

Since $B'\sqcup\{0\}$ is a circuit and $B$ is not dual defect,
then generic polynomials $\varphi\in\C[B]$ and $\psi\in\C[B']$
have singular points by Lemma \ref{dualdef2}, thus their sum
$\varphi+\psi$, which is a generic polynomial in $\C[B\sqcup
B']$, has a singular point as well, hence $B\sqcup B'$ is not
dual defect by the same lemma. $\quad\Box$

\begin{lemma}\label{ic2} If a finite set $B\subset\R^m$ does not contain an
iterated circuit, then it is contained in a union of two parallel
hyperplanes.
\end{lemma}
\begin{exa} If $B$ does not contain a circuit, it does not imply that $B$ is contained in a union of two parallel
hyperplanes. The simplest example is the same as the previous one.
\end{exa}
\textsc{Proof.} Without loss of generality we may assume that
$0\in B$, and pick a maximal subset $B'\subset B$, such that

1) $0\in B'$, and

2) $B'$ is an iterated circuit as a subset of its vector span $L$.

The image $\widetilde B$ of the set $B$ under the projection
along $L$ consists of at most $k-\dim L + 1$ points, otherwise
$B$ were not maximal. Thus $\widetilde B$ is contained in the
union of two parallel hyperplanes, and so does $B$. $\quad\Box$

\vspace{0.3cm} \textsc{Proof of Proposition \ref{pthin}.} If $A$
is dual defect, then it does not contain an iterated circuit by
Proposition \ref{lmonot} and Lemma \ref{ic1}, thus it is
contained in the union of two parallel hyperplanes by Lemma
\ref{ic2}. $\quad\Box$

In particular, we have proved the easy half of the following
conjecture.
\begin{conj}\label{conjdualdef2} A finite set $A\subset\Z^k$, whose convex hull is
$k$-dimensional, is dual defect if and only if it does not
contain an iterated circuit.
\end{conj}

Recall that $D_A$ is a polynomial in the indeterminate
coefficients $\varphi_a$ of the polynomial $\varphi(y)=\sum_{a\in
A}\varphi_ay^a$.
\begin{lemma} \label{lnontriv} For every non-dual defect $A\subset\Z^m$ and every $a\in A$,
the discriminant $D_A$ has positive degree as
a polynomial of $\varphi_a$, unless it is a constant.
\end{lemma}
\textsc{Proof.} Since dual defectiveness is by definition
invariant with respect to parallel translations, it is enough to
prove the statement for $a=0$, which follows from Lemma
\ref{dualdef2}. $\quad\Box$

\subsection{Discriminants $D_{A_0,\ldots,A_l}$ and $B_{A_0,\ldots,A_l}$.}\label{sgkzdiscrimci}

\textsc{Discriminant $D_{A_0,\ldots,A_l}$.} Let $A_0,\ldots,A_l,\;
l\leqslant k,$ be finite sets in $\Z^k$, and let
$\Sigma_{A_0,\ldots,A_l}\subset\C[A_0]\oplus\ldots\oplus\C[A_l]$
be the set of all collections of polynomials
$(\varphi_0,\ldots,\varphi_l)$, such that their differentials are
linearly dependent at some point of the set $\{y\in\CC^k\; |\;
\varphi_0(y)=\ldots=\varphi_l(y)=0\}$. The union of the
codimension 1 components of the closure
$\overline{\Sigma_{A_0,\ldots,A_l}}$ is defined by the equation
$G=0$ for a certain square-free polynomial $G$ on
$\C[A_0]\oplus\ldots\oplus\C[A_l]$.

\begin{defin} \label{defgkzdiscrci} The polynomial $G$ is called the \textit{reduced $A$-discriminant
of degree} $(A_0,\ldots,A_l)$ and is denoted by
$D^{red}_{A_0,\ldots,A_l}.$
\end{defin}

We assume that $l\leqslant k$, because the codimension of
$\Sigma_{A_0,\ldots,A_l}$ is greater than 1 otherwise (however,
one can still study the tropicalization of
$\Sigma_{A_0,\ldots,A_l}$ instead of its Newton polyhedron, see
\cite{dfs} and \cite{st} for details).

If $l=k$, then $D^{red}_{A_0,\ldots,A_{k}}$ is the reduced
version of the resultant $R_{A_0,\ldots,A_{k}}$ (see Definition
\ref{defgkzres}); if $l=0$, then $D^{red}_{A_0}$ is the reduced
version of the discriminant $D_{A_0}$ (see Definition
\ref{defgkzdiscr}). In both of these cases, the discriminant set
$\Sigma_{A_0,\ldots,A_l}$is irreducible, and a combinatorial way
to verify $\codim\Sigma_{A_0,\ldots,A_l}=1$ is known (see
Corollary \ref{ldegree}(2) for $l=0$ and Proposition
\ref{essent}(1) for $l=k$). In general, however, unlike in these
two special cases, the set $\Sigma_{A_0,\ldots,A_l}$ may be not
irreducible and not even of pure dimension. Thus, both its
codimension 1 part $\{D^{red}_{A_0,\ldots,A_l}=0\}$ and its
higher codimension part
$\Sigma_{A_0,\ldots,A_l}\setminus\{D^{red}_{A_0,\ldots,A_l}=0\}$
may be non-empty for $0<l<k$ (see Example \ref{exanonequidim}
below). Nevertheless, restricting our attention to the
codimension 1 components, it turns out possible to express
$D^{red}_{A_0,\ldots,A_{l}}$ in terms of $A$-discriminants by
means of the Cayley trick, See Theorem \ref{univcayley} below.

\vspace{0.3cm} \textsc{Dual defectiveness.}
\begin{defin}\label{dualdefl} A collection of sets
$A_0,\ldots,A_l$ is said to be \textit{dual defect}, if the
closure of the set $\Sigma_{A_0,\ldots,A_l}$ is not a
hypersurface.
\end{defin}
One sufficient condition for non-dual-defectiveness is as follows:
\begin{utver}\label{equidim} If neither of the sets $A_0,\ldots,A_l$ is contained
in an affine hyperplane, and at least one of them is not dual
defect, then the collection $A_0,\ldots,A_l$ is not dual defect.
\end{utver}
Proof is given in Subsection \ref{sproofequidim} below, as well as
a more refined condition (Proposition \ref{equidim2}), which is
presumably a criterion. Note that neither of the conditions of
this statement can be omitted in general, as the following
examples demonstrate:
\begin{exa}\label{exanonequidim}
If $k=2$ and $A_0=A_1$ is the the set of vertices of the standard
2-dimensional simplex (which is dual defect), then
$\Sigma_{A_0,A_1}$ has codimension 2.

If $k=2,\; A_0=\{0,1,2\}\times\{0\}$ and
$A_1=\{0,1\}\times\{0,1\}$ (neither of these sets is dual defect,
but $\dim A_0<2$), then $\overline{\Sigma_{A_0,A_1}}$ has two
components. The first one consists of all pairs of polynomials of
the form $\Bigl(c(x-a)^2,\;
b_{11}xy+b_{01}x+b_{10}y+b_{00}\Bigr)$, and has codimension 1.
Another one consists of all pairs of polynomials of the form
$\Bigl(c_1(x-a_1)(x-a_2),\; c_2(x-a_1)(y-b)\Bigr)$, and has
codimension 2.
\end{exa}

\vspace{0.3cm} \textsc{Discriminant $B_{A_0,\ldots,A_l}$.} We
also consider another possible definition of discriminant, such
that the corresponding version of the non-dual defectiveness
assumption is weaker than the conventional one. Let $W$ be the
set of all collections
$(y,\varphi_0,\ldots,\varphi_l)\in\CC^k\times\C[A_0]\times\ldots\times\C[A_l]$,
such that $\varphi_0(y)=\ldots=\varphi_l(y)=0$. Let
$S_{A_0,\ldots,A_l}\subset\C[A_0]\times\ldots\times\C[A_l]$ be
the minimal (closed) set, such that the projection
$W\to\C[A_0]\times\ldots\times\C[A_l]$ is a fiber bundle outside
of $S$.
\begin{defin} \label{defbifuniv} The set $S_{A_0,\ldots,A_l}$ is called the \textit{bifurcation set}.
The collection $A_0,\ldots,A_l$ is said to be
\textit{$B$-nondegenerate}, if $S_{A_0,\ldots,A_l}$ is a
hypersurface. In this case the equation of $S_{A_0,\ldots,A_l}$
is denoted by $B_{A_0,\ldots,A_l}$ and is called the
\textit{bifurcation discriminant}.
\end{defin}
In contrast to the discriminant $D^{red}_{A_0,\ldots,A_l}$, the
bifurcation discriminant takes ``singularities of the system
$\varphi_0=\ldots=\varphi_l=0$ at infinity" into account.
\begin{exa} In the notation of the previous example ($k=2,\;
A_0=\{0,1,2\}\times\{0\},\; A_1=\{0,1\}\times\{0,1\}$), despite
the discriminant set $\overline{\Sigma_{A_0,A_1}}$ is not of pure
dimension, the bifurcation set
$S_{A_0,A_1}\supset\overline{\Sigma_{A_0,A_1}}$ is a hypersurface
that consists of five components. A generic point in each of
these components is as follows (the codimension 2 component of
the discriminant set is the intersection of the last two
components):

$c(x-a)^2,\; b_{11}xy+b_{10}y+b_{01}x+b_{00}$;

$a_1x^2-a_2x,\; b_{11}xy+b_{10}y+b_{01}x+b_{00}$;

$a_1x-a_2,\; b_{11}xy+b_{10}y+b_{01}x+b_{00}$;

$c_1(x-a_1)(x-a_2),\; c_2(xy-a_1y+b_{01}x+b_{00})$;

$c_1(x-a_1)(x-a_2),\; c_2(b_{11}xy+b_{10}y+x-a_1)$.

In particular, we have $$B_{A_0,A_1}=D^{red}_{A_0,A_1}\cdot
D^{red}_{\{(0,0)\},A_1}\cdot D^{red}_{\{(2,0)\},A_1}\cdot
D^{red}_{A_0,\{0,1\}\times\{1\}}\cdot
D^{red}_{A_0,\{0,1\}\times\{0\}}.$$ These observations generalize
as follows.
\end{exa}
\begin{conj} \label{conj1} All collections are
$B$-nondegenerate.
\end{conj}
\begin{utver} \label{conj1sim} If the convex hulls of finite sets $A_0,\ldots,A_l$ in $\Z^k$ have the
same dual fan, then the collection $A_0,\ldots,A_l$ is
$B$-nondegenerate.
\end{utver}
In particular, if a collection consists of one set, then it is
$B$-nondegenerate, in contrast to dual-defectiveness. The proof
is given in Subsection \ref{sproofconj1sim}.
\begin{lemma}\label{bifdiscr} The bifurcation discriminant
$B_{A_0,\ldots,A_l}$ equals the least common multiple of the
discriminants $D_{A'_0,\ldots,A'_l}$ over all compatible
collections of faces $A'_0\subset A_0,\ldots,A'_l\subset
A_l$.\end{lemma} We omit the proof, since it follows by
definitions. In the next subsection, we explicitly decompose the
discriminants $D_{A'_0,\ldots,A'_l}$ into irreducible factors and
compute the desired least common multiple (Corollary
\ref{bifcayley}).

\subsection{Cayley trick.}\label{scayley}

Let $e_0,\ldots,e_l$ be the standard basis in $\Z^{l+1}$. For
$J\subset\{0,\ldots,l\}$, denote the set $\bigcup\limits_{j\in J}
A_j\times\{e_j\}$ by $A_J\subset\Z^k\oplus\Z^{l+1}$. We identify
the space $\bigoplus_{j\in J} \C[A_j]$ with the space $\C[A_J]$ by
identifying a collection of polynomials $\varphi_j\in\C[A_j],\;
j\in J,$ on the complex torus $\CC^n$ with the polynomial
$\sum_{j\in J}\lambda_j\varphi_j$ on $\CC^k\times\CC^{l+1}$, where
$\lambda_0,\ldots,\lambda_l$ are the standard coordinates on
$\CC^{l+1}$. This identification allows us to regard the
discriminant $D^{red}_{A_J}$ as a polynomial on the space
$\C[A_0]\oplus\ldots\oplus\C[A_l]$.

\begin{theor}[Cayley trick]\label{univcayley} The discriminant $D^{red}_{A_0,\ldots,A_l}$ equals
the product of the discriminants $D^{red}_{A_J}$ over all subsets
$J\subset\{0,\ldots,l\}$, such that $\codim J\leqslant\codim J'$
for every $J'\supset J$ (recall that $\codim J$ stands for the
difference $\dim\sum_{j\in J}A_j-|J|$).\end{theor}

For $l=k$, this is proved in \cite{gkz}; the only admissible $J$
is $\{0,\ldots,k\}$ in this case. If neither of $A_i$ is
contained in an affine hyperplane, then $\{0,\ldots,k\}$ is the
only admissible $J$ as well. Lemma \ref{bifdiscr} leads to a
similar formula for the bifurcation discriminant:
\begin{sledst} \label{bifcayley} The bifurcation  discriminant $B_{A_0,\ldots,A_l}$
equals the product of the discriminants
$D^{red}_{A'_{\{j_1,\ldots,j_p\}}}$ over all collections of
compatible faces $A'_{j_1}\subset A_{j_1},\ldots,A'_{j_p}\subset
A_{j_p}$ that can be extended to a collection of compatible faces
$A'_0\subset A_0,\ldots,A'_l\subset A_l$ such that
$\dim\sum_{j\in J}A'_j - \dim\sum_iA'_{j_i}\geqslant |J|-p$ for
every $J\supset\{j_1,\ldots,j_p\}$.
\end{sledst}

The proof of Theorem \ref{univcayley} is given at the end of this
subsection and is based on the following construction.

\begin{defin} We define $\Sigma_{\{j_0,\ldots,j_q\}}$ as the set of all
collections\newline
$(\varphi_0,\ldots,\varphi_l)\subset\C[A_0]\oplus\ldots\oplus\C[A_l]$,
such that
\newline
1) $\varphi_0(x)=\ldots=\varphi_l(x)=0$ for some $x\in\CC^k$, and
\newline
2) $\lambda_{j_0}{\rm d}\varphi_{j_0}(x)+\ldots+\lambda_{j_q}{\rm
d}\varphi_{j_q}(x)=0$ for some
$(\lambda_{j_0},\ldots,\lambda_{j_q})\in\CC^{q+1}$.\end{defin}

We have two options for each set $\Sigma_J$:
\begin{lemma} \label{cayleyalt}
A) If $A_J$ is not dual defect, and $\codim J\leqslant\codim J'$
for every $J'\supset J$, then the closure of $\Sigma_J$ is
defined by the equation $D^{red}_{A_J}=0$.\newline B) Otherwise,
$\codim\Sigma_J>1$.
\end{lemma}

The proof of this lemma is given below and is based on the
following important fact:
\begin{lemma}[\cite{kh77}]\label{bkh1} Generic polynomials $\psi_i\in\C[A_i],;
i=0,\ldots,l$, have a common root in $\CC^n$ if an only if $\dim
A_{i_1}+\ldots+A_{i_q}\geqslant q$ for every sequence $0\leqslant
i_1<\ldots<i_q\leqslant l$.
\end{lemma}
This fact is mentioned as obvious in \cite{kh77}, but we prefer
to give a proof for the sake of completeness.

\textsc{Proof.} If $l\geqslant k$ then the statement is obvious,
because, on one hand, generic polynomials $\psi_0,\ldots,\psi_l$
have no common roots and, on the other hand, we have $\dim
A_0+\ldots+A_l<l+1$. The case $l<k-1$ can be reduced to the case
$l=k-1$ by introducing arbitrary finite sets
$A_{l+1},\ldots,A_{k-1}$, whose convex hulls are $k$-dimensional,
and considering generic polynomials $\psi_i\in\C[A_i],\;
i=0,\ldots,k-1$. Finally, if $l=k-1$, then the number of common
roots of generic polynomials $\psi_0,\ldots,\psi_{k-1}$ equals
$k!$ times the mixed volume of the convex hulls of the sets
$A_0,\ldots,A_{k-1}$, which is non-zero under the assumption of
the lemma by Lemma \ref{zeromv}. $\quad\Box$

\vspace{0.3cm} \textsc{Proof of Lemma \ref{cayleyalt}.} First,
note that the image of the set $\Sigma_J$ under the natural
projection $\C[A_0]\oplus\ldots\oplus\C[A_l]\to\bigoplus_{j\in
J}\C[A_j]$ is contained in $\Sigma_{A_J}$. In particular, if
$\codim\Sigma_{A_J}>1$, then the set $\Sigma_J$ satisfies both of
the statements (A) and (B) of Lemma \ref{cayleyalt},
independently of $\codim J'$ for $J'\supset J$. Thus, we can
assume that $\codim\Sigma_{A_J}=1$, i.e. $A_J$ is not dual defect.

Consider the vector space $L\subset\R^n$, parallel to the affine
span of the sum $\sum_{j\in J}A_j$, denote the projection
$\R^k\to\R^k/L$ by $p_L$, and consider the torus $T_L=\{ x\, |\,
x^a=1 \mbox{ for } a\in L\}\subset\CC^k$. Then, for a generic
polynomial in $\C[A_j]$, its restriction to $T_L$ is a generic
polynomial in $\C[p_LA_j]$. Thus, by Lemma \ref{bkh1}, generic
polynomials $\varphi_i\in\C[A_i],\; i\notin J,$ have a common
zero on a torus $c\cdot T_L,\; c\in\CC^k$, if and only if every
subset $I\subset\{0,\ldots,l\}\setminus J$ satisfies inequality
$\dim \sum_{i\in I}p_L(A_i)\geqslant|I|$, i.e. $\codim I\cup
J\geqslant\codim J$.

\vspace{0.3cm} \textsc{Proof of Part A.} For
$J=\{j_0,\ldots,j_q\}$ and a polynomial
$\lambda_{j_0}\varphi_{j_0}+\ldots+\lambda_{j_q}\varphi_{j_q})\in\C[A_J]$,
the set
$$\{x\in\CC^k\; |\; \varphi_{j_0}(x)=\ldots=\varphi_{j_q}(x)=0 \mbox{ and } \lambda_{j_0}{\rm d}\varphi_{j_0}(x)+\ldots+\lambda_{j_q}{\rm
d}\varphi_{j_q}(x)=0 $$ $$\mbox{ for some }
(\lambda_{j_0},\ldots,\lambda_{j_q})\in\CC^{q+1}\}$$ is non-empty
and preserved under multiplication by elements of $T_L$, thus it
contains a torus $c\cdot T_L,\; c\in\CC^k$. Thus, generic
polynomials $\varphi_i\in\C[A_i],\; i\notin J,$ have a common
zero on it under the assumption of Part A. Thus, every fiber of
the projection $\Sigma_J\to\Sigma_{A_J}$ is Zariski open in the
corresponding fiber of the ambient projection
$\C[A_0]\oplus\ldots\oplus\C[A_l]\to\bigoplus_{j\in J}\C[A_j]$.
Thus, $\codim\Sigma_J=\codim\Sigma_{A_J}=1$.

\vspace{0.3cm} \textsc{Proof of Part B.} For
$J=\{j_0,\ldots,j_q\}$ and a generic polynomial
$\lambda_{j_0}\varphi_{j_0}+\ldots+\lambda_{j_q}\varphi_{j_q}$ in
$\C[A_J]$, the set
$$\{x\in\CC^k\; |\; \varphi_{j_0}(x)=\ldots=\varphi_{j_q}(x)=0 \mbox{ and } \lambda_{j_0}{\rm d}\varphi_{j_0}(x)+\ldots+\lambda_{j_q}{\rm
d}\varphi_{j_q}(x)=0 $$ $$\mbox{ for some }
(\lambda_{j_0},\ldots,\lambda_{j_q})\in\CC^{q+1}\}$$ consists of
finitely many tori $c_s\cdot T_L,\; c_s\in\CC^k,\; s=1,\ldots S$,
because $A_J$ is not dual defect. Thus, generic polynomials
$\varphi_i\in\C[A_i],\; i\notin J,$ have no common zero on it
under the assumption of Part B. Thus, a generic fiber of the
projection $\Sigma_J\to\Sigma_{A_J}$ has codimension 1 or greater
in the corresponding fiber of the ambient projection
$\C[A_0]\oplus\ldots\oplus\C[A_l]\to\bigoplus_{j\in J}\C[A_j]$.
Thus, $\codim\Sigma_J\geqslant\codim\Sigma_{A_J}>1$. $\quad\Box$

\vspace{0.3cm} \textsc{Proof of Theorem \ref{univcayley}.} The
desired set $\Sigma_{A_0,\ldots,A_l}$ is the union of the sets
$\Sigma_J$ over all $J\subset\{0,\ldots,l\}$. By Lemma
\ref{cayleyalt}, the product of the discriminants, mentioned in
the formulation of the theorem, vanishes at the union of
codimension 1 components of the closure
$\overline{\Sigma_{A_0,\ldots,A_l}}$. Since all non-constant
discriminants in this product are irreducible and distinct by
Lemmas \ref{lirred} and \ref{lnontriv}, then this product is a
square-free polynomial and thus equals
$D^{red}_{A_0,\ldots,A_l}$. $\quad\Box$

\subsection{Proof of Proposition \ref{equidim}.}\label{sproofequidim}

The proof relies upon Lemma \ref{cayleyalt} and relevant notation
from the previous subsection. We split Proposition \ref{equidim}
into the following two lemmas.
\begin{lemma}\label{eql1} If $A_0$ is not dual defect, and $\dim A_i=k$ for every
$i=0,\ldots,l$, then $A_{\{0,\ldots,l\}}$ is not dual defect.
\end{lemma}
\begin{lemma} \label{eql2} If $A_{\{0,\ldots,l\}}$ is not dual defect, and $\dim A_i=k$ for every
$i=0,\ldots,l$, then $\Sigma_{A_0,\ldots,A_l}$ is irreducible of
codimension 1.
\end{lemma}
Proof of both of these lemmas is given below and is based on the
following fact.
\begin{lemma} \label{eql3} If $\dim A=k,\; \varphi\in\C[A],\; \varphi(y_0)=0,$ and a vector
$v\in\C^n$ is close enough to ${\rm d}\varphi(y_0)$, then there
exists $\tilde\varphi\in\C[A]$ near $\varphi$, such that
$\tilde\varphi(y_0)=0$ and ${\rm d}\tilde\varphi(y_0)=v$.
\end{lemma}
\textsc{Proof.} One can readily verify this statement if $A$ is
of cardinality $k+1$ and $\varphi=0$. Thus, if $A_0\subset A$ is
the set of vertices of a $k$-dimensional simplex, then there
exists a small $\psi\in\C[A_0]$ such that ${\rm
d}\psi(y_0)=v-{\rm d}\varphi(y_0)$, and we can set
$\tilde\varphi=\varphi+\psi$. $\quad\Box$

\vspace{0.3cm} \textsc{Proof of Lemma \ref{eql1}.} Since dual
defectiveness is preserved by parallel translations, we can assume
without loss of generality that $0\in A_0$. By Lemma
\ref{ldimsing}, a generic polynomial $\varphi_0\in\Sigma_{A_0}$
has an isolated singular point $y_0$. By Lemma \ref{eql3},
generic polynomials $\varphi_i\in\C[A_i]$, such that
$\varphi_i(y_0)=0$, define a non-degenerate complete intersection
$\varphi_1=\ldots=\varphi_l=0$, passing through $y_0$ and
transversal to the hypersurface $\varphi_0=0$ in a punctured
neighborhood of $y_0$. Then, for generic
$\tilde\varphi_i\in\C[A_i]$ near $\varphi_i,\; i=0,\ldots,l$,
there exists a small number $\varepsilon$ such that the
non-degenerate complete intersection
$\tilde\varphi_1=\ldots=\tilde\varphi_l=0$ is tangent to the
smooth hypersurface $\tilde\varphi_0=\varepsilon$ at some point
near $y_0$, which implies that
$(\tilde\varphi_0-\varepsilon,\tilde\varphi_1,\ldots,\tilde\varphi_l)\in\Sigma_{A_0,\ldots,A_l}$.
Thus, $\codim\Sigma_{A_0,\ldots,A_l}=1$ near the point
$(\varphi_0,\ldots,\varphi_l)$.

Since $\codim\Sigma_{A_0,\ldots,A_l}=1$ at some point of this set,
$$\Sigma_{A_0,\ldots,A_l}=\bigcup_{J\subset\{0,\ldots,l\}}\Sigma_J,$$
$$\codim\Sigma_J>1 \mbox{ for } J\varsubsetneq \{0,\ldots,l\} \mbox{ by Lemma \ref{cayleyalt}(B),}$$
$$\mbox{ and } \Sigma_{\{0,\ldots,l\}}=\Sigma_{A_{\{0,\ldots,l\}}}, $$
we have $\codim \Sigma_{A_{\{0,\ldots,l\}}}=1$, thus
$A_{\{0,\ldots,l\}}$ is not dual defect. $\quad\Box$

\vspace{0.3cm} \textsc{Proof of Lemma \ref{eql2}.} For every
collection
$(\varphi_0,\ldots,\varphi_l)\in\Sigma_{A_0,\ldots,A_l}$, there
exists a point $y_0$, such that
$\varphi_0(y_0)=\ldots=\varphi_l(y_0)=0$ and ${\rm
d}\varphi_0(y_0),\ldots,{\rm d}\varphi_l(y_0)$ are linearly
dependent. There exist vectors $v_0,\ldots,v_l$ near ${\rm
d}\varphi_0(y_0),\ldots,{\rm d}\varphi_l(y_0)$, such that
$\sum_{j=0}^l\lambda_jv_j=0$ with $\lambda_j\ne 0$ for
\textit{every} $j=0,\ldots,l$. Thus, by Lemma \ref{eql3}, there
exists a collection $(\tilde\varphi_0,\ldots,\tilde\varphi_l)$
near $(\varphi_0,\ldots,\varphi_l)$, such that
$$\sum_{j=0}^l\lambda_j{\rm d}\tilde\varphi_j(y_0)=0$$
with $\lambda_j\ne 0$ for every $j=0,\ldots,l$. This means that
$(\lambda_0,\ldots,\lambda_l,y_0)\in\Sigma_{A_{\{0,\ldots,l\}}}$.
Thus, $\Sigma_{A_0,\ldots,A_l}$ contains the irreducible
codimension 1 set $\Sigma_{A_{\{0,\ldots,l\}}}$ and is contained
in its closure, which completes the proof. $\quad\Box$

\vspace{0.3cm} \textsc{Refinement of Proposition \ref{equidim}.}
Consider finite sets $A_0,\ldots,A_l$ in $\Z^k$ and a subset
$J\subset\{0,\ldots,l\}$.
\begin{defin} $J$ is said to be \textit{reliable}, if $A_J$ is
not dual defect, and $\codim J\leqslant\codim J'$ for every
$J'\supset J$ (recall that $\codim J$ stands for the difference
$\dim\sum_{j\in J}A_j-|J|$).
\end{defin}
Note that $J$ is reliable if and only if $\Sigma_J$ is a
hypersurface.
\begin{defin} Vectors $v_0,\ldots,v_l$ are said to be $J$-linearly
dependent, if $\sum_{j\in J}\lambda_j v_j=0$ with $\lambda_j\ne
0$ for every $j\in J$.
\end{defin}
For a set $A\subset\R^k$, denote the space of all linear
functions, whose restrictions to $A$ are constant, by
$A^{\bot}\subset(\R^k)^*$. The following refined version of
Proposition \ref{equidim} is presumably a criterion for non-dual
defectiveness of a collection of sets.
\begin{utver} \label{equidim2} A collection of finite sets $A_0,\ldots,A_l$ in
$\Z^k$ is not dual defect, if $\codim J\geqslant -1$ for every
$J\subset\{0,\ldots,l\}$, and there exists a $J$-linearly
dependent collection $\tilde v_0\in A_0^{\bot},\ldots,\tilde
v_l\in A_l^{\bot}$ with a reliable $J$ in every neighborhood of
every linearly dependent collection of vectors $v_0\in
A_0^{\bot},\ldots,v_l\in A_l^{\bot}$.
\end{utver}

We can formulate Lemma \ref{eql3} with no assumption $\dim A=k$ as
follows. Denote the natural identification $T_y\CC^n\to
T_{(1,\ldots,1)}\CC^n$ by $e_y$.
\begin{lemma} \label{l2sim} Suppose that $\varphi\in\C[A],\; \varphi(y)=0,$ and $v\in\C^n$ is
near ${\rm d}\varphi(y)$. Then every neighborhood of $\varphi$
contains $\tilde\varphi\in\C[A]$ with $\tilde\varphi(y)=0$ and
${\rm d}\tilde\varphi(y)=v$ if and only if $e_y(v)\in A^{\bot}$.
\end{lemma}
\textsc{Proof of Proposition \ref{equidim2}.} For a collection
$(\varphi_0,\ldots,\varphi_l)\subset\Sigma_{A_0,\ldots,A_l}$,
choose a point $y$, such that
$\varphi_0(y)=\ldots=\varphi_l(y)=0$ and the differentials ${\rm
d}\varphi_0(y),\ldots,{\rm d}\varphi_l(y)$ are linearly
dependent, and choose nearby vectors $v_0,\ldots,v_l$ to be
$J$-linearly dependent for a reliable $J$. By Lemma \ref{l2sim},
we have $\tilde\varphi_i(y)=0$ and ${\rm d}\tilde\varphi_i(y)=v_i$
for some $\tilde\varphi_i$ near $\varphi_i$, thus the collection
$(\tilde\varphi_0,\ldots,\tilde\varphi_l)$ is contained in the
hypersurface $\Sigma_J$. $\quad\Box$

\subsection{Proof of Proposition \ref{conj1sim}.}\label{sproofconj1sim}
See Subsection \ref{scayley} for the definition of $A_J$ and
$\Sigma_J$ for $J\subset\{0,\ldots,l\}$.

We prove the statement by induction on $k$. For a subset
$A'\subset A$, there is a natural projection from $\C[A]$ to
$\C[A']$ that assigns the polynomial $\sum_{a\in A'}c_ay^a$ to a
polynomial $\sum_{a\in A}c_ay^a$. For subsets $A'_0\subset
A_0,\ldots,A'_l\subset A_l$, we denote the preimage of the set
$B\subset\C[A'_0]\oplus\ldots\oplus\C[A'_l]$ under this
projection by $\widetilde B$.

The bifurcation set $S_{A_0,\ldots,A_l}$ contains the set
$\widetilde S_{A'_0,\ldots,A'_l}$ for every collection of
compatible faces $A'_0\subset A_0,\ldots,A'_l\subset A_l$. The
difference of $S_{A_0,\ldots,A_l}$ and
$\bigcup_{A'_0,\ldots,A'_l}\widetilde S_{A'_0,\ldots,A'_l}$ is
contained in the set $\Sigma_{A_0,\ldots,A_l}$, which is the union
of irreducible sets $\Sigma_J,\; J\subset\{0,\ldots,l\}$. Thus,
we can reformulate Proposition \ref{conj1sim} as follows: every
set of the form $\Sigma_J$ is contained in the closure of a set
of the form $\Sigma_J$ or $\widetilde S_{A'_0,\ldots,A'_l}$, whose
codimension is 1.

For a generic point $(\varphi_0,\ldots,\varphi_l)\in\Sigma_J$, we
have the following three cases: \newline 1) The set of all
$y\in\CC^k$, such that ${\rm d}\varphi_0(y),\ldots,{\rm
d}\varphi_l(y)$ are linearly dependent, has positive dimension.
Then it contains a germ of a curve $\CC\to\CC^k$, whose leading
term is $(c_1t^{\gamma_1},\ldots,c_kt^{\gamma_k})\ne{\rm const}$.
The covector $(\gamma_1,\ldots,\gamma_k)\ne 0$, as a function on
$A_i$, attains its maximum at some proper face $A'_i\subsetneq
A_i$. Thus, $(\varphi_0,\ldots,\varphi_l)\in\widetilde
S_{A'_0,\ldots,A'_l}$, which is a hypersurface by induction.
\newline 2) $J=\{0,\ldots,l\}$, and there exists an isolated point
$y\in\CC^k$, such that ${\rm d}\varphi_0(y),\ldots,{\rm
d}\varphi_l(y)$ are linearly dependent. Then we can choose
vectors $v_i$ near ${\rm d}\varphi_i(y)$, such that the
non-trivial linear combination of $v_0,\ldots,v_l$ is unique and
has all non-zero coefficients. By Lemma \ref{eql3}, we can
perturb the collection $(\varphi_0,\ldots,\varphi_l)$ into
$(\tilde\varphi_0,\ldots,\tilde\varphi_l)\in\Sigma_J$, with
$\tilde\varphi_i(y)=0$ and ${\rm d}\tilde\varphi_i(y)=v_i$. Thus,
the polynomial
$\sum_{i=0}^l\lambda_i\tilde\varphi_i\in\Sigma_{A_{J}}$ has an
isolated line of singular points in its zero set, and
$\Sigma_{A_{J}}=\Sigma_J$ is a hypersurface by Lemma
\ref{ldimsing}. \newline  3) In the general case, in the same way
as above, we can perturb the collection
$(\varphi_0,\ldots,\varphi_l)$ into
$(\tilde\varphi_0,\ldots,\tilde\varphi_l)\in\Sigma_{\{0,\ldots,l\}}$.
But $\Sigma_{\{0,\ldots,l\}}$ is either a hypersurface itself, or
is contained in a hypersurface $\widetilde S_{A'_0,\ldots,A'_l}$
(see Cases 1 and 2). $\quad\Box$

\subsection{Proof of GKZ decomposition formula.}\label{sproofdecomp}

To deduce Gelfand-Kapranov-Zelevinsky's decomposition formula
(Proposition \ref{gkzdecomp}) from the relative
Kouchnirenko-Bernstein formula, we reformulate the definition of
the $A$-determinant as follows (we use notation and facts from
Subsection \ref{smultci}).

\vspace{0.3cm} \textsc{Geometric characterization of $A$-resultant
and $A$-determinant.} Denote the projection
$\CC^k\times\C[A_0]\oplus\ldots\oplus\C[A_k]\to\C[A_0]\oplus\ldots\oplus\C[A_k]$
by $p$. Let $R_i$ be the tautological polynomial on
$\CC^k\times\C[A_0]\oplus\ldots\oplus\C[A_k]$ that assigns the
number $\varphi_i(y)$ to a point
$(y,\varphi_0,\ldots,\varphi_k)\in\CC^k\times\C[A_0]\oplus\ldots\oplus\C[A_k]$.
Denote the convex hull of $A_i$ by $B_i$, and the dual fan of
$B_0+\ldots+B_k$ by $\Sigma$. Then $s_{B_i}\cdot R_i$ extends to
a section $\widetilde R_i$ of the line bundle $I_{B_i}$ on the
product $\T^{\Sigma}\times\C[A_0]\oplus\ldots\oplus\C[A_k]$ (see
Subsection \ref{srelbernst} for the notation $I_B$ and $s_B$), and
we can reformulate the definition of the resultant
$R_{A_0,\ldots,A_k}$ as follows.
\begin{lemma}\label{gmdet1}
$[R_{A_0,\ldots,A_k}=0]=p_*[\widetilde R_0=\ldots=\widetilde
R_k=0]$.
\end{lemma}
We also need a similar description for the $A$-determinant.
Denote the projection $\CC^k\times\C[A]\to\C[A]$ by $p$. Let
$S_0$ be the tautological polynomial on $\CC^k\times\C[A]$ that
assigns the number $\varphi(y)$ to a point
$(y,\varphi)\in\CC^k\times\C[A]$, and let $S_i$ be
$y_i\frac{\partial S_0}{\partial y_i}$, where $y_1,\ldots,y_k$
are the standard coordinates on $\CC^k$. Denote the convex hull
of $A$ by $B$ and its dual fan by $\Sigma$. Then $s_{B}\cdot S_i$
extends to a section $\widetilde S_i$ of the line bundle $I_{B}$
on the product $\T^{\Sigma}\times\C[A]$, and we can reformulate
the definition of the $A$-determinant as follows.
\begin{lemma} \label{lgmdet} $[E_A=0]=p_*[\widetilde
S_0=\ldots=\widetilde S_k=0]$.
\end{lemma} \textsc{Proof.} Let $j$ be the inclusion
$\C[A]\hookrightarrow\underbrace{C[A]\oplus\ldots\oplus\C[A]}_{k+1}$
that assigns the collection
$(\varphi,y_1\frac{\partial\varphi}{\partial
y_1},\ldots,y_k\frac{\partial\varphi}{\partial y_k})$ to every
$\varphi\in\C[A]$, and denote the induced inclusion
$\CC^k\times\C[A]\hookrightarrow\CC^k\times\underbrace{C[A]\oplus\ldots\oplus\C[A]}_{k+1}$
by the same letter $j$. Then $$p_*[\widetilde
S_0=\ldots=\widetilde S_k=0]=p_*j^*[\widetilde
R_0=\ldots=\widetilde R_k=0]=$$ $$=j^*p_*[\widetilde
R_0=\ldots=\widetilde R_k=0]=j^*[R_{A_0,\ldots,A_k}=0]=[E_A=0],$$
where the last two equalities are by Lemma \ref{gmdet1} and by
definition of the $A$-determinant respectively. $\quad\Box$

\vspace{0.3cm} \textsc{Proof of Proposition \ref{gkzdecomp}.} In
the notation of Lemma \ref{lgmdet}, represent the complete
intersection $[\widetilde S_0=\ldots=\widetilde S_k=0]$ as a
linear combination of irreducible varieties $\sum a_iV_i$. For
every $V_i$, let $\T_i$ be the minimal orbit of the toric variety
$\T^{\Sigma}$, such that $V_i$ is contained in the closure of
$\T_i\times\C[A]$. For an arbitrary face $A'$ of the set $A$, we
denote the corresponding orbit of $\T^{\Sigma}$ by $\T_{A'}$, and
denote the sum $\sum_{i\; |\; \T_i=\T_{A'}}a_iV_i$ by $V_{A'}$.
To prove the equality $[E_A=0]=\sum_{A'\subset
A}c^{A',A}[D_{A'}=0]$ (which is exactly the statement of
Proposition \ref{gkzdecomp}), it is enough to prove the following
lemma.
\begin{lemma} $p_*(V_{A'})=c^{A',A}[D_{A'}=0]$.
\end{lemma}
\textsc{Proof.} Denote $\dim A'$ by $l$, and choose a generic
real $(k+1)\times(k+1)$-matrix $M$, whose first $l+1$ rows
generate the vector span of the set $\{1\}\times
A'\subset\Z\oplus\Z^k$. Denote the entries of the product
$M\cdot(\widetilde S_0,\ldots,\widetilde S_k)^T$ by
$Z_0,\ldots,Z_k$ (recall that they are sections of the line
bundle $I_{B}$). The desired equality is a corollary of the
following facts:

1)
$[E_A=0]=p_*\bigl([Z_0=\ldots=Z_{l-1}=0]\cap[Z_l=\ldots=Z_k=0]\bigr)$
by Lemma \ref{lgmdet}.

2) The orbit $\T_{A'}$ is a component of the complete intersection
$[Z_0=\ldots=Z_{l-1}=0]$ of multiplicity $c^{A',A}$ by Lemma
\ref{orbmult}.

3) The complete intersection $[Z_l=\ldots=Z_k=0]$ intersects the
orbit $\T_{A'}$ transversally, and
$p_*([Z_l=\ldots=Z_k=0]\cap\T_{A'})=[D_{A'}=0]$ by definition of
the $A$-discriminant. $\quad\Box$

\section{Eliminants.}

In the first subsection, we formulate a local version of
elimination theory in the context of Newton polyhedra (the global
version is presented in \cite{ekh}), i.e. we study the Newton
polyhedron and leading coefficients of the equation of a
projection of a complete intersection which is defined by
equations with given Newton polyhedra and generic leading
coefficients. The main result is Theorem \ref{elimth}, the proof
is given in Subsection \ref{sproofelimth}. In Subsection
\ref{srdiscrim}, we specialize this to $A$-determinants.

\subsection{Elimination theory.}\label{selimth}

\textsc{Notation.} Let $\tau\subsetneq(\R^n)^*$ be a convex
$n$-dimensional rational polyhedral strictly convex cone (i.e. a
cone that does not contain a line), denote its dual cone
$\{v\in\R^n\; |\; \gamma(v)\geqslant 0 \mbox{ for }
\gamma\in\tau\}$ by $\tau^{\vee}$, and consider the corresponding
affine toric variety $\T^{\tau}=\spec\C[\tau^{\vee}]$ with the
maximal torus $\CC^n$ and the vertex $O$. Note that $\tau^{\vee}$
is unbounded. A germ of a meromorphic function on $\T^{\tau}$
near $O$ with no poles in the maximal torus can be represented as
a power series $f(x)=\sum_{b\in\Z^n}c_bx^b$ for $x\in\CC^n$, and
the convex hull of the set $\{b\; |\; c_b\ne 0\}+\tau^{\vee}$ is
called \textit{the Newton polyhedron} $\Delta_f$ of $f$. The
union of all bounded faces of $\Delta_f$ is called \textit{the
Newton diagram} $\partial\Delta_f$, and the coefficients $c_b,\;
b\in\partial\Delta_f$, are called the \textit{leading
coefficients} of the germ $f$.

\vspace{0.3cm} \textsc{Eliminant.} Let $A_0,\ldots,A_k$ be finite
sets in $\Z^k$. For $i=0,\ldots,k$ and $a\in A_i$, let $f_{a,i}$
be a germ of a meromorphic function on the toric variety
$\T^{\tau}$ with no poles in the maximal torus. We define the
germs of functions $F_0,\ldots,F_k$ on $\T^{\tau}\times\CC^k$ by
the formula
$$F_i(x,y)=\sum_{a\in A_i} f_{a,i}(x)y^a \mbox{ for }
x\in\T^{\tau},\; y\in\CC^k,$$  note that $F_i(x,\cdot)\in\C[A_i]$,
and denote the number
$R_{A_0,\ldots,A_k}\Bigl(F_0(x,\cdot),\ldots,F_k(x,\cdot)\Bigr)$
by $R_{F_0,\ldots,F_k}(x)$ (see Subsection \ref{sadeterm} for the
definition of the resultant $R_{A_0,\ldots,A_k}$).

\begin{defin} The function $R_{F_0,\ldots,F_k}$ on the toric variety
$\T^{\tau}$ is called \textit{the eliminant of the projection of
the complete intersection $F_{0}=\ldots=F_k=0$ to} $\T^{\tau}$.
\end{defin}

Geometric meaning of the function $R_{F_0,\ldots,F_k}$ is as
follows: if leading coefficients of the functions $f_{a,i}$ are
in general position, and the image of the complete intersection
$F_0=\ldots=F_k=0$ under the projection
$\T^{\tau}\times\CC^k\to\T^{\tau}$ has codimension 1, then the
closure of this image is the zero locus of $R_{F_0,\ldots,F_k}$
(see the beginning of the next subsection for details). Under
these assumptions, the Newton polyhedron of $R_{F_0,\ldots,F_k}$
does not depend on coefficients of the functions $f_{a,i}$, but
only on their Newton polyhedra. Theorem \ref{elimth} below solves
the following problem:

\begin{quote}
Express the Newton polyhedron and leading coefficients of the
eliminant $R_{F_0,\ldots,F_k}$ in terms of the Newton polyhedra
and leading coefficients of the functions $f_{a,i}$, provided that
leading coefficients are in general position, and describe this
condition of general position explicitly.
\end{quote}

\textsc{Condition of general position.} We denote the Newton
polyhedron of $f_{a,i}$ by $\Delta_{a,i}$, and define the Newton
polyhedron $\Delta_i$ of the function $F_i$ as the convex hull of
the set $\bigcup_{a\in
A_i}\Delta_{a,i}\times\{a\}\subset\R^n\oplus\R^k$; then $F_i(z)$
can be represented as a power series
$\sum_{b\in\Delta_i}c_{b,i}z^b$ for $z\in\CC^n\times\CC^k$. If
$\Gamma$ is a face of $\Delta_i$, then we denote the function
$\sum_{b\in\Gamma}c_{b,i}z^b$ by $F_i^{\Gamma}$.
\begin{defin} \label{genposres}
The leading coefficients of the functions $F_0,\ldots,F_k$ are
said to be \textit{in general position}, if, for every collection
of compatible bounded faces $\Gamma_{i}\subset\Delta_{i},\;
i=0,\ldots,k$ (see Definition \ref{defcompat}), such that the
restriction of the projection $\R^n\oplus\R^k\to\R^k$ to
$\Gamma_{0}+\ldots+\Gamma_{k}$ is injective, the system of
polynomial equations $F_0^{\Gamma_0}=\ldots=F_k^{\Gamma_k}=0$ has
no solutions in $\CC^n\times\CC^k$.
\end{defin}
This condition is obviously satisfied for generic leading
coefficients of the functions $f_{a,i}$. Note that this condition
is slightly weaker than the one in \cite{ekh}. For example, if a
face $\Gamma_i\subset\Delta_i$ is contained in a fiber of the
projection $\R^n\oplus\R^k\to\R^k$, then we do not impose any
assumptions on the leading coefficients
of $F_i$, corresponding to internal integer points of $\Gamma_i$. 
This slight difference is important for our purpose (cf.
\cite{gp}), see the proof of Proposition \ref{elimthunmix} below.

\vspace{0.3cm} \textsc{Elimination theorem.} We recall that the
Minkowski integral $\int\Delta\subset\R^n$ of a polyhedron
$\Delta\subset\R^n\oplus\R^k$ is defined in Section 1. For
bounded faces $\Gamma_i\subset\Delta_i,\, i=0,\ldots,k$, we
denote the intersection of $A_i$ with the image of $\Gamma_i$
under the projection $\R^n\oplus\R^k\to\R^k$ by
$\widetilde\Gamma_i$, and the value
$R_{\widetilde\Gamma_0,\ldots,\widetilde\Gamma_k}\Bigl(F^{\Gamma_0}_0(x,\cdot),\ldots,F^{\Gamma_k}_k(x,\cdot)\Bigr)$
by $R_{F^{\Gamma_0}_0,\ldots,F^{\Gamma_k}_k}(x)$ for $x\in\CC^n$,
then the Laurent polynomial
$R_{F^{\Gamma_0}_0,\ldots,F^{\Gamma_k}_k}$ on $\CC^n$ depends
only on leading coefficients of the functions $f_{a,i}$.

\begin{theor}\label{elimth} 1) The Newton polyhedron of the eliminant
$R_{F_0,\ldots,F_k}$ is contained in the mixed fiber polyhedron
$$\MF\bigl(\Delta_0,\ldots,\Delta_k\bigr).$$
These two polyhedra coincide if and only if the leading
coefficients of $F_0,\ldots,F_k$ are in general position in the
sense of Definition \ref{genposres}.
\newline
2) For every face $\Gamma$ of the polyhedron
$\MF\bigl(\Delta_0,\ldots,\Delta_k\bigr)\subset\R^n\subset\R^n\oplus\R^k$,
$$R_{F_0,\ldots,F_k}^{\Gamma}=\prod_{\Gamma_0,\ldots,\Gamma_k}R_{F^{\Gamma_0}_0,\ldots,F^{\Gamma_k}_k},$$
where the collection $(\Gamma_0,\ldots,\Gamma_k)$ runs over all
collections of faces of the polyhedra $\Delta_0,\ldots,\Delta_k$,
such that $\Gamma,\Gamma_0,\ldots,\Gamma_k$ are compatible.
\end{theor} The proof of a global version of this theorem is given in \cite{ekh}, but
cannot be extended to the local case word by word.

\subsection{Proof of Theorem \ref{elimth}.}\label{sproofelimth}

\textsc{Geometric characterization of eliminant.} We can describe
the geometric meaning of the eliminant $R_{F_0,\ldots,F_k}$ as
follows (we assume that all functions $f_{a,i},\; a\in A_i$, are
holomorphic for simplicity). Denote the convex hull of $A_i$ by
$B_i$, and the dual fan of $B_0+\ldots+B_k$ by $\Sigma$. Then the
product $s_{B_i}\cdot F_i$ extends to a section $\widetilde F_i$
of the line bundle $I_{B_i}$ on the product
$\T^{\Sigma}\times\T^{\tau}$ (we use the notation $I_B$, $s_B$ and
$m(f_1\cdot\ldots\cdot f_k\cdot V)$ introduced in Subsection
\ref{srelbernst}). Denote the projection
$\T^{\Sigma}\times\T^{\tau}\to\T^{\tau}$ by $p$.
\begin{lemma} \label{lgeommeanres} For a germ of a curve
$C\subset\T^{\tau}$ near the origin,
$$m(R_{F_0,\ldots,F_k}\cdot C)=m\Bigl(\widetilde F_0\cdot\ldots\cdot\widetilde F_k\cdot{p^{(-1)}(C)}\Bigr).$$
In particular, both parts of the equality make sense
simultaneously.
\end{lemma}
\textsc{Proof.} We first consider the special case
$R_{F_0,\ldots,F_k}=R_{A_0,\ldots,A_k}$. In this case
$\T^{\tau}=\C[A_0]\oplus\ldots\oplus\C[A_k]$, and the function
$F_i$ equals the \textit{tautological} function $R_i$ on
$\CC^k\times\C[A_0]\oplus\ldots\oplus\C[A_k]$, that maps a point
$(x,\varphi_0,\ldots,\varphi_k)$ to $\varphi_i(x)$. Accordingly,
we denote the section $\widetilde F_i$ by $\widetilde R_i$ in
this case.

1) If $R_{F_0,\ldots,F_k}=R_{A_0,\ldots,A_k}$, and $C$ intersects
the set $R_{A_0,\ldots,A_k}=0$ transversally, then the statement
follows by definition of the $A$-resultant (Definition
\ref{defgkzres}).

2) If $R_{F_0,\ldots,F_k}=R_{A_0,\ldots,A_k}$, and $C$ intersects
the set $R_{A_0,\ldots,A_k}=0$ properly, then we can perturb $C$
so that it intersects the set $R_{A_0,\ldots,A_k}=0$
transversally at a finitely many points, which reduces the
statement to the case (1).

3) In general, consider the map
$\mathcal{F}:\T^{\tau}\to\C[A_0]\oplus\ldots\oplus\C[A_k]$ that
assigns the collection of polynomials
$\bigl(F_0(x,\cdot),\ldots,F_k(x,\cdot)\bigr)$ to every point
$x\in\T^{\tau}$. Accordingly, denote the induced map
$\T^{\Sigma}\times\T^{\tau}\to\T^{\Sigma}\times\C[A_0]\oplus\ldots\oplus\C[A_k]$
by $(\id,\mathcal{F})$. Then we have
$$R_{F_0,\ldots,F_k}=R_{A_0,\ldots,A_k}\circ\mathcal{F}\; \mbox{ and } \;\widetilde F_i=(\id,\mathcal{F})^*\widetilde R_i,$$
which reduces the general case to the case (2). $\quad\Box$

\vspace{0.3cm} \textsc{Proof of Theorem \ref{elimth}.} First, we
can assume without loss of generality that all germs $f_{a,i},\;
a\in A_i$, are holomorphic. Indeed, if we multiply every
$f_{a,i}$ by a monomial $x^{c_i}$, where $c_i\in\tau^{\vee}$ is
far enough from the boundary of $\tau^{\vee}$, then all functions
$f_{a,i}$ become holomorphic, while the statement of Theorem
\ref{elimth} does not change because of the homogeneity of the
$(A_0,\ldots,A_k)$-resultant.

\vspace{0.3cm} \textsc{Proof of Part 1.} For an arbitrary
positive integer linear function $l$ on the cone $\tau^{\vee}$,
let $(l_1,\ldots,l_n)$ be its differential, pick generic complex
numbers $c_1,\ldots,c_n$, and consider the corresponding germ of
a monomial curve $C:\C\to\T^{\tau{\vee}}$ defined by the formula
$C(t)=(c_1t^{l_1},\ldots,c_nt^{l_n})\in\CC^n$ for $t\ne 0$. Then
the minimal value of $l$ on the Newton polyhedron of the
eliminant $R_{F_0,\ldots,F_k}$ equals the intersection number
$m(R_{F_0,\ldots,F_k}\cdot C)$, which, by Lemma
\ref{lgeommeanres}, equals
$$m\Bigl(\widetilde F_0\cdot\ldots\cdot\widetilde F_k\cdot{p^{(-1)}(C)}\Bigr).$$

Denote the convex hull of $A_i$ by $B_i$, the projection
$\R^k\oplus\R^n\to\R^k\oplus\R$ along $\ker l\subset\R^n$ by
$\pi_l$, and the restriction of the germ $F_i$ to the toric
variety $p^{(-1)}(C)$ by $G_i$. If the leading coefficients of
the functions $F_0,\ldots,F_k$ are in general position in the
sense of Definition \ref{genposres}, and the exponents
$l_1,\ldots,l_n$ are generic in the sense that the restriction of
the projection $\pi_l$ to $\Delta_0+\ldots+\Delta_k$ is
one-to-one over bounded faces of its image, then the Newton
polyhedron of $G_i$ equals $\pi_l(\Delta_i)$, and the leading
coefficients of $G_i$ are in general position in the sense of
Definition \ref{genposbernst}. Thus, by Theorem \ref{relbernst},
we have
$$m\Bigl(\widetilde F_0\cdot\ldots\cdot\widetilde
F_k\cdot{p^{(-1)}(C)}\Bigr)=(k+1)!\MV\bigl((\pi_l\Delta_0,B_0\times\R_+),\ldots,(\pi_l\Delta_k,B_k\times\R_+)\bigr).$$

Thus, if $\mathcal{N}$ is the desired Newton polyhedron of the
eliminant $R_{F_0,\ldots,F_k}$, then, for every positive integer
linear function $l$ on the cone $\tau^{\vee}$, we have
$$\min l|_{\mathcal{N}}=(k+1)!\MV\bigl((\pi_l\Delta_0,B_0\times\R_+),\ldots,(\pi_l\Delta_k,B_0\times\R_+)\bigr).$$
This is exactly the formula for the support function of the mixed
fiber polyhedron $\MF(\Delta_0,\ldots,\Delta_k)$, see Proposition
\ref{supportmixfib}.

\begin{rem} Suppose that, on the contrary, the condition of
general position of Definition \ref{genposres} is not satisfied
for some compatible faces $\Gamma_0,\ldots,\Gamma_k$ of the
polyhedra $\Delta_0,\ldots,\Delta_k$. Then the Newton polyhedron
of the eliminant $R_{F_0,\ldots,F_k}$ is strictly smaller than
the mixed fiber polyhedron $\MF(\Delta_0,\ldots,\Delta_k)$.

Namely, pick the face $\Gamma$ of the polyhedron
$\MF(\Delta_0,\ldots,\Delta_k)$, compatible with
$\Gamma_0,\ldots,\Gamma_k$. Consider a linear function $l$ that
attains its minimum on $\Gamma$ as a function on
$\MF(\Delta_0,\ldots,\Delta_k)$. Then, in the notation of the
proof of Part 1, the leading coefficients of the functions
$G_0,\ldots,G_k$ are not in general position in the sense of
Definition \ref{genposbernst}, thus $$\min
l|_{\mathcal{N}}>(k+1)!\MV\bigl((\pi_l\Delta_0,B_0\times\R_+),\ldots,(\pi_l\Delta_k,B_0\times\R_+)\bigr),$$
thus, the face $\Gamma$ is not contained in the Newton polyhedron
of the eliminant $R_{F_0,\ldots,F_k}$.
\end{rem}

\vspace{0.3cm} \textsc{Proof of Part 2.} First, suppose that the
condition of general position of Definition \ref{genposres} is
not satisfied for some faces $\Gamma_0,\ldots,\Gamma_k$ of the
polyhedra $\Delta_0,\ldots,\Delta_k$, compatible with the face
$\Gamma$. Then the corresponding factor in the right hand side of
the desired equality vanishes. On the other hand, by the remark
above, the left hand side vanishes as well.

Suppose that, on the contrary, the condition of general position
of Definition \ref{genposres} is satisfied for all collections of
faces $\Gamma_0,\ldots,\Gamma_k$ of the polyhedra
$\Delta_0,\ldots,\Delta_k$, compatible with the face $\Gamma$.
Then the desired equality is proved in \cite{ekh}. Note that the
proof in \cite{ekh} is written in the global setting, with a
complex torus instead of the toric variety $\T^{\tau}$, and
polynomials instead of analytic functions on it. However, one can
readily verify that the same proof remains valid in the local
setting as well. $\quad\Box$

\subsection{Reach discriminants and their Newton polyhedra.} \label{srdiscrim}

The following version of the discriminant of a projection is not
what we promised to study in the introduction; nevertheless, it
allows us to reduce the study of discriminants of projections to
elimination theory. Let $A$ be a finite set in $\Z^k$, and let
$f_{a}$ be a germ of a meromorphic function on the toric variety
$\T^{\tau}$ for every $a\in A$. Define the germ of a function $F$
on $\T^{\tau}\times\CC^k$ by the formula
$$F(x,y)=\sum_{a\in A} f_{a}(x)y^a \mbox{ for }
x\in\T^{\tau},\; y\in\CC^k,$$ and
denote the number $E_A\bigl(F(x,\cdot)\bigr)$ by $E_F(x)$ for
every $x\in\T^{\tau}$ near the origin (see Subsection
\ref{sadeterm} for the definition of the $A$-determinant $E_A$).
\begin{defin}
\label{defrdiscr} The function $E_{F}$ on the toric variety
$\T^{\tau}$ is called the \textit{reach discriminant of the
projection of the hypersurface $F=0$ to} $\T^{\tau}$.
\end{defin}
We denote the Newton polyhedron of $f_a$ by $\Delta_a$, and
define the Newton polyhedron $\Delta$ of the function $F$ as the
convex hull of the set $\bigcup_{a\in
A}\Delta_{a}\times\{a\}\subset\R^n\oplus\R^k$; then $F(z)$ can be
represented as a power series $\sum_{b\in\Delta}c_{b}z^b$ for
$z\in\CC^n\times\CC^k$. For any $\Gamma\subset\R^n\oplus\R^k$, we
denote the function $\sum_{b\in\Gamma}c_{b}z^b$ by $F^{\Gamma}$.
\begin{defin} \label{genpos}
The leading coefficients of the functions $f_{a},\; a\in A,$ are
said to be \textit{in general position}, if, for every bounded
face $\Gamma\subset\Delta$, such that the restriction of the
projection $\R^n\oplus\R^k\to\R^k$ to $\Gamma$ is injective, $0$
is a regular value of the Laurent polynomial $F^{\Gamma}$.
\end{defin}
Obviously, this condition is satisfied for generic leading
coefficients of the functions $f_{a}$.

\begin{utver} \label{elimthunmix} \label{lprodtrunc} 1) If the leading coefficients of the
functions $f_a,\; a\in A$ are in general position in the sense of
Definition \ref{genpos}, then the Newton polyhedron of $E_F$
equals $\int \Delta$.

2) For every bounded face $\Gamma$ of the polyhedron $\int
\Delta\subset\R^n\subset\R^n\oplus\R^k$,
$$E_{F}^{\Gamma}=\prod E_{F^{\Gamma'}},$$ where $\Gamma'$ runs
over all compatible with $\Gamma$ bounded faces $\Gamma'$ of the
polyhedron $\Delta\subset\R^n\oplus\R^k$, such that the image of
$\Gamma'$ under the projection $\R^n\oplus\R^k\to\R^k$ has the
same dimension as $A$.
\end{utver}
\textsc{Proof.} We can assume that $\dim A=k$ without loss of
generality (otherwise, we can dehomogenize the function $F$).
Consider $k+1$ generic linear combinations of the functions
$F,y_1\frac{\partial F}{\partial y_1},\ldots,y_k\frac{\partial
F}{\partial y_k}$, where $y_1,\ldots,y_k$ are the standard
coordinates on the torus $\CC^k$. We denote these linear
combinations by $F_0,\ldots,F_k$, and note that $\Delta$ is the
Newton polyhedron of each of these functions (while it is not
always the Newton polyhedron of the functions $y_i\frac{\partial
F}{\partial y_i}$). General position for the leading coefficients
of the functions $f_a,\; a\in A$, in the sense of Definition
\ref{genpos} implies general position for the leading
coefficients of the functions $F_0,\ldots,F_k$ in the sense of
Definition \ref{genposres}, thus the statement of Proposition
\ref{elimthunmix} follows from Theorem \ref{elimth} for the
functions $F_0,\ldots,F_k$. $\quad\Box$

\section{Discriminants of hypersurfaces.}

In this section, we study the Newton polyhedron and leading
coefficients of the discriminant of a projection of an analytic
hypersurface, whose Newton polyhedron is given and whose leading
coefficients are in general position.

In the first subsection we give an ``algebraic" definition of the
discriminant, and clarify its geometric meaning in Propositions
\ref{gm11} and \ref{gm12}; the proof of these facts occupies
Subsections \ref{sgenmaps} and \ref{sproofprop1}. In Subsection
\ref{snewtdiscrim}, we study the Newton polyhedron (Theorem
\ref{mainth}) and leading coefficients (Proposition
\ref{thleadcoefdiscr}) of the discriminant. These results are
proved in the last subsection.


\subsection{Discriminants of hypersurfaces.} \label{sdiscrim}

Let $A$ be a finite set in $\Z^k$, let $\tau\subset(\R^n)^*$ be a
convex $n$-dimensional rational polyhedral cone that does not
contain a line, and let $f_{a}$ be a germ of a meromorphic
function on the affine toric variety $\T^{\tau}$ for every $a\in
A$. Define the germ of a function $F$ on $\T^{\tau}\times\CC^k$ by
the formula
$$F(x,y)=\sum_{a\in A} f_{a}(x)y^a \mbox{ for }
x\in\T^{\tau},\; y\in\CC^k,$$ note that $F(x,\cdot)\in\C[A]$, and
denote the number $D_{A}\Bigl(F(x,\cdot)\Bigr)$ by $D_{F}(x)$
(see Subsection \ref{sadiscrim} for the definition of the
$A$-discriminant $D_A$).

\begin{defin} \label{defdiscr} The function $D_{F}$ on the toric variety
$\T^{\tau}$ is called \textit{the discriminant of the projection
of the hypersurface $F=0$ to} $\T^{\tau}$. \end{defin}

The discriminant has the expected geometric meaning if the leading
coefficients are in general position. Namely, denote the set
$\{x\in\CC^n\; |\; D_F(x)=0\}$ by $Z(F)$, and consider the set
$\Sigma(F)$ of all $x\in\CC^n$ such that $0$ is a singular value
of the polynomial $F(x,\cdot)$ on $\CC^k$.
\begin{utver}\label{gm11} Suppose that the Newton polyhedra
of the functions $f_{a},\; a\in A,$ are given, and the leading
coefficients of these functions are in general position. If $A$
is not dual defect (for example, if $A$ satisfies assumptions of
Proposition \ref{pthin}), then
$\overline{\Sigma(F)}=\overline{Z(F)}$ in $\CC^n$, otherwise
$\codim \Sigma(F)>1$.
\end{utver}
This statement can be extended from the maximal torus $\CC^n$ to
the toric variety $\T^{\tau}$ as follows. For a face $\theta$ of
the cone $\tau^{\vee}\subset\R^n$, define $A(\theta)$ as the set
of all $a$ such that the Newton polyhedron of $f_a$ intersects
$\theta$. Consider the set $\Sigma_0(F)$ of all $x\in\T^{\tau}$
such that $0$ is a singular value of the polynomial $F(x,\cdot)$
on $\CC^k$.
\begin{utver} \label{gm12} Suppose that the functions $f_a,\; a\in A,$ are holomorphic,
their Newton polyhedra are such that $A(\theta)\ne\varnothing$ for
every codimension 1 face $\theta\subset\tau^{\vee}$, and their
leading coefficients are in general position. Then \newline 1)
the union of all codimension 1 components of
$\overline{\Sigma_0(F)}$ equals $\overline{Z(F)}$ in $\T^{\tau}$.
\newline 2) If, in addition, $A$ is not dual defect and $\dim
A(\theta)>\dim A+\dim\theta-n$ for every $\theta\ne\tau^{\vee}$,
then $\overline{\Sigma_0(F)}=\overline{Z(F)}$ (in particular,
$\overline{\Sigma_0(F)}$ is a hypersurface).
\end{utver}

Note that $\overline{Z(F)}$ is contained in the zero set of the
discriminant $D_F$ on $\T^{\tau}$, but is smaller in general (even
under the assumptions of the proposition). The equality of
Proposition \ref{gm12}(2) may turn into the strict inequality
$\overline{\Sigma_0(F)}\varsubsetneq\overline{Z(F)}$ in the case
of arbitrary leading coefficients of the functions $f_{a}$, and
even this inequality may be not valid if $\dim A(\theta)<\dim
A+\dim\theta-n$ for some $\theta$. One can readily observe
corresponding examples in the simplest non-trivial case $n=1,\;
A=\{0,1,2\}\subset\Z^1$; a more refined example with $\dim
A(\theta)=\dim A+\dim\theta-n$ and $\Sigma(F)$ not of pure
dimension is given at the end of Subsection \ref{sproofprop1}.

\subsection{Maps with generic leading coefficients.}\label{sgenmaps}

To prove the above statements, we need the following
\begin{utver}
\label{lmian2} Let $h_1,\ldots,h_p$ be either
\newline
1) Laurent polynomials on the complex torus $\CC^n$, or
\newline
2) germs of meromorphic functions on an affine toric variety
$(\T^{\tau},O)$ with no poles in the maximal torus $\CC^n$.

In both cases, consider the map $h=(h_1,\ldots,h_p):\CC^n\to\C^p$.

If $S\subset\C^p$ is an arbitrary algebraic set of codimension
$s$, the Newton polyhedra of the functions $h_i$ are given, and
the leading coefficients of these functions are in general
position, then the set $h^{(-1)}(S)$ has the same codimension $s$.
\end{utver}
Note that, in both settings, $h$ is defined as a map from the
torus $\CC^n$, rather than from the toric variety $(\T^{\tau},O)$,
and, in particular, $h^{(-1)}(S)\subset\CC^n$. If $S'\subset\C^p$
is a constructible set (i.e. if it is obtained by applying the
operations of union, intersection and subtraction to algebraic
sets), then, applying Proposition \ref{lmian2} to its closure
$\overline{S'}$ and to the closure of the difference
$\overline{S'}\setminus S'$, one gets
\begin{sledst}\label{lmian3} If $S'\subset\C^p$ is a constructible set, then,
under the assumptions of Proposition \ref{lmian2}, the closure of
$h^{(-1)}(S')$ equals $h^{(-1)}(\overline{S'})$.
\end{sledst}
To prove Proposition \ref{lmian2}, we reduce it to the following
lemma.
\begin{lemma} \label{lmian} Let $g_1,\ldots,g_p$ be Laurent polynomials on
$\CC^k$, whose coefficients are germs of meromorphic functions on
an affine toric variety $(\T^{\tau},O)$ with no poles in its
maximal torus $\CC^n$. If $Q\subset\CC^n\times\CC^k$ is an
arbitrary algebraic set, the Newton polyhedra of the functions
$g_i$ are given, and the leading coefficients of these functions
are in general position, then the set $\{g_1=\ldots=g_p=0\}$
intersects $Q$ properly near
$\{O\}\times\CC^k\subset\T^{\tau}\times\CC^k$.
\end{lemma}

\begin{defin} For a covector $\gamma\in(\Z^n)^*$ and an analytic function $g(y)=\sum_{a\in\Z^n}c_ay^a$ on $\CC^n$,
the $\gamma$-\textit{truncation} $g^{\gamma}$ is defined to be the
last non-zero sum in the sequence of sums $\sum_{a\; |\;
\gamma(a)=k}c_ay^a,\; k\in\Z$, provided that these sums are equal
to 0 for large $k$.

For a covector $\gamma\in(\Z^n)^*$ and an ideal $I$ in $\C[\Z^n]$,
the $\gamma$-\textit{truncation} $I^{\gamma}$ is the ideal,
generated by $\gamma$-\textit{truncations} of all elements of
$I$. The $\gamma$-\textit{truncation} $Q^{\gamma}$ of an
algebraic variety $Q\in\CC^n$ is defined to be the zero locus of
the $\gamma$-{truncation} of its ideal.
\end{defin}

\textsc{Proof of Lemma \textbf{\ref{lmian}}.} Choose any covector
$\gamma\in(\R^n\times\R^k)^*$ that takes only negative values on
the closed cone $\tau^{\vee}\times\{0\}\subset\R^n\times\R^k$. By
the Bertini-Sard theorem, the set
$\{g^{\gamma}_1=\ldots=g^{\gamma}_p=0\}$ intersects $Q^{\gamma}$
properly under an appropriate assumption of general position for
the leading coefficients of the functions $g_1,\ldots,g_p$. Since
the set of all possible varieties of the form $Q^{\gamma}$ is
finite (see e.g. \cite{sturmf96} or \cite{kaz}), one can choose
the latter assumption of general position to be independent of
$\gamma$. Under this assumption, the set
$\{g^{\gamma}_1=\ldots=g^{\gamma}_p=0\}$ intersects $Q^{\gamma}$
properly for every $\gamma$, thus the same holds for
$\{g_1=\ldots=g_p=0\}$ and $Q$ near the set
$\{O\}\times\CC^k\subset\T^{\tau}\times\CC^k$ (see e.g.
\cite{sturmf96} or \cite{kaz}). $\quad\Box$

In the same way, if $Q$ is smooth, one can prove that
$\{g_1=\ldots=g_p=0\}$ intersects $Q$ transversally, if leading
coefficients are in general position.

\vspace{0.3cm} \textsc{Proof of Proposition
\textbf{\ref{lmian2}}.} Denote the standard coordinates on $\C^p$
by $y_1,\ldots,y_p\,$, consider an arbitrary subdivision
$\{1,\ldots,p\}=I\sqcup J$, and denote the torus
$\{(y_1,\ldots,y_p)\; |\; y_i\ne 0 \mbox{ for } i\in I \mbox{ and
} y_j=0 \mbox{ for } j\in J\}$ by $\CC^I$, and apply the
following Lemma \ref{lmian} to the set
$Q=\bigl(S\cap\CC^I\bigr)\times\T^{\tau}$ and functions
$g_i=\left\{\begin{smallmatrix}h_i-y_i \mbox { \scriptsize for }
i\in I\\h_i \mbox { \scriptsize for } i\in
J\end{smallmatrix}\right.$ on the toric variety
$\CC^I\times\T^{\tau}$. $\quad\Box$

\subsection{Proof of Propositions \ref{gm11} and \ref{gm12}.}\label{sproofprop1}

\textsc{Proof of Proposition \textbf{\ref{gm11}}.} Recall that
$\Sigma\subset\C[A]$ is the set of all $\varphi$ such that
$\varphi(y)={\rm d}\varphi(y)=0$ for some $y\in\CC^k$. Its
closure $\overline{\Sigma}$ is a hypersurface and is defined by
the equation $D_A=0$ (otherwise identically $D_A=1$ and
$Z(F)=\varnothing$ by definition).

Let $\mathcal{F}:\T^{\tau}\to\C[A]$ be the map that assigns the
polynomial $F(x,\cdot)$ to a point $x\in\T^{\tau}$. We can
express the desired sets $\Sigma(F)$ and $Z(F)=\{D_F=0\}$ in
terms of this map:
$$\Sigma(F)=\mathcal{F}^{(-1)}(\Sigma)\cap\CC^n,$$
$$D_F=D_A\circ\mathcal{F}.$$ Proposition \ref{gm11} now follows by Corollary \ref{lmian2}
with $(h_1,\ldots,h_p)=\mathcal{F}$ and $S'=\Sigma$. $\quad\Box$

We can also formulate a refinement of Proposition \ref{gm11} with
multiplicities taken into account (the proof follows the same
lines but requires more technical details; we omit it, since we
do not need this refinement in what follows). Restrict the
projection $p:\T^{\tau}\times\CC^k\to \T^{\tau}$ to the regular
locus of the set $\{F=0\}$ and denote the singular locus of this
restriction by $S$.
\begin{utver}(See Subsection \ref{smultci} for the
notation.)\newline $p_*(S)$ equals $[D_{\overline F}=0]$ on the
complex torus $\CC^n$, and does not contain codimension 1 orbits
of the toric variety $\T^{\tau}$, if the leading coefficients of
the functions $f_{a},\; a\in A, $ are in general position.
\end{utver}

\vspace{0.3cm} \textsc{Proof of Proposition \textbf{\ref{gm12}}.}
For a closed regular subvariety $M\subset\CC^k$, define
$\Sigma(M)\subset\C[A]$ as the set of all Laurent polynomials
$\varphi$ in $\C[A]$, such that the set $\{y\in M\; |\;
\varphi(y)={\rm d}\varphi(y)=0\}$ has at least one isolated
point. Such definition implies the following properties of
$\Sigma(M)$ (in contrast to $\Sigma$):
\newline
1) for every $\varphi\in\Sigma(M)$, at least one of the local
components of $\Sigma(M)$ near $\varphi$ is closed.
\newline
2) The set $\Sigma(M)$ has codimension $1+\codim M+\dim A-k$ at
all of its points. We need the following corollary of (1) and (2):
\newline 3) The set
$\mathcal{F}^{(-1)}\bigl(\Sigma(M)\bigr)$ has codimension at most
$1+\codim M+\dim A-k$ at all of its points. We denote the latter
set by $\Sigma(M,F)$.

For every face $\theta\subset\tau^{\vee}$, we can choose a closed
regular subvariety $M_{\theta}\subset\CC^k$ such that
$\Sigma(M_{\theta})\cap\C[A(\theta)]$ is dense in
$\Sigma\cap\C[A(\theta)]$. Let $T_{\theta}$ be the orbit of the
variety $\T^{\tau}$, corresponding to the face $\theta$, then the
restriction of the map $\mathcal{F}$ to this orbit is a map
$\mathcal{F}_{\theta}:T_{\theta}\to\C[A(\theta)]$. Proposition
\ref{lmian2} for $(h_1,\ldots,h_p)=\mathcal{F}_{\theta},\;
S=\overline{\Sigma(M)}_{\theta}\cap\C[A(\theta)]$ and Corollary
\ref{lmian3} for $(h_1,\ldots,h_p)=\mathcal{F}_{\theta},\;
S={\Sigma(M)}_{\theta}\cap\C[A(\theta)]$ imply the following:
\newline 4) If the leading coefficients of the functions $f_a,\;a\in A,$
are in general position, then the set $\Sigma(M,F)\cap
T_{\theta}$ is dense in $\Sigma_0(F)\cap T_{\theta}$, and its
codimension in $T_{\theta}$ is equal to $1+\codim M+\dim
A(\theta)-k$, which is greater than $1+\codim M+\dim
A-k+\dim\theta-n$ under the assumption of Proposition \ref{gm12}.

Since, by (3) and (4), the codimension of $\Sigma(M,F)\cap
T_{\theta}$ in the toric variety $\T^{\tau}$ is greater than the
codimension $\Sigma(M,F)$ in $\T^{\tau}$ at every point of
$\Sigma(M,F)\cap T_{\theta}$, then $\Sigma(M,F)\cap T_{\theta}$
is contained in the closure of $\Sigma(M,F)\setminus T_{\theta}$.
The inclusions
$$\overline{\Sigma_0(F)\cap T_{\theta}}=\overline{\Sigma(M,F)\cap
T_{\theta}}\subset\overline{\Sigma(M,F)\setminus
T_{\theta}}\subset\overline{\Sigma(F)}\subset\overline{Z(F)}$$
prove Proposition \ref{gm12}. $\quad\Box$

\begin{exa} Note that the inclusion $\overline{\Sigma_0(F)\cap
T_{\theta}}\subset\overline{\Sigma(F)}$ and the statement of
Proposition \ref{gm12} may fail, if $\dim A(\theta)=\dim
A+\dim\theta-n$. For example, let $\T^{\tau}$ be the space $\C^3$
with coordinates $x,y,z$, and let $T_{\theta}$ be the the torus
$\{x\ne 0,\; y\ne 0,\; z=0\}$. Pick generic linear functions
$l_1,l_2,l_3,m_1,m_2$ of the variables $x$ and $y$, and choose
the face $A$ and the functions $f_a,\; a\in A$ as follows (each
function $f_a(x,y,z)$ is written near the corresponding point
$a\in A$):
\begin{center}
\noindent\includegraphics[width=7cm]{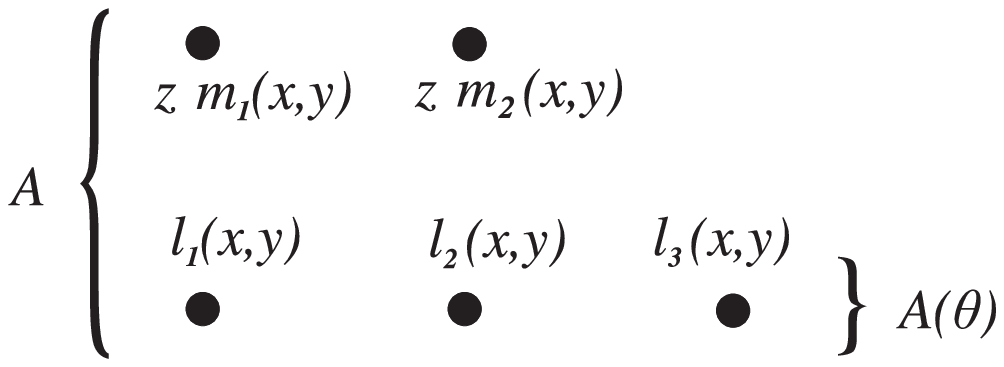} \end{center} Then
the sets $\overline{\Sigma_0(F)\cap T_{\theta}}$ and
$\overline{\Sigma(F)}\cap T_{\theta}$ are given by the equations
$$l_1^2-l_1\, l_3=0\quad \mbox{ and }\quad l_1\, m_1^2-l_2\, m_1\, m_2+l_3\, m_2^2=0$$
in $x$ and $y$, and hence do not intersect.
\end{exa}

\subsection{Newton polyhedra of
discriminants of hypersurfaces.} \label{snewtdiscrim} If the
Newton polyhedra $\Delta_a$ of functions $f_{a},\, a\in
A\subset\Z^k$, on the affine toric variety $\T^{\tau}$ are given,
and the leading coefficients of these functions satisfy a certain
condition of general position, then the Newton polyhedron of the
discriminant $D_{F}$, where $F(x,y)=\sum_{a\in A} f_{a}(x)y^a$,
depends only on the Newton polyhedra of these functions, not on
the coefficients. Theorem \ref{mainth} and Proposition
\ref{thleadcoefdiscr} below solve the following problem:

\begin{quote}
Express the Newton polyhedron and leading coefficients of the
discriminant $D_{F}$ in terms of the Newton polyhedra and leading
coefficients of the functions $f_{a}$, provided that leading
coefficients are in general position.
\end{quote}

\textsc{Newton polyhedron of the discriminant.} The Minkowski
integral $\int\Delta\subset\R^n$ of a polyhedron
$\Delta\subset\R^n\oplus\R^k$ and combinatorial Euler obstructions
$e^{A',A}$ and Milnor numbers $c^{A',A}$ are introduced in Section
1. Recall that a face of a set $A\subset\Z^k$ is the intersection
of $A$ with a face of its convex hull, and $\dim A$ is the
dimension of its convex hull. For a face $A'$ of the set
$A\subset\Z^k$, we denote the convex hull of the union
$\bigcup_{a\in A'}\Delta_{a}\times\{a\}$ by $\Delta(A')$; this is
an unbounded face of $\Delta=\Delta(A)$.


\begin{theor} \label{mainth}
1) If the Newton polyhedra of the functions $f_a$ are given, and
the leading coefficients of these functions are in general
position in the sense of Definition \ref{genpos}, then the Newton
polyhedron of the discriminant $D_F$ equals
$$\mathcal{N}^A_{\Delta}=\sum_{A'\subset A}e^{A',\, A}\cdot \int \Delta(A'),$$
where $A'$ runs over all faces of $A$, including $A'=A$.

2) If the leading coefficients of the functions $f_a$ are
arbitrary, then the Newton polyhedron of the discriminant $D_F$
is contained in $\mathcal{N}^A_{\Delta}$.
\end{theor}

The proof is given in Subsection \ref{sproofth1}, an example of
application is given in Subsection \ref{losungex}. Note that
Definition \ref{genpos} is not the weakest possible condition of
general position for this theorem. The weakest one (very
complicated to verify though) can be extracted from Proposition
\ref{thleadcoefdiscr} below, see the subsequent discussion.

\vspace{0.3cm} \textsc{Leading coefficients of the discriminant.}
Let $\Gamma$ be a face of the polyhedron $\Delta$, denote the set
of all $a\in A$, such that $\Gamma$ intersects $\Delta_a\times
a$, by $A_{\Gamma}$, and denote the minimal face of $A$,
containing $A_{\Gamma}$, by $\bar A_{\Gamma}$. For a bounded face
$\widetilde\Gamma$ of another polyhedron, define the number
$d^{\widetilde\Gamma,\Gamma}$ as
$\sum_{\Gamma'}c^{A_{\Gamma},A_{\Gamma'}}\cdot e^{\bar
A_{\Gamma'},A}$, where $\Gamma'$ runs over all bounded compatible
with $\widetilde\Gamma$ faces of $\Delta$ such that
$\Gamma'\supset\Gamma$ and $\dim A_{\Gamma'}=\dim \bar
A_{\Gamma'}$. Defining the value $D_{F^{\Gamma}}(x)$ as
$R_{A_{\Gamma}}\Bigl(F^{\Gamma}(x,\cdot)\Bigr)$ for $x\in\CC^n$,
the Laurent polynomial $D_{F^{\Gamma}}$ on $\CC^n$ depends only on
leading coefficients of the functions $f_a,\, a\in A_{\Gamma}$.
\begin{utver} \label{thleadcoefdiscr} 1) For every bounded face $\widetilde\Gamma$ of the expected Newton polyhedron
$\mathcal{N}^A_{\Delta}$ of the discriminant $D_F$,
$$D^{\widetilde\Gamma}_{F}=\prod_{\Gamma}( D_{ F^{\Gamma}
})^{d^{\widetilde\Gamma,\, \Gamma}},$$ where $\Gamma$ runs over
all bounded faces of the polyhedron $\Delta$.

2) Every factor $( D_{ F^{\Gamma} })^{d^{\widetilde\Gamma,\,
\Gamma}}$ in the right hand side of this equality is a polynomial
(i.e. $d^{\widetilde\Gamma,\, \Gamma}\geqslant 0$ whenever the
polynomial  $D_{ F^{\Gamma} }$ is of positive degree).
\end{utver}
In particular, the leading coefficients of the discriminant $D_F$
only depend on those of the functions $f_a,\; a\in A$. The proof
is given in Subsection \ref{sproofth2}.

An assumption of general position for leading coefficients would
be redundant in this statement: if the leading coefficients of the
functions $f_a,\; a\in A,$ are degenerate enough, then both parts
of the equality become identically zero simultaneously. In
particular, the Newton polyhedron of $D_F$ equals
$\mathcal{N}^A_{\Delta}$ (i.e. $D_F$ has non-zero coefficients of
the monomials, corresponding to the vertices of
$\mathcal{N}^A_{\Delta}$) if and only if the following condition
is satisfied for all $\widetilde\Gamma$ and $\Gamma$: if
$\widetilde\Gamma$ is a vertex, and $d^{\widetilde\Gamma,\,
\Gamma}>0$, then $D_{ F^{\Gamma} }$ is not identically zero. Note
that the coefficient $d^{\widetilde\Gamma,\, \Gamma}$ is
complicated to compute, and no simple combinatorial criterion for
its positivity (i.e. for divisibility $D^{\widetilde\Gamma}_{F}$
by $D_{ F^{\Gamma} }$) is known. See, for example, Theorem 15 in
\cite{cat} for one important special case.

\vspace{0.3cm} \textsc{Degree of $A$-discriminants.} 
The Gelfand-Kapranov-Zelevinsky discriminant $D_A$ is a special
case of the discriminant $D_F$ with $\T^{\tau}=\C[A]$. In this
case, the discriminant is homogeneous, and its Newton polyhedron
$\mathcal{N}$ is contained in the space $\R^A$, whose integer
lattice consists of monomials of the form $\prod_{a\in A}
c_a^{\lambda_a}$ in coefficients of the indeterminate polynomial
$\sum_{a\in A}c_ax^a\in\C[A]$. We consider $\lambda_a,\, a\in A$,
as a system of coordinates on $\R^A$, denote $\sum_a \lambda_a$
by $l$, and note that the degree of the discriminant $D_A$ is the
minimal value of $l$ on $\mathcal{N}$.

Computing $\mathcal{N}$ by Theorem \ref{mainth} and then $\min
l|_{\mathcal{N}}$ by Proposition \ref{supportmixfib} in this
case, we get the following formula for the degree of $D_A$. Let
$A'$ be a face of $A$, and let $M\subset\R^k$ be the vector space,
parallel to the affine span of $A'$. We choose the volume form
$\mu$ on $M$ such that the volume of $M/(M\cap\Z^k)$ equals
$(\dim M)!$, and denote the $\mu$-volume of the convex hull of
$A'$ by $\Vol A'$.
\begin{defin}\label{defdegree}
Define the number $\deg A$ as the sum $$\sum\limits_{A'\subset
A}e^{A',\, A}\cdot (\dim A'+1)\Vol(A')$$ over all faces $A'$ of
the set $A$, including $A'=A$.
\end{defin}
\begin{sledst}[\cite{tak}]\label{ldegree} 1) $\deg D_A=\deg A$. \newline 2) $D_A$
is a constant if and only if $\,\deg A=0$.
\end{sledst}
\begin{rem} A useful generalization of this fact for discriminants of
higher codimension is proved in \cite{tak}, based on a totally
different technique (the Ernstr\"om formula \cite{ernst}).
Amazingly, that technique also ends up with Euler obstructions of
toric varieties, which suggests that the two techniques could be
unified. In particular, it would be interesting to find a common
generalization of Theorem \ref{mainth} and the Ernstr\"om
formula, which would, for instance, compute the tropicalization
of the dual of an arbitrary projective variety $V$ in terms of
Euler obstructions of truncations of $V$.\end{rem}

\subsection{Proof of Theorem \ref{mainth} and Proposition \ref{thleadcoefdiscr}.} \label{sproofth1} \label{sproofth2}

Applying the Gelfand-Kapranov-Zelevinsky decomposition
(Proposition \ref{gkzdecomp}) to the definition of reach
discriminant (Definition \ref{defrdiscr}), we have the following
relation between discriminants and reach discriminants (see
Definition \ref{defrdiscr}) for every face $A'\subset A$:
$$E_{F^{\Delta(A')}}=\prod_{A''}( D_{ F^{\Delta(A'')} })^{c^{A'',\, A'}},\eqno{(*_{A'})}$$
where $A''$ runs over all faces of $A'$, including $A''=A'$.
Inverting the formulas $(*_{A'})$ by induction on the dimension
of $A'$, we have the following relations:
$$D_{F^{\Delta(A')}}=\prod_{A''}( E_{ F^{\Delta(A'')} })^{e^{A'',\, A'}},\eqno{(**_{A'})}$$
where $A''$ runs over all faces of $A'$, including $A''=A'$.

In more detail, assume that we have already obtained the formulas
$(**_{B})$ for all faces $B$ of dimension less than $p$. Then, for
every $A'$ of dimension $p$, rewriting the formula $(*_{A'})$ as
$D_{F^{\Delta(A')}}=E_{F^{\Delta(A')}}\cdot\prod_{A''\ne A'}(
D_{F^{\Delta(A'')} })^{-c^{A'',\, A'}}$, expressing
$D_{F^{\Delta(A'')}}$ in terms of $E_{F^{\Delta(A'')}}$ by the
formulas $(**_{B})$ in the right hand side, and collecting similar
multipliers, we obtain the formula $(**_{A'})$.

Informally speaking, if we consider the formal logarithm of the
formulas $(*_{A'})$ to pass to the additive notation instead of
the multiplicative one, then the vector of logarithms $\ln
E_{F^{\Delta(A')}},\, A'\subset A$, equals the matrix with
entries $c^{A'',\, A'}$ times the vector of logarithms $\ln
D_{F^{\Delta(A')}},\, A'\subset A$. Inverting the matrix, we
obtain the logarithm of the formulas $(**_{A'})$.

In particular, $D_A(\varphi)= \prod_{A'}
\Bigl(E_{A'}(\varphi^{A'})\Bigr)^{e^{A',\, A}}$.

\vspace{0.3cm} \textsc{Proof of Theorem} \textbf{\ref{mainth}}.
To prove Part 1, apply Proposition \ref{elimthunmix}(1) and the
identity $(**_A)$ above.

To prove Part 2, consider a monomial curve $C:\C\to\T^{\tau}$,
corresponding to an arbitrary positive integer linear function
$\gamma:\tau^{\vee}\to\R_{>0}$; by definition,
$C(t)=(h_1t^{\gamma_1},\ldots,h_nt^{\gamma_n})\in\CC^n$, where
$\gamma_1,\ldots,\gamma_n$ are the coefficients of the linear
function $\gamma$, and the coefficients $h_1,\ldots,h_n$ are
generic. Then $D_{F}\bigl(C(t)\bigr)$ is a germ of a meromorphic
function of one variable $t$, and the order of its zero or pole
is equal to the minimal value of the function $\gamma$ on the
Newton polyhedron of $D_F$. Since this order of zero depends
upper-semicontinuously on $F$, so does the Newton polyhedron of
$D_F$. $\quad\Box$


\vspace{0.3cm} \textsc{Proof of Proposition
\textbf{\ref{thleadcoefdiscr}}, Part \textbf{1}.} For every face
$A'$ of the set $A$, choose a face $\widetilde\Gamma(A')$ of the
polyhedron $\int\Delta(A')$, such that these faces
$\widetilde\Gamma(A')$ together with $\widetilde\Gamma$ form a
compatible collection. Represent the discriminant $D_F$ as a
product of reach discriminants by formula $(**_A)$ above, then
$$D^{\widetilde\Gamma}_{F}=\prod_{A'}( E^{\widetilde\Gamma(A')}_{ F^{\Delta(A')} })^{e^{A',\, A}}.$$
By Proposition \ref{elimthunmix}(2), represent every truncated
reach discriminant $E^{\widetilde\Gamma(A')}_{ F^{\Delta(A')} }$
in the right hand side as a product of reach discriminants of the
form $E_{F^{\Gamma'}},\; \Gamma'\subset\Delta(A')$. Finally,
represent each of these reach discriminants as a product of
discriminants by formula $(*_{A_{\Gamma'}})$ above. $\quad\Box$

\vspace{0.3cm} \textsc{Proof of Proposition
\textbf{\ref{thleadcoefdiscr}}, Part \textbf{2}.} If the set
$A_{\Gamma}$ is dual defect, then identically $D_{ F^{\Gamma}
}=1$, and the sign of the exponent $d^{\widetilde\Gamma,\,
\Gamma}$ in the right hand side of the identity in the statement
of Proposition \ref{thleadcoefdiscr}(1) is not important.
Otherwise, we have
\begin{lemma} If the set $A_{\Gamma}$ is not dual defect, then $d^{\widetilde\Gamma,\, \Gamma}\geqslant 0$ for
every bounded face
$\widetilde\Gamma\subset\mathcal{N}^A_{\Delta}$.
\end{lemma}
\textsc{Proof.} Let $A_{\Gamma}$ consist of points
$a_1,\ldots,a_N$, and let $D_{A_{\Gamma}}(u_1,\ldots,u_N)$ be the
value of the discriminant $D_{A_{\Gamma}}$ at the polynomial
$\sum_i u_i t^{a_i}\in\C[A_{\Gamma}]$. Consider both sides of the
identity in the statement of Proposition \ref{thleadcoefdiscr}(1)
as polynomials of the leading coefficients of the functions
$f_a,\; a\in A$, with a fixed value of the variable $x\in\CC^n$.
Then the discriminant $D_{ F^{\Gamma} }$ equals
$D_{A_{\Gamma}}(l_1,\ldots,l_N)$, where $l_i$ is a non-zero
linear combination of leading coefficients of the function
$f_{a_i}$. Since $D_{A_{\Gamma}}$ is a power of a homogeneous
irreducible polynomial that non-trivially depends on all the
variables $u_1,\ldots,u_N$ (by Lemmas \ref{lirred} and
\ref{lnontriv}), so does $D_{ F^{\Gamma}
}=D_{A_{\Gamma}}(l_1,\ldots,l_N)$: it is a power of a homogeneous
irreducible polynomial that non-trivially depends on all the
coefficients of $F^{\Gamma}$ and does not depend on other leading
coefficients of the functions $f_a,\, a\in A$.

Thus, if $d^{\widetilde\Gamma,\, \Gamma}$ were negative for some
$\Gamma$, then the factor $( D_{ F^{\Gamma}
})^{d^{\widetilde\Gamma,\, \Gamma}}$ could not be cancelled by
other multipliers in the right hand side of the equality in the
statement of Proposition \ref{thleadcoefdiscr}(1), and the right
hand side were a rational function with a non-trivial
denominator. But this is impossible because the left hand side is
a polynomial. $\quad\Box$

\section{Discriminants of complete intersections.}

In this section, we generalize the results of the previous two
sections to discriminants of projections of analytic complete
intersections.

The discriminant is defined in the first subsection. In Subsection
\ref{snewtdiscrci}, we reduce the study of its Newton polyhedron
and leading coefficients to the case of projections of
hypersurfaces, which is studied in the previous section. In some
important special cases (Subsections \ref{snewtdiscrci2} and
\ref{snewtdiscrci3}), this leads to an explicit answer. An
example of computation of such answer is given in the last
subsection. We also consider an alternative definition of the
discriminant in Subsection \ref{sfibration}.


\subsection{Discriminants of complete intersections.}\label{sdiscrimci}
Let $\tau\subset(\R^n)^*$ be a convex $n$-dimensional rational
polyhedral cone that does not contain a line, and let
$A_0,\ldots,A_l,\; l\leqslant k,$ be finite sets in $\Z^k$. For
every $i=0,\ldots,l,\; a\in A_i,\;$ let $f_{a,i}$ be a germ of a
meromorphic function on the toric variety $\T^{\tau}$ with no
poles in the maximal torus. Define the germ of a function $F_i$
on $\T^{\tau}\times\CC^k$ by the formula
$$F_i(x,y)=\sum_{a\in A_i} f_{a,i}(x)y^a \mbox{ for }
x\in\T^{\tau},\; y\in\CC^k,$$ and denote the number
$D^{red}_{A_0,\ldots,A_l}\Bigl(F_0(x,\cdot),\ldots,F_l(x,\cdot)\Bigr)$
by $D^{red}_{F_0,\ldots,F_l}(x)$ (see Subsection
\ref{sgkzdiscrimci} for the definition of the discriminant
$D^{red}_{A_0,\ldots,A_l}$).

\begin{defin} The germ of the function $D^{red}_{F_0,\ldots,F_l}$ on the affine toric variety
$(\T^{\tau},O)$ is called \textit{the discriminant of the
projection of the complete intersection $F_0=\ldots=F_l=0$ to}
$\T^{\tau}$.
\end{defin}

The discriminant has the expected geometric meaning if the leading
coefficients are in general position. Namely, denote the set
$\{x\in\CC^n\; |\; D^{red}_{F_0,\ldots,F_l}(x)=0\}$ by
$Z(F_0,\ldots,F_l)$, and consider the set
$\Sigma(F_0,\ldots,F_l)$ of all $x\in\CC^n$, such that
$(0,\ldots,0)$ is a singular value of the map
$\Bigl(F_0(x,\cdot),\ldots,F_l(x,\cdot)\Bigr):\CC^k\to\C^{l+1}$.
\begin{utver} \label{gm21} Suppose that the Newton polyhedra
of the functions $f_{a,i},\; a\in A_i,$ are given, and the leading
coefficients of these functions are in general position. Then
\newline 1) the union of codimension 1 components of the closure
$\overline{\Sigma(F_0,\ldots,F_l)}$ equals
$\overline{Z(F_0,\ldots,F_l)}$. \newline 2) If, in addition, the
collection $A_0,\ldots,A_k$ is not dual defect (for example, if it
satisfies assumptions of Proposition \ref{equidim}), then
$$\overline{\Sigma(F_0,\ldots,F_l)}=\overline{Z(F_0,\ldots,F_l)},$$
and, in particular, $\Sigma(F_0,\ldots,F_l)$ is a hypersurface.
\end{utver}
This can be extended from the maximal torus $\CC^n$ to the toric
variety $\T^{\tau}$ in the same way as Proposition \ref{gm11}
(see Propositon \ref{gm12} for the notation). Consider the set
$\Sigma_0(F_0,\ldots,F_l)$ of all $x\in\T^{\tau}$, such that
$(0,\ldots,0)$ is a singular value of the map
$\Bigl(F_0(x,\cdot),\ldots,F_l(x,\cdot)\Bigr):\CC^k\to\C^{l+1}$.
\begin{utver}\label{gm22} Suppose that the functions $f_{a,i},\; a\in A_i,$ are holomorphic,
their Newton polyhedra are such that $A_j(\theta)\ne\varnothing$
for every codimension 1 face $\theta\subset\tau^{\vee}$ and
$j=0,\ldots,l$, and their leading coefficients are in general
position. Then \newline 1) the union of all codimension 1
components of $\overline{\Sigma_0(F_0,\ldots,F_l)}$ equals
$\overline{Z(F_0,\ldots,F_l)}$.
\newline 2) If, in addition,
the collection $A_0,\ldots,A_k$ is not dual defect, and $\dim
\sum_j A_j(\theta)>\dim L+\dim\theta-n$ for every
$\theta\ne\tau^{\vee}$, then
$$\overline{\Sigma_0(F_0,\ldots,F_l)}=\overline{Z(F_0,\ldots,F_l)},$$
and, in particular, $\Sigma_0(F_0,\ldots,F_l)$ is a hypersurface.
\end{utver}
\noindent Since the proof of these facts is the same as for
Propositions \ref{gm11} and \ref{gm12}, with the exception of
more complicated notation coming from $l>0$, we omit it.

\subsection{Newton polyhedra of discriminants of complete intersections.}\label{snewtdiscrci}

The study of the Newton polyhedron and leading coefficients of
the discriminant $D^{red}_{F_0,\ldots,F_l}$ can be reduced to the
case $l=0$ (which is studied in the previous section) by the
Cayley trick, which represents $D^{red}_{F_0,\ldots,F_l}$ as a
product of discriminants of the form $D^{red}_{G_J}$ for linear
combinations $G_J=\sum_{j\in J}\lambda_jF_j$ with indeterminate
coefficients $\lambda_j$, where $J$ runs over certain subsets of
$\{0,\ldots,l\}$.

More precisely, define the function $G_J$ on
$\T^{\tau}\times\CC^k\times\CC^{l+1}$ as $\sum_{j\in
J}\lambda_jF_j$, where $\lambda_0,\ldots,\lambda_l$ are
coordinates on $\CC^{l+1}$. Let $e_0,\ldots,e_l$ be the standard
basis in $\Z^{l+1}$, denote the set $\bigcup\limits_{j\in J}
A_j\times\{e_j\}$ by $A_J\subset\Z^k\oplus\Z^{l+1}$. Then, for
every $x\in\T^{\tau}$, the polynomial $G_J(x,\cdot)$ on
$\CC^k\times\CC^{l+1}$ is contained in $\C[A_J]$, and the
discriminant $D^{red}_{G_J}$ is defined by the formula
$D^{red}_{G_J}(x)=D^{red}_{A_J}\Bigl(G_J(x,\cdot)\Bigr)$.

These discriminants are related to the desired one as follows.
For every $J\subset\{0,\ldots,l\}$, denote the difference
$\dim\sum_{j\in J}A_j-|J|$ by $\codim J$.

\begin{theor}[Cayley trick]\label{cayley} The discriminant $D^{red}_{F_0,\ldots,F_l},\; 0<l<k,$ equals the product of the discriminants
$D^{red}_{G_J}$ over all subsets $J\subset\{0,\ldots,l\}$, such
that $\codim J\leqslant\codim J'$ for every $J'\supset J$.
\end{theor}
This is Theorem \ref{univcayley} in the new notation. To describe
the Newton polyhedra and leading coefficients of the discriminants
$D^{red}_{G_J}$, and therefore those of
$D^{red}_{F_0,\ldots,F_l}$, we can apply Theorem \ref{mainth} and
Proposition \ref{thleadcoefdiscr} to the functions $G_J$ under an
appropriate condition of general position for their leading
coefficients (see Definition \ref{genpos}). In many cases, the
result of this computation can be written as an explicit formula
for the Newton polyhedron of $D^{red}_{F_0,\ldots,F_l}$; see, for
example, Theorem \ref{losung} below. The aforementioned condition
of general position can be formulated in terms of leading
coefficients of the functions $F_0,\ldots,F_l$ as follows.

We denote the Newton polyhedron of $f_{a,i}$ by $\Delta_{a,i}$,
and define the Newton polyhedron $\Delta_i$ of the function $F_i$
as the convex hull of the set $\bigcup_{a\in
A_i}\Delta_{a,i}\times\{a\}\subset\R^n\oplus\R^k$; then $F_i(z)$
can be represented as a power series
$\sum_{b\in\Delta_i}c_{b,i}z^b$ for $z\in\CC^n\times\CC^k$. If
$\Gamma$ is a face of $\Delta_i$, then we denote the function
$\sum_{b\in\Gamma}c_{b,i}z^b$ by $F_i^{\Gamma}$.
\begin{defin} \label{genposci}
The leading coefficients of the functions $f_{a,i},\; a\in A_i,$
are said to be \textit{in general position}, if, for every
sequence $0\leqslant i_1<\ldots<i_q\leqslant l$ and every
collection of compatible bounded faces
$\Gamma_{i_j}\subset\Delta_{i_j}$ (see Definition
\ref{defcompat}), such that the restriction of the projection
$\R^n\oplus\R^k\to\R^k$ to $\Gamma_{i_1}+\ldots+\Gamma_{i_q}$ is
injective, $(0,\ldots,0)$ is a regular value of the polynomial map
$(F_{i_1}^{\Gamma_{i_1}},\ldots,F_{i_q}^{\Gamma_{i_q}}):\CC^{n+k}\to\C^q$.
\end{defin}

\subsection{The case of analogous Newton polyhedra.}\label{snewtdiscrci2}
Under some additional assumptions on the sets $A_0,\ldots,A_l$,
the Cayley trick allows to explicitly compute the Newton
polyhedron of the discriminant of a complete intersection as
follows.
\begin{defin} Let $A'$ be a face of a
finite set $A\subset\R^k$, and let $L$ be the vector subspace in
$\R^k$, parallel to the affine span of $A'$. The
$A'$-\textit{link of} $A$ is a (non-convex) polyhedron $\widetilde
A\setminus\widetilde A'\subset\R^k/L$, where $\widetilde A$ and
$\widetilde A'$ are the convex hulls of the images of the sets 
$A$ and $A\setminus A'$ under the projection $\R^k\to\R^k/L$.
\end{defin}
\begin{defin} \label{defan} Finite sets $A$ and $B$ in $\R^k$ are said to be
\textit{analogous}, if there is a one-to-one correspondence
between the posets of their faces, such that, for every pair of
corresponding faces $A'\subset A$ and $B'\subset B$, the
$A'$-link of $A$ equals the $B'$-link of $B$ up to a parallel
translation (in particular, the affine spans of $A'$ and $B'$ are
parallel to the same subspace $L\subset\R^k$).
\end{defin}
\begin{exa} 1) If $A=B$, then $A$ and $B$ are analogous. \newline
2) If $P$ and $Q$ are analogous integer polyhedra (i.e. their dual
fans coincide) and $k\in\Z$ is large enough, then the sets of
integer points in $kP$ and $kQ$ are analogous. Note that those
sets are not necessary analogous for $k=1$. For example, the two
sets on the picture in Subsection \ref{seulobstr} have different
links of their vertical faces.
\end{exa}

Recall that the standard basis in $\Z^{l+1}$ is denoted by
$e_0,\ldots,e_l$.
\begin{lemma} \label{euljoin} If $A_0,\ldots,A_l$ in $\Z^k$ are analogous, then, for every
collection of corresponding faces $A'_0\subset
A_0,\ldots,A'_l\subset A_l$, \newline 1)
$e^{A'_0,A_0}=\ldots=e^{A'_l,A_l}$, \newline 2)
$A'=\bigcup_{i=0}^{l'} A'_i\times\{e_i\}$ is a face of
$A=\bigcup_{i=0}^l A_i\times\{e_i\}$, and $e^{A',A}=e^{A'_0,A_0}$.
\end{lemma}
\textsc{Proof.} Part 1 and Part 2 for $l'=l$ follow by the fact,
that the Euler obstruction $e^{B',B}$ depends on the $B'$-link of
$B$ only (by definition, see Subsection \ref{seulobstr}). Since
$e^{B',B}$ is the Euler obstruction of the $B$-toric variety at a
point of its $B'$-orbit (Theorem \ref{thmt}), and since Euler
obstruction is a local topological invariant, then $e^{A',A}$
does not depend on $l'$, and it is enough to prove Part 2 for
$l'=l$. $\quad\Box$

In the notation of Subsection \ref{sdiscrimci}, let $M$ be the
lattice, generated by pairwise differences of points of the set
$A_0+\ldots+A_l$, and let $\Delta_0,\ldots,\Delta_l$ be the
Newton polyhedra of the functions $F_0,\ldots,F_l$. Recall that
we denote the mixed fiber polyhedron of polyhedra
$P_1,\ldots,P_q$ by the monomial $P_1\cdot\ldots\cdot P_q$.
\begin{theor}\label{losung} 1) If the leading coefficients of the
functions $f_{a,i},\; a\in A_i,$ are in general position in the
sense of Definition \ref{genposci}, the sets $A_0,\ldots,A_l$ are
analogous and not contained in an affine hyperplane, then the
Newton polyhedron of the discriminant $D^{red}_{F_0,\ldots,F_l}$
equals
$$\mathcal{N}\; =\; \frac{1}{|\Z^k/M|}\sum_{A'_0,\ldots,A'_l} \qquad e^{A'_0,A_0}
\sum_{a_0>0,\ldots,a_l>0\atop a_0+\ldots+a_l=\dim A'_0+1}
\Delta_0(A'_0)^{a_0}\cdot\ldots\cdot\Delta_l(A'_l)^{a_l},$$ where
the collection $(A'_0,\ldots,A'_l)$ runs over all collections of
corresponding faces $A'_0\subset A_0,\ldots,A'_l\subset A_l$,
including $A'_0=A_0,\ldots,A'_l=A_l$. \newline 2) If the leading
coefficients of the functions $f_{a,i}$ are arbitrary, then the
Newton polyhedron of the discriminant $D^{red}_{F_0,\ldots,F_l}$
is contained in $\mathcal{N}$.
\end{theor}
An example of application is given in Subsection \ref{losungex}.

\textsc{Proof.} By Theorem \ref{cayley}, we have
$D^{red}_{F_0,\ldots,F_l}=D^{red}_{F_{\{0,\ldots,l\}}}$. Thus,
the Newton polyhedron of $D^{red}_{F_0,\ldots,F_l}$ is $|\Z^k/M|$
times smaller than the Newton polyhedron of
$D_{F_{\{0,\ldots,l\}}}$. We compute the latter one by Theorem
\ref{mainth}, and simplify the answer by Lemmas \ref{euljoin} and
\ref{chimf}(2). $\quad\Box$

\subsection{The case of branched coverings and higher additivity.}\label{snewtdiscrci3}
The assumptions of Theorem \ref{losung} can be significantly
relaxed, especially for large $l$. We illustrate this for $l=k$
(elimination theory) and for $l=k-1$ (the projection of
$F_0=\ldots=F_l=0$ onto $\T^{\tau}$ is typically a branched
covering in this case).

Let $\sigma_m(t_1,\ldots,t_l)$ be the symmetric function
$\sum_{a_1>0,\ldots,a_l>0\atop a_1+\ldots+a_l=m} t_1^{a_1}\ldots
t_l^{a_l}$ of formal variables.
\begin{lemma}[Higher additivity] $$\sigma_m(t_0+\tilde t_0,t_1,\ldots,t_l) =
\sum_{\mu=1}^{\infty}
\sigma_m(\underbrace{t_0,\ldots,t_0}_{\mu},\underbrace{\tilde
t_0,\ldots,\tilde t_0}_{\mu-1},t_1,\ldots,t_l)+$$
$$+\sigma_m(\underbrace{t_0,\ldots,t_0}_{\mu-1},\underbrace{\tilde
t_0,\ldots,\tilde
t_0}_{\mu},t_1,\ldots,t_l)+2\sigma_m(\underbrace{t_0,\ldots,t_0}_{\mu},\underbrace{\tilde
t_0,\ldots,\tilde t_0}_{\mu},t_1,\ldots,t_l).$$
\end{lemma}
Note that there are only finitely many non-zero terms (those for
$2\mu+k\leqslant m$) in the right hand side. For $m=l+1$, the
identity degenerates to $(t_0+\tilde t_0)t_1\ldots t_l =
t_0t_1\ldots t_k+\tilde t_0t_1\ldots t_l$. The proof is standard.

Let $\mathcal{M}_{\tau^{\vee}}(A_0)$ be the semigroup of all pairs
of the form $(A,\Delta)$, such that the finite set $A\subset\Z^k$
is analogous to $A_0\subset\Z^k$, the polyhedron $\Delta$ is in
$\mathcal{M}_{\tau^{\vee}}$, and its image under the projection
$\R^n\oplus\R^k\to\R^k$ equals the convex hull of $A$ (this is a
semigroup with respect to Minkowski addition of finite sets and
polyhedra).
\begin{defin} \label{hmfp} The \textit{higher mixed fiber polyhedron} is the
collection of symmetric functions
$$\HP:\underbrace{\mathcal{M}_{\tau^{\vee}}(A_0)\times\ldots\times\mathcal{M}_{\tau^{\vee}}(A_0)}_{l+1}\to\mathcal{M}_{\tau^{\vee}}(0)$$
for $l\geqslant 0$, such that
$$\mbox{ 1) } \HP\Bigl( \underbrace{(A,\Delta),\ldots,(A,\Delta)}_{l+1}\Bigr)=\sum_{A'\subset A} e^{A',A} {\dim A'\choose l} \int \Delta(A'),$$
for every $(A,\Delta)\in\mathcal{M}_{\tau^{\vee}}(A_0)$ with $\dim
A=k$ ($A'$ runs over all faces of $A$ of dimension $l$ or
greater), and $\HP\Bigl(
(A,\Delta),\ldots,(A,\Delta)\Bigr)=\tau^{\vee}$ for $\dim A<k$;
$$\mbox{ 2) } \HP(t_0+\tilde t_0,t_1,\ldots,t_l) =
\sum_{\mu=1}^{\infty}
\HP(\underbrace{t_0,\ldots,t_0}_{\mu},\underbrace{\tilde
t_0,\ldots,\tilde t_0}_{\mu-1},t_1,\ldots,t_l)+$$
$$+\HP(\underbrace{t_0,\ldots,t_0}_{\mu-1},\underbrace{\tilde
t_0,\ldots,\tilde
t_0}_{\mu},t_1,\ldots,t_l)+2\HP(\underbrace{t_0,\ldots,t_0}_{\mu},\underbrace{\tilde
t_0,\ldots,\tilde t_0}_{\mu},t_1,\ldots,t_l)$$ for all pairs
$t_0,\tilde t_0,t_1\ldots,t_l$ in
$\mathcal{M}_{\tau^{\vee}}(A_0)$.
\end{defin}
By induction on $k-l$, these conditions uniquely define the
function $\HP$ (at the base of the induction, for $l=k$, we have
the definition of the mixed fiber polyhedron). On the other hand,
by the lemma stated above, the polyhedron $\HP\Bigl(
(A_0,\Delta_0),\ldots,(A_l,\Delta_l)\Bigr)=$ $$=
\sum_{A'_0,\ldots,A'_l} \qquad e^{A'_0,A_0}
\sum_{a_0>0,\ldots,a_l>0\atop a_0+\ldots+a_l=\dim A'_0+1}
\Delta_0(A'_0)^{a_0}\cdot\ldots\cdot\Delta_l(A'_l)^{a_l},$$ with
$(A'_0,\ldots,A'_l)$ running over all collections of
corresponding faces $A'_0\subset A_0,\ldots,$ $A'_l\subset A_l$,
satisfies Definition \ref{hmfp}. In particular, $\HP=0$ for $l>k$.
We can now formulate Theorem \ref{losung} as follows.
\begin{theor} \label{thhmfp} If $A_0,\ldots,A_l$ are analogous finite sets in $\R^k$,
and $\sum_i A_i\times\{1\}$ generates $\Z^k\oplus\Z^1$, then, in
the notation of Subsection \ref{sdiscrimci}, the Newton
polyhedron of the discriminant $D^{red}_{F_0,\ldots,F_l}$ equals
$\HP\Bigl( (A_0,\Delta_0),\ldots,(A_l,\Delta_l)\Bigr)$.
\end{theor}
Unexpectedly, as soon as we formulate Theorem \ref{losung} in
this form, it can be generalized to non-analogous collections
$A_0,\ldots,A_l$ in some cases (examples are Propositions
\ref{branched0} and \ref{branched} below), which motivates the
following question:
\begin{quote} To what extent one can relax the assumption that
the arguments of $\HP$ are analogous in Definition \ref{hmfp}, so
that the higher mixed fiber polyhedron still exists and Theorem
\ref{thhmfp} remains valid?
\end{quote}

In Subsection \ref{snewtdiscrci}, we computed the Newton
polyhedron of the discriminant $D^{red}_{F_0,\ldots,F_l}$ of
functions $F_0,\ldots,F_l$, whose leading coefficients are in
general position. We denote this Newton polyhedron by
$\mathcal{N}^{A_0,\ldots,A_l}_{\Delta_0,\ldots,\Delta_l}$, where
$\Delta_0,\ldots,\Delta_l$ are the Newton polyhedra of the
functions $F_0,\ldots,F_l$. The following description of
$\mathcal{N}^{A_0,\ldots,A_k}_{\Delta_0,\ldots,\Delta_k}$ is
equivalent to Theorem \ref{losung} for $l=k$ and analogous sets
$A_0,\ldots,A_k$, but is valid for arbitrary sets
$A_0,\ldots,A_k$.
\begin{utver} \label{branched0} Suppose that $l=k$, and the lattice $\Z^k$ is generated by pairwise differences
of elements of $A_i$ for every $i=0,0',1,\ldots,k$. Then \newline
Additivity:
$\mathcal{N}^{A_0+A_{0'},A_1\ldots,A_k}_{\Delta_0+\Delta_{0'},\Delta_1\ldots,\Delta_k}=
\mathcal{N}^{A_0,\ldots,A_k}_{\Delta_0,\ldots,\Delta_k}
+\mathcal{N}^{A_{0'},A_1\ldots,A_k}_{\Delta_{0'},\Delta_1\ldots,\Delta_k}$.
\newline Unmixed case: If $\Delta_0=\ldots=\Delta_k$, and $A_0$ is not contained in a hyperplane, then
$\mathcal{N}^{A_0,\ldots,A_k}_{\Delta_0,\ldots,\Delta_k}=\int\Delta_0$.
\end{utver}
This is just another formulation of Theorem \ref{elimth}(1). We
generalize this proposition to the case $l=k-1$ as follows.
\begin{defin} Finite sets $A\subset\Z^k$ and $V\subset(\Z^k)^*$ are said to be \textit{compatible}, if
\newline 1) $V$ contains the primitive external normal covector to
every codimension 1 face of the convex hull of $A$, \newline 2)
for every linear function $v\in V$, the maximal and the next to
the maximal values of $v$ on $A$ differ by 1, \newline 3) pairwise
differences of elements of $A$ generate $\Z^k$.
\end{defin}
Nota that, if $A_1$ and $A_2$ are compatible with the same $V$, it
does not imply that $A_1$ and $A_2$ are analogous: the simplest
example is $A_1=\{(0,0),(0,1),(1,0)\}$ and
$A_2=\{(0,0),(0,-1),(-1,0)\}$.
\begin{utver} \label{branched} Suppose that $l=k-1$, and all the sets $\sum\limits_{i\in I} A_i,\,
I\subset\{0,0',1,\ldots,k-1\}$, are compatible with the same set
$V\in(\Z^k)^*$. Then \newline higher additivity:
$$\mathcal{N}^{A_0+A_{0'},A_1\ldots,A_{k-1}}_{\Delta_0+\Delta_{0'},\Delta_1\ldots,\Delta_{k-1}}=
\mathcal{N}^{A_0,\ldots,A_{k-1}}_{\Delta_0,\ldots,\Delta_{k-1}}
+\mathcal{N}^{A_{0'},A_1\ldots,A_{k-1}}_{\Delta_{0'},\Delta_1\ldots,\Delta_{k-1}}+
2\mathcal{N}^{A_0,A_{0'},A_1\ldots,A_{k-1}}_{\Delta_0,\Delta_{0'},\Delta_1\ldots,\Delta_{k-1}};$$
\newline unmixed case: If $\Delta_0=\ldots=\Delta_{k-1}$, and $A_0$ has
codimension 0, then
$$\mathcal{N}^{A_0,\ldots,A_{k-1}}_{\Delta_0,\ldots,\Delta_{k-1}}=k\int\Delta_0-\sum_{A'}\int\Delta_0(A'),$$
where $A'$ runs over all codimension 1 faces of $A_0$.
\end{utver}
In the same way as for $l=k$, these two identities are enough to
computate
$\mathcal{N}^{A_0,\ldots,A_{k-1}}_{\Delta_0,\ldots,\Delta_{k-1}}$.
The proof of Proposition \ref{branched} follows the same lines as
for Theorem \ref{losung}; the computations are not affected by the
fact that the sets $A_0,\ldots,A_{k-1}$ may be not analogous
under the assumptions above.

\begin{rem}Another way to prove additivity in Proposition \ref{branched} is to consider functions
$F_i$ with generic leading coefficients and Newton polyhedra
$\Delta_i$ for $i=0,0',1,\ldots,k-1$, and a function $F$ with
generic leading coefficients and the Newton polyhedron
$\Delta_0+\Delta_{0'}$. Then, as $F$ tends to the product
$F_0F_{0'}$, the discriminant $D_{F,F_1,\ldots,F_{k-1}}$ tends to
the product
$D_{F_0,F_1,\ldots,F_{k-1}}D_{F_{0'},F_1,\ldots,F_{k-1}}(D_{F_0,F_{0'},F_1,\ldots,F_{k-1}})^2$,
as the following picture illustrates for $l=k-1=0$.\end{rem}
\begin{center}
\noindent\includegraphics[width=12cm]{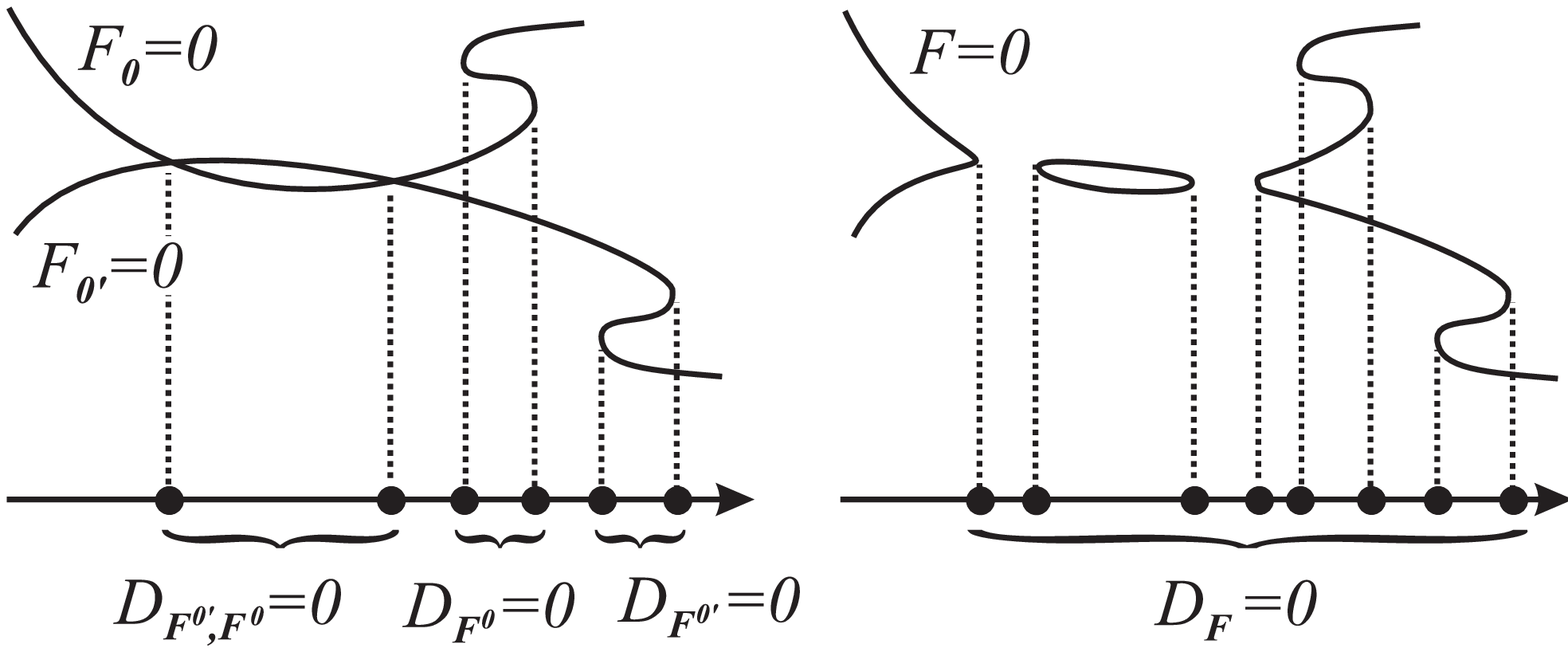} \end{center}

We omit a detailed proof of Proposition \ref{branched}, because
its only purpose is to provide a motivation for the question,
formulated after Theorem \ref{thhmfp}, but neither of the two
mentioned ideas of the proof seem relevant to answer this
question.

\subsection{Bifurcation sets and their Newton polyhedra.} \label{sfibration}

In this subsection we study one more problem, similar to the one
studied in the first part of this section. Namely, in the notation
of Section \ref{sdiscrimci}, we study the minimal (closed) subset
$S_{F_0,\ldots,F_l}\subset\T^{\tau}$, such that the restriction
of the projection $\T^{\tau}\times\CC^k\to\T^{\tau}$ to
$\{F_0=\ldots=F_l=0\}$ is a fiber bundle outside of
$S_{F_0,\ldots,F_l}$. It is called the \textit{bifurcation set}
of the projection. Note that, in contrast to
$\{D^{red}_{F_0,\ldots,F_l}=0\}$, the set $S_{F_0,\ldots,F_l}$
takes into account ``singularities at infinity" of fibers of the
projection $\{F_0=\ldots=F_l=0\}\to\T^{\tau}$. For details and
examples, see Subsection \ref{sgkzdiscrimci}  where the universal
case of this problem is studied.

We are interested in the Newton polyhedron of the equation of
$S_{F_0,\ldots,F_l}$, under the assumption that the Newton
polyhedra of $F_0,\ldots,F_l$ are given, and their leading
coefficients are in general position. If
$\Delta_0,\ldots,\Delta_l$ are the Newton polyhedra of the
functions $F_0,\ldots,F_l$, and the leading coefficients of these
functions are in general position, then we denote the Newton
polyhedron of the discriminant
$D^{red}_{\lambda_0F_0+\ldots+\lambda_lF_l}$ by
$\mathcal{N}^{A_0*\ldots*A_l}_{\Delta_0*\ldots*\Delta_l}$ (see
Theorem \ref{mainth} for its computation).
\begin{theor} \label{newtloctriv} If the collection $A_0,\ldots,A_l$ is $B$-nondegenerate (see Subsection \ref{sgkzdiscrimci}),
and the leading coefficients of the functions $f_{a,i},\; a\in
A_i,$ are in general position in the sense of Definition
\ref{genposci}, then
\newline 1) the bifurcation set $S_{F_0,\ldots,F_l}$ is a
hyersurface.
\newline 2)
Assigning appropriate positive multiplcities to the components of
the hypersurface $S_{F_0,\ldots,F_l}\subset\T^{\tau}$ outside the
maximal torus, it becomes a Cartier divisor, and the Newton
polyhedron of its equation equals
$$\sum_{A'_{j_1},\ldots,A'_{j_p}}\mathcal{N}^{A'_{j_1}*\ldots*A'_{j_p}}_{\Delta(A'_{j_1})*\ldots*\Delta(A'_{j_p})},$$ where
$A'_{j_1}\subset A_{j_1},\ldots,A'_{j_p}\subset A_{j_p}$ runs over
all collections of compatible faces that can be extended to a
collection of compatible faces $A'_0\subset
A_0,\ldots,A'_l\subset A_l$ such that $\dim\sum_{j\in J}A'_j -
\dim\sum_iA'_{j_i}\geqslant |J|-p$ for every
$J\supset\{j_1,\ldots,j_p\}$.
\end{theor}
Note that the first assumption in this statement can be omitted,
if Conjecture \ref{conj1} is valid.

The proof of Theorem \ref{newtloctriv} is based on the following
idea: define the function $H_{F_0,\ldots,F_l}$ on $\T^{\tau}$ as
$H_{F_0,\ldots,F_l}(x)=B_{A_0,\ldots,A_l}\Bigl(F_0(x,\cdot),\ldots,F_l(x,\cdot)\Bigr)$
(see Subsection \ref{sgkzdiscrimci} for the definition of the
discriminant $B_{A_0,\ldots,A_l}$). Then
$\{H_{F_0,\ldots,F_l}=0\}$ is the equation of
$S_{F_0,\ldots,F_l}$ (the proof follows the same lines as the
proof of Propositions \ref{gm11} and \ref{gm12}), and can be
expressed in terms of discriminants by Corollary \ref{bifcayley}.
\begin{sledst}\label{losung2} If
each of $A_0,\ldots,A_l$ is the set of integer lattice points in
a Delzant polytope, these Delzant polytopes have the same dual
fan, and the leading coefficients of the functions $f_{a,i},\;
a\in A_i,$ are in general position in the sense of Definition
\ref{genposci}, then $e^{A'_i,A_i}=(-1)^{\dim A_i-\dim A'_i}$ for
every face $A'_i$, and the Newton polyhedron of the equation of
$S_{F_0,\ldots,F_l}$ equals
$$\sum_{A'_0,\ldots,A'_l} \sum_{a_0>0,\ldots,a_l>0\atop
a_0+\ldots+a_l=k+1}
\Delta_0(A'_0)^{a_0}\cdot\ldots\cdot\Delta_l(A'_l)^{a_l},$$ where
$A'_0,\ldots,A'_l$ runs over all collections of compatible faces
of the sets $A_0,\ldots,A_l$, including
$A'_0=A_0,\ldots,A'_l=A_l$.
\end{sledst}
This is a corollary of Theorem \ref{newtloctriv}, Theorem
\ref{losung} and Proposition \ref{conj1sim}.

\subsection{Example and computability questions.}\label{losungex}

\textsc{Example.} Consider the first coordinate in the torus
$\CC^3$ as the ``height function", and the first coordinate line
in $\R^3$ as the ``vertical" line. We will compute the number of
critical points of the restriction of the height function to the
curve $f=g=0$ and to the surface $f=0$ for generic equations $f$
and $g$ with a given Newton polyhedron $\Delta\subset\R^3$. We
can compute these numbers by the global version of Theorems
\ref{losung} and \ref{mainth} respectively with $k=2$ and $n=1$,
because the number of critical points of the height function is
the degree of the discriminant of the projection onto the
vertical coordinate line. The answer is as follows.

Let $S_1,\ldots,S_n$ be the areas of vertical faces of $\Delta$,
let $l_1,\ldots,l_n$ be the lengths of its vertical edges, and
let $d_1,\ldots,d_n,e_1,\ldots,e_n$ be the Euler obstructions of
$\Delta$ at these faces and edges. In this notation, the number
of critical points of the restriction of the height function to
the curve $f=g=0$ equals
$$12\Vol\Delta+2\sum_i d_i S_i.$$ The number of critical points of
the restriction of the height function to the surface $f=0$ equals
$$6\Vol\Delta+2\sum_i d_i S_i+\sum_i e_i l_i.$$
To explain the first coefficient in the first of these answers
informally, note that the desired critical points are solutions to
the following system of equations:
$$f=g=\det\begin{pmatrix}
  \frac{\partial f}{\partial y} & \frac{\partial f}{\partial z} \\
  \frac{\partial g}{\partial y} & \frac{\partial g}{\partial z}
\end{pmatrix}=0.$$
The Newton polyhedron of the first two equations is denoted by
$\Delta$, thus the Newton polyhedron of the last equation
``approximately" equals $2\Delta$. Thus, if the
Kouchnirenko-Bernstein formula were applicable to this system of
the equations, then it would have approximately
$6\MV(\Delta,\Delta,2\Delta)=12\Vol(\Delta)$ solutions. Although
this illustrates why the coefficient in the first answer equals
$12$, neither the Newton polyhedron of the last solution equals
$2\Delta$ in general, nor are the equations generic with respect
to their Newton polyhedra. Thus, such a straightforward way to
count critical points would be irrelevant.

\vspace{0.3cm} \textsc{ Computability questions.} Since Euler
obstructions of polyhedra may be negative, many statements and
computations above involve subtraction of polyhedra. The
difference of polyhedra $A$ and $B$ is by definition the solution
of the equation $B+X=A$. It does not always exists (e.g. the
difference of polygons $P$ and $Q$ exists if and only if, for
every edge of $Q$, we can find a longer or equal edge of $P$ with
the same external normal). If the difference $A-B$ exists, then
it is unique by the following reason.

Recall that the \textit{support function}
$A(\cdot):(\R^m)^*\to\R\sqcup\{-\infty\}$ of a polyhedron
$A\subset\R^m$ is defined by the equality $A(l)=\min_{a\in
A}l(a)$. If the difference of support functions
$A(\cdot)-B(\cdot)$ is not concave, then the Minkowski difference
$A-B$ does not exist, otherwise $A-B$ can be reconstructed from
its support function, which equals $A(\cdot)-B(\cdot)$.

Thus, when computing Minkowski linear combinations of mixed fiber
polyhedra that appear throughout the paper, it is reasonable to
encode polyhedra with their support functions. Then Minkowski
summation and subtraction is substituted with summation and
subtraction of support functions, and mixed fiber polyhedra can
be computed by means of \cite{sy} (where the corner locus of the
support function of the mixed fiber polyhedron
$\MF(P_0,\ldots,P_k)$ is computed in terms of the corner loci of
the support functions of the arguments $P_0,\ldots,P_k$), or by
means of Proposition \ref{supportmixfib} (in Appendix, this form
of the answer is also represented as the mixed volume of certain
bounded virtual polyhedra, i.e. the tropical intersection number
of the corner loci of their support functions).

For instance, denote the \textit{support face} of the polyhedron
$P\subset\R^m$, at which a linear function $\mu\in(\R^m)^*$
attains its minimum, by $P^{\mu}$, and discuss the following
problem regarding the Newton polyhedron $\mathcal{N}\subset\R^n$
of the discriminant that was discussed in the introduction:

\begin{quote} Given the Newton polyhedra $\Delta_0,\ldots,\Delta_l$ and a linear
function $\mu\in(\R^n)^*$ that attains its minimum at a vertex of
$\mathcal{N}$, compute the coordinates of the vertex
$\mathcal{N}^{\mu}$.
\end{quote}
We restrict our attention to a coordinate function that we denote
by $\lambda$, and compute the coordinate
$\mathcal{N}^{\mu}(\lambda)$ as follows:
\newline
1) Since $(P\pm Q)^{\mu}(\lambda)=P^{\mu}(\lambda)\pm
Q^{\mu}(\lambda)$, the theorem in the introduction represents
$\mathcal{N}^{\mu}(\lambda)$ as a linear combination of values
$R^{\mu}(\lambda)$, where $R$ runs over mixed fiber polyhedra of
the form $\MF(\Delta_{i_0},\ldots,\Delta_{i_k})$;
\newline
2) Since the support function of the mixed fiber polyhedron
$\MF(P_0,\ldots,P_k)\subset\R^n$ for a linear function
$\mu\in(\R^n)^*$ equals the Minkowski sum of mixed fiber polyhedra
$\MF(P_0^{\mu'},\ldots,P_k^{\mu'})$ over all $\mu'$ whose
restriction to $\R^n$ equals $\mu$ (see \cite{mcm} or
\cite{ekh}), we can represent
$\MF(\Delta_{i_0},\ldots,\Delta_{i_k})^{\mu}(\lambda)$ as the sum
of $\MF(\Delta^{\mu'}_{i_0},\ldots,\Delta^{\mu'}_{i_k})(\lambda)$;
\newline
3) The latter value can be computed by Proposition
\ref{supportmixfib}.

In \cite{dfs}, the same problem (of finding the support vertex of
a given linear function) for the Newton polyhedron of the
discriminant $D_A$ is solved in another way, which has two
advantages: it is positive (i.e. the algorithm is based on
formulas that do not involve subtraction) and it will works for
tropicalizations of discriminant sets of higher codimension (i.e.
when $A$ is dual defect). Thus, it would be useful to generalize
the technique of \cite{dfs} to our setting.

\section{Appendix (\cite{mmj}): Mixed fiber bodies.} The
notion of the mixed fiber polytope is a natural generalization of
the mixed volume and the Minkowski integral (see \cite{bs} or
Definition \ref{def1}). It is closely related to elimination
theory, see \cite{ekh}, \cite{st}, or Theorem \ref{elimth}. The
existence of mixed fiber polytopes was predicted in \cite{mcd}
and proved in \cite{mcm}. One can extend the notion of the mixed
fiber polytope to convex bodies by continuity. We present a
direct proof of the existence of mixed fiber bodies, which does
not exploit the reduction to polytopes by continuity (see the
proof of Theorem \ref{thmixfib}). It is based on an explicit
formula \ref{ap1}.$(*)$ for the support function of a mixed fiber
body. Other applications of this formula include the proof of
Theorem \ref{elimth} and a certain monotonicity property for
mixed fiber bodies (Theorem \ref{thmonot1}). For simplicity, we
discuss bounded convex bodies here, although this restriction can
be easily omitted (see Section \ref{smixedfib} for unbounded mixed
fiber polyhedra).

\subsection{Mixed fiber bodies}\label{ap1} Let $L$ and $M$ be real vector
spaces of dimension $l$ and $m$ respectively, and let $\mu$ be a
volume form on $M$. Denote the projections of $L\oplus M$ to $L$
and $M$ by $u$ and $v$ respectively. Let $\Delta\subset L\oplus
M$ be a convex body, i.~e. a compact set, which contains all the
line segments connecting any pair of its points. For a convex
body and a point $a\in M$, denote the \textit{fiber}
$u\bigl(\Delta\cap v^{(-1)}(a)\bigr)$ of $\Delta$ by $\Delta_a$.
Recall that the \textit{support function} $B(\cdot):L^*\to\R$ of
a convex body $B\subset L$ is defined as
$B(\gamma)=\max\limits_{b\in B}\langle\gamma,b\rangle$ for every
covector $\gamma\in L^*$.
\begin{defin} \label{def1} For a convex body $\Delta\subset L\oplus M$, its \textit{Minkowski integral} is
the convex body $B\subset L$, such that its support function
equals the integral of the support functions of the fibers
$\Delta_a$, where $a$ runs over $v(\Delta)$:
$$B(\gamma)={\textstyle\int\limits}_{v(\Delta)}\Delta_a(\gamma)\, \mu\, \mbox{ for every } \gamma\in L^*.$$
The Minkowski integral is denoted by $\int \Delta\, \mu$.
\end{defin}
This definition is slightly different from the original one (see
\cite{bs}). We discuss this difference in Section \ref{ap4}.

Denote the set of all convex bodies in a real vector space $K$ by
$\mathcal{C}(K)$. This set is a semigroup with respect to the
Minkowski summation $A+B=\{a+b\, |\, a\in A, b\in B\}$.
\begin{theor} \label{thmixfib} There exists a unique symmetric Minkowski-multilinear
map $\MF_{\mu}:\underbrace{\mathcal{C}(L\oplus
M)\times\ldots\times \mathcal{C}(L\oplus M)}_{m+1}\to
\mathcal{C}(L)$, such that\newline
$\MF_{\mu}(\Delta,\ldots,\Delta)=\int\Delta\, \mu$ for every
convex body $\Delta\subset L\oplus M$.
\end{theor}

\begin{defin} The convex body $\MF_{\mu}(\Delta_0,\ldots,\Delta_m)$ is
called the \textit{mixed fiber body} of bodies
$\Delta_0,\ldots,\Delta_m$. \end{defin}

It is quite easy to see that the Minkowski integral is a
homogeneous polynomial, and, thus, admits such a polarization in
the class of virtual convex bodies (see Definition
\ref{defvirt}). The fact that this polarization gives actual
convex bodies rather than virtual convex bodies is the most
important part of the assertion.

\textsc{Proof of Theorem \ref{thmixfib}.}
\begin{defin} The \textit{shadow volume} $S_{\mu}(B)$ of a convex body
$B\subset \R\oplus M$ is the integral $\int_{B_0}\varphi\, \mu$,
where $B_0$ is the projection of $B$ to $M$, and $\varphi$ is the
maximal function on $B_0$ such that its graph is contained in $B$
(in other words, $\varphi(a)=\max_{(t,a)\in B} t$ for every $a\in
B_0$).\end{defin} One can reformulate the definition of the
Minkowski integral as follows. For a covector $\gamma\in L^*$ and
a convex body $\Delta\subset L\oplus M$, denote the image of
$\Delta$ under the projection $(\gamma,\id):L\oplus M\to\R\oplus
M$ by $\Gamma_{\Delta}(\gamma)$.
\begin{lemma}\label{reformul} The value of the support function
of the Minkowski integral $\int\Delta\, \mu$ at a covector
$\gamma\in L^*$ equals the shadow volume of the body
$\Gamma_{\Delta}(\gamma)\subset\R\oplus M$.
\end{lemma} This lemma implies that, instead of constructing mixed fiber bodies, it is enough to construct
 the mixed shadow volume in the following sense.
\begin{theor} \label{mixshad} There exists a
unique symmetric Minkowski-multilinear function $\MS_{\mu}$ of
$m+1$ convex bodies in $\R\oplus M$ such that
$\MS_{\mu}(B,\ldots,B)=S_{\mu}(B)$ for every convex body
$B\subset\R\oplus M$.
\end{theor} The proof is given in
Section \ref{ap2} and is based on an explicit formula for the
function $\MS_{\mu}$, which is later used in the proof of Theorem
\ref{thmixfib} (see Lemma \ref{lmonotapp}).
\begin{defin} The number $\MS_{\mu}(B_0,\ldots,B_m)$ is called
the \textit{mixed shadow volume} of convex bodies
$B_0,\ldots,B_m\subset\R\oplus M$.
\end{defin} Note that mixed shadow volume is a special case of mixed volume
of pairs (see Proposition \ref{msmv}).

If the existence of mixed fiber bodies is proved, then Lemma
\ref{reformul} and Theorem \ref{mixshad} imply that the value of
the support function of the mixed fiber body
$\MF_{\mu}(\Delta_0,\ldots,\Delta_m)$  at a covector $\gamma\in
L^*$ equals the mixed shadow volume
$\MS_{\mu}\bigl(\Gamma_{\Delta_0}(\gamma),\ldots,\Gamma_{\Delta_m}(\gamma)\bigr)$.
We reverse this argument, using the following fact. Recall that a
function $f:\R^k\to\R$ is said to be \textit{positively
homogeneous}, if $f(ta)=tf(a)$ for all $a\in\R^k$ and $t\geqslant
0$.
\begin{theor}\label{thconv} For any collection of convex bodies
$\Delta_0,\ldots,\Delta_m\subset L\oplus M$, the expression
$\MS_{\mu}\bigl(\Gamma_{\Delta_0}(\gamma),\ldots,\Gamma_{\Delta_m}(\gamma)\bigr)$
is a convex positively homogeneous function of a covector
$\gamma\in L^*$.
\end{theor} The proof is given in Section \ref{ap3}.

Theorem \ref{thconv} implies that, for any collection of convex
bodies $\Delta_0,\ldots,\Delta_m\subset L\oplus M$, the expression
$\MS_{\mu}\bigl(\Gamma_{\Delta_0}(\gamma),\ldots,\Gamma_{\Delta_m}(\gamma)\bigr)$
defines the support function of some convex body $B\in L$:
$$\MS_{\mu}\bigl(\Gamma_{\Delta_0}(\gamma),\ldots,\Gamma_{\Delta_m}(\gamma)\bigr)=B(\gamma)\, \mbox{ for every } \gamma\in L^*.\eqno(*)$$
This body $B$ satisfies the definition of the mixed convex body of
$\Delta_0,\ldots,\Delta_m$ by Lemma \ref{reformul} and Theorem
\ref{mixshad}. $\Box$

\subsection{Mixed shadow volume. Proof of Theorem \ref{mixshad}.}\label{ap2}
The shadow volume of a convex body $\Delta\subset\R\oplus M$ is
equal to the volume of some virtual convex body $[\Delta]$
associated with $\Delta$ (see Definition \ref{defshad} and Lemma
\ref{vssv} below). Since the correspondence $\Delta\to[\Delta]$ is
Minkowski-linear, one can define the mixed shadow volume of
convex bodies $\Delta_0,\ldots,\Delta_m\subset\R\oplus M$ as the
mixed volume of virtual bodies $[\Delta_0],\ldots,[\Delta_m]$,
which implies the existence of the mixed shadow volume. To
formulate this in detail, recall the definition of a virtual
convex body.

\textit{The Grothendieck group} $\Lambda_G$ of a commutative
semigroup $\Lambda$ with the cancellation law ($a+c=b+c\;
\Rightarrow\; a=b$) is the group of formal differences of elements
from $\Lambda$. In more detail, it is the quotient of the set
$\Lambda\times\Lambda$ by the equivalence relation
$(a,b)\sim(c,d) \Leftrightarrow a+d=b+c$, with operations
$(a,b)+(c,d)=(a+c,\, b+d)$ and $-(a,b)=(b,a)$. The map, which
carries every $a\in\Lambda$ to $(a+a,a)\in\Lambda_G$, is an
inclusion $\Lambda\hookrightarrow \Lambda_G$. An element of the
form $(a+a,a)\in \Lambda_G$ is said to be \textit{proper} and is
usually identified with $a\in \Lambda$. Under this convention,
one can write $(a,b)=a-b$.
\begin{defin} \label{defvirt} \textit{The group
of virtual bodies} in a real vector space $K$ is the Grothendieck
group of the semigroup of convex bodies in $K$ (with respect to
the operation of Minkowski summation).
\end{defin}
The classical operation of taking the mixed volume can be
extended to virtual bodies by linearity. This extension is
unique, but fails to be increasing (for example,
$\MV(-A,A)>\MV(-A,2A)$ for a convex polygon $A$). \begin{defin}
Let $\mu$ be a translation invariant volume form on a real vector
space $K$ of dimension $n$. \textit{The mixed volume} $\MV_{\mu}$
is the symmetric Minkowski-multilinear function of $n$ virtual
bodies in $K$, such that $\MV_{\mu}(\Delta,\ldots,\Delta)$ equals
the volume of $\Delta$ in the sense of the form $\mu$ for every
convex body $\Delta\subset K$.
\end{defin}
\begin{defin} For the difference
$\Delta$ of two convex bodies $\Delta_1$ and $\Delta_2$ in $K$,
the \textit{support function} $\Delta(\cdot):K^*\to\R$ is
$\Delta(\gamma)=\Delta_1(\gamma)-\Delta_2(\gamma)$.
\end{defin}
One can reformulate the definition of the group of virtual bodies
more explicitely as follows. A function $f:K\to\R$ is called
\textit{a DC function} if it can be represented as the difference
of two convex functions.
\begin{lemma} The map which carries every virtual body $\Delta$ to its support function $\Delta(\cdot)$
is an isomorphism between the group of virtual bodies in a real
vector space $K$ and the group of positively homogenious DC
functions on the dual space $K^*$.
\end{lemma}
\textsc{Proof of Theorem \ref{mixshad}.} Let $M$ be an
$m$-dimensional real vector space. Denote the ray $\{(t,0)\, |\,
t\leqslant 0\}\subset\R\oplus M$ by $l_-$, and denote the
half-space $\{(t,x)\, |\, t\geqslant a\}\subset\R\oplus M$ by
$H_a$.
\begin{defin} \label{defshad}
\textit{The shadow} $[\Delta]$ of a convex body
$\Delta\subset\R\oplus M$ is the difference of the convex bodies
$(\Delta+l_-)\cap H_a$ and $l_-\cap H_a$, where $a$ is a negative
number  such that $H_a\supset\Delta$.
\end{defin}
This definition does not depend on the choice of $a$, because one
can reformulate it in terms of support functions as follows. For
every covector $\gamma=(t,\gamma_0)\in(\R\oplus M)^*$, denote the
covector $(\max\{0,t\},\gamma_0)\in(\R\oplus M)^*$ by $[\gamma]$.
Then $[\Delta](\gamma)=\Delta([\gamma])$ for every $\gamma$.

Denote the unit volume form on $\R$ by $dt$. The function which
assigns the number
$\MV_{dt\wedge\mu}\bigl([\Delta_0],\ldots,[\Delta_m]\bigr)$ to
every collection of convex bodies $\Delta_0,\ldots,\Delta_m$ in
$\R\oplus M$ is symmetric Minkowski-multilinear by Lemma
\ref{ladd} below and assigns the shadow volume $S_{\mu}(\Delta)$
to the collection $(\Delta,\ldots,\Delta)$ for every convex body
$\Delta$ by Lemma \ref{vssv} below. Thus, it satisfies the
definition of the mixed shadow volume. Its uniqueness follows from
Lemma \ref{expl} below. $\Box$
\begin{lemma} \label{ladd} $[\Delta_1+\Delta_2]=[\Delta_1]+[\Delta_2]$.
\end{lemma}
\textsc{Proof.}
$[\Delta_1+\Delta_2](\gamma)=\Delta_1([\gamma])+\Delta_2([\gamma])=[\Delta_1](\gamma)+[\Delta_2](\gamma)$
for every covector $\gamma$. $\Box$
\begin{lemma} \label{vssv} $S_{\mu}(\Delta)=\MV_{dt\wedge\mu}([\Delta],\ldots,[\Delta]).$
\end{lemma}
\textsc{Proof.} If $\Delta\subset H_0$ in the notation of
Definition \ref{defshad}, then the shadow $[\Delta]$ equals the
convex body $(\Delta+l_-)\cap H_0$ by definition, and both
$S_{\mu}(\Delta)$ and
$\MV_{dt\wedge\mu}([\Delta],\ldots,[\Delta])$ equal the volume of
$(\Delta+l_-)\cap H_0$. One can reduce the statement of Lemma
\ref{vssv} to this special case: substitute an arbitrary body
$\Delta$ by a shifted body $\Delta+\{v\}\subset H_0$, where $v$ is
a vector of the form $(s,0)\in\R\oplus M$, and note that both
$S_{\mu}(\Delta)$ and
$\MV_{dt\wedge\mu}([\Delta],\ldots,[\Delta])$ increase by $s$
times the volume of the projection of $\Delta\subset\R\oplus M$
to $M$. The latter fact follows by the definition for the shadow
volume $S_{\mu}(\Delta)$, and follows from Lemma \ref{transf1}
for the mixed volume
$\MV_{dt\wedge\mu}([\Delta],\ldots,[\Delta])$. $\Box$
\begin{lemma}\label{transf1}
$\MV_{dt\wedge\mu}([\Delta_0+\{(s,0)\}],[\Delta_1],\ldots,[\Delta_m])=$
$$=\MV_{dt\wedge\mu}([\Delta_0],\ldots,[\Delta_m])+\frac{s}{m+1}\MV_{\mu}(B_1,\ldots,B_m),$$
where the convex body $B_j\subset M$ is the projection of
$\Delta_j\subset\R\oplus M$ to $M$ and $s$ is a positive
number.\end{lemma} \textsc{Proof.}
$\MV_{dt\wedge\mu}([\{(s,0)\}],[\Delta_1],\ldots,[\Delta_m])=\frac{s}{m+1}\MV_{\mu}(B_1,\ldots,B_m)$,
which is a corollary of the following well known formula (one can
consider $N=\R\oplus M$ and $L=\R\times\{0\}$). Let
$A_1,\ldots,A_n$ be convex bodies in an $n$-dimensional real
vector space $N$, suppose that $A_1,\ldots,A_l$ are contained in
an $l$-dimensional subspace $L\subset N$, and denote the
projection $N\to N/L$ by $p$. Then
$n!\MV_{\mu\wedge\mu'}(A_1,\ldots,A_n)=(n-l)!\MV_{\mu}(pA_{l+1},\ldots,pA_n)\cdot
l!\MV_{\mu'}(A_1,\ldots,A_l),$ where $\mu$ and $\mu'$ are volume
forms on $N/L$ and $L$. $\Box$
\begin{lemma}\label{expl} $$\MS_{\mu}(\Delta_0,\ldots,\Delta_m)=\frac{1}{(m+1)!}\sum_{0\leqslant i_1<\ldots<i_p\leqslant
m}(-1)^{m+1-p}S_{\mu}(\Delta_{i_1}+\ldots+\Delta_{i_p}).$$
\end{lemma}
\textsc{Proof.} To prove the identity $$n!
l(a_1,\ldots,a_n)=\sum_{1\leqslant i_1<\ldots<i_p\leqslant
n}(-1)^{n-p}l(a_{i_1}+\ldots+a_{i_p},\ldots,a_{i_1}+\ldots+a_{i_p})$$
for every symmetric multilinear function $l$, open the brackets
and collect like terms in the right hand side. $\Box$

The following lemma describes how the mixed shadow volume changes
under translation and dilatation of arguments along the line
$\R\times\{0\}\subset\R\oplus M$.
\begin{lemma} \label{transf} 1) Let $D_s:\R\oplus M\to\R\oplus M$
be a dilatation along $\R\times\{0\}$, i.~e. $D_s(t,x)=(st,x)$ for
all $t\in\R$ and $x\in M$. Then
$\MS_{\mu}(D_s\Delta_0,\ldots,D_s\Delta_m)=s\MS_{\mu}(\Delta_0,\ldots,\Delta_m)$
for every non-negative $s$.
\newline 2) Let $T_s:\R\oplus M\to\R\oplus M$ be a translation,
$T_s(t,x)=(t+s,x)$ for all $t\in\R$ and $x\in M$. Then
$\MS_{\mu}(T_s\Delta_0,\Delta_1,\ldots,\Delta_m)=\MS_{\mu}(\Delta_0,\ldots,\Delta_m)+\frac{s}{m+1}\MV_{\mu}(B_1,\ldots,B_m)$,
where the convex body $B_j\subset M$ is the projection of
$\Delta_j\subset\R\oplus M$ to $M$.
\end{lemma} In
particular, the mixed shadow volume is not translation
invariant.\newline \textsc{Proof.} Part 1 follows from the
definition of shadows and the equality\newline
$\MV_{dt\wedge\mu}(D_s\Delta_0,\ldots,D_s\Delta_m)=s\MV_{dt\wedge\mu}(\Delta_0,\ldots,\Delta_m)$
for convex bodies $\Delta_0,\ldots,\Delta_m$. Part 2 follows
Lemma \ref{transf1}. $\Box$

\subsection{Convexity of mixed shadow volume. Proof of Theorem
\ref{thconv}.}\label{ap3}
\begin{defin} \textit{A set with a convexity structure} is a pair
$(U, C)$, where $U$ is an arbitrary set and $C$ is an arbitrary
map $U\times(0,1)\times U\to U$.
\end{defin}
\begin{exa} The pair $(\R^k, C_{\R^k})$, where
$C_{\R^k}(u,t,v)=tu+(1-t)v$, is a set with a convexity structure.
\end{exa}
\begin{defin}
Let $(U,C)$ and $(V,D)$ be two sets with convexity structures, and
let $\leqslant$ be a partial order on $V$. A map $f:U\to V$ is
said to be \textit{a convex map from $(U,C)$ to
$(V,D,\leqslant)$,} if $f\Bigl(C(u,t,v)\Bigr) \leqslant
D\Bigl(f(u),t,f(v)\Bigr)$ for all triples $(u,t,v)\in
U\times(0,1)\times U$.
\end{defin}
\begin{exa} For a map from $(\R^k, C_{\R^k})$ to $(\R^1, C_{\R^1},
\leqslant)$, this definition coincides with the classical one.
\end{exa}
\begin{lemma}[\textmd{tautological}] \label{lcompos} If maps
$f:(U,C)\to(V,D,\leqslant)$ and $g:(V,D)\to(W,E,\leqslant)$ are
convex, and $g$ is increasing, then their composition is convex.
\end{lemma}
We apply this lemma to prove Theorem \ref{thconv} as follows: the
map which assigns the mixed shadow volume
$\MS_{\mu}\bigl(\Gamma_{\Delta_0}(\gamma),\ldots,\Gamma_{\Delta_m}(\gamma)\bigr)$
to every covector $\gamma\in L^*$, and whose convexity we wish to
prove, can be represented as a composition of simpler maps (see
the diagram $(**)$ below), whose convexity and monotonicity are
almost obvious. The proof of their convexity and monotonicity
occupies the rest of this subsection, and implies the convexity
for the map of Theorem \ref{thconv} by Lemma \ref{lcompos} (see
the end of this subsection for details).

For a convex body $B\subset M$, let $\mathcal{C}(B)$ be the set
of all convex bodies $\Delta\subset\R\oplus M$ such that the
projection of $\Delta$ to $M$ equals $B$. Introduce the shadow
and the Minkowski convexity structures $C_S$ and $C_M$ on
$\mathcal{C}(B)$ as follows. Consider convex bodies
$\Delta_i=\{(t,a)\, |\, a\in B,\,
t\in[\psi_i(a),\varphi_i(a)]\},\; i=1, 2,$ in $\mathcal{C}(B)$,
where $\varphi_i$ and $-\psi_i$ are concave functions on $B$.
Then, by definition,
$$C_S(\Delta_1,\alpha,\Delta_2)=\Bigl\{(t,a)\, |\, a\in B,\, t\in[\alpha\psi_1(a)+(1-\alpha)\psi_2(a),\, \alpha\varphi_1(a)+(1-\alpha)\varphi_2(a)]\Bigr\},$$
$$ C_M(\Delta_1,\alpha,\Delta_2)=\alpha\Delta_1+(1-\alpha)\Delta_2,\quad \Delta_1\leqslant\Delta_2 \Leftrightarrow
\varphi_1\leqslant\varphi_2.$$

For convex bodies $\Delta_0,\ldots,\Delta_m\subset L\oplus M$,
denote the projection of $\Delta_j$ to $M$ by $B_j$, and consider
the maps
\begin{multline*}{ (L^*,\,C_{L^*})
\xrightarrow{(\Gamma_{\Delta_0},\ldots,\Gamma_{\Delta_m})}
\bigl(\mathcal{C}(B_0)\times\ldots\times\mathcal{C}(B_m),\,C_S,\,\leqslant\bigr)
 \xrightarrow{(\id,\ldots,\id)}\\ \to \bigl(\mathcal{C}(B_0)\times\ldots\times\mathcal{C}(B_m),\,C_M,\,\leqslant\bigr)
 \xrightarrow{\MS_{\mu}}(\R,\,C_{\R},\,\leqslant),\qquad(**)}\end{multline*} where
the map $\id$ sends every convex body to itself.

\begin{lemma} \label{lgamma} If a convex body $B\subset M$ is the projection of a
convex body $\Delta\subset L\oplus M$, then the map
$\Gamma_{\Delta}:(L^*,\,C_{L^*})\to\bigl(\mathcal{C}(B),\,C_S,\,\leqslant\bigr)$
is convex.
\end{lemma} \textsc{Proof.} If $\dim M=0$, then the body
$\Gamma_{\Delta}(\gamma)\subset\R\oplus\{0\}$ is the segment
$[-\Delta(-\gamma),\Delta(\gamma)]$ for every $\gamma\in L^*$, and
the convexity of $\Gamma_{\Delta}$ follows from the convexity of
the support function $\Delta(\cdot)$. One can reduce the general
statement to this special case, because the convexity of the map
$\Gamma_{\Delta}:(L^*,\,C_{L^*})\to\bigl(\mathcal{C}(B),\,C_S,\,\leqslant\bigr)$
is equivalent to the convexity of the maps
$\Gamma_{\Delta_a}:(L^*,\,C_{L^*})\to\bigl(\mathcal{C}(\{0\}),\,C_S,\,\leqslant\bigr)$
for all points $a\in B$ (see the first paragraph of Appendix for
the definition of the fiber $\Delta_a$ of the body $\Delta$).
$\Box$
\begin{lemma} \label{lid} If $\Delta_1$ and
$\Delta_2$ are convex bodies from $\mathcal{C}(B)$, and
$\alpha\in(0,1)$, then $C_S(\Delta_1,\alpha,\Delta_2)\leqslant
C_M(\Delta_1,\alpha,\Delta_2)$.
\end{lemma} \textsc{Proof.} Every point $(t,a)$ of the left hand side can be represented as $\bigl(\alpha s_1+(1-\alpha)s_2,a\bigr)$, where
$(s_1,a)\in\Delta_1$ and $(s_2,a)\in\Delta_2$. Thus, it equals
$\alpha(s_1,a)+(1-\alpha)(s_2,a)$, which is contained in the right
hand side. $\Box$
\begin{lemma} \label{lmonotapp} If $\Delta_j^0\in C(B_j)$, $\Delta_j^1\in C(B_j)$ and $\Delta_j^0\leqslant\Delta_j^1$ for $j=0,\ldots,m$,
then
$\MS_{\mu}(\Delta^0_0,\ldots,\Delta^0_m)\leqslant\MS_{\mu}(\Delta^1_0,\ldots,\Delta^1_m)$.
\end{lemma}
\textsc{Proof.} Denote the convex body $\{0\}\times
B_j\subset\R\oplus M$ by $\tilde B_j$, and denote the half-space
$\R_{\geqslant 0}\times M$ by $H$. Shifting the bodies
$\Delta_j^i$ and using part 2 of Lemma \ref{transf}, we can
assume without loss of generality, that $\Delta_j^i\subset H$ for
all $i=1,2$ and $j=0,\ldots,m$. Under this assumption, the
shadows $[\Delta_j^i]$ of bodies $\Delta_j^i$ (see Definition
\ref{defshad}) are convex hulls $\conv(\Delta_j^i\cup \tilde
B_j)$. In particular, the shadows $[\Delta_j^i]$ are convex
bodies, and $[\Delta_j^0]\subset[\Delta_j^1]$. Since the mixed
shadow volume equals the mixed volume of shadows, the inequality
$\MS_{\mu}(\Delta^0_0,\ldots,\Delta^0_m)\leqslant\MS_{\mu}(\Delta^1_0,\ldots,\Delta^1_m)$
follows from the monotonicity of the mixed volume of convex
bodies. $\Box$

\textsc{Proof of Theorem \ref{thconv}.} Positive homogeneity
follows from part 1 of Lemma \ref{transf}. To prove convexity,
represent the map which assigns the mixed shadow volume
$\MS_{\mu}\bigl(\Gamma_{\Delta_0}(\gamma),\ldots,\Gamma_{\Delta_m}(\gamma)\bigr)$
to every covector $\gamma\in L^*$, as a composition of simpler
maps $(**)$. These three maps are convex and increasing.

Namely, the convexity of the map
$\Gamma_{\Delta_j}:(L^*,\,C_{L^*})\to\bigl(\mathcal{C}(B_j),\,C_S,\,\leqslant\bigr)$
is proved in Lemma \ref{lgamma}. The increasing monotonicity of
$\id:\bigl(\mathcal{C}(B_j),\,C_S,\,\leqslant\bigr)\to(\mathcal{C},\,C_M,\,\leqslant)$
is tautological, the convexity follows from Lemma \ref{lid}. The
convexity of the mixed shadow volume follows from its linearity,
and the increasing monotonicity of the mixed shadow volume is the
statement of Lemma \ref{lmonotapp}.

Since these maps are convex and increasing, their composition is
also convex by Lemma \ref{lcompos}. $\Box$

\subsection{Remarks.}\label{ap4}

\textsc{Mixed volumes of pairs.} Mixed shadow volume is a special
case of mixed volume of pairs (see Definition \ref{defrelmixvol}).
Let $p$ be the projection $\R\oplus\R^k\to\{0\}\times\R^k$, and
let $l_-$ be the ray $\{(t,0,\ldots,0)\, |\, t\leqslant
0\}\subset\R\oplus\R^k$. For a convex body
$\Delta\subset\R\oplus\R^k$, denote the pair
$\bigl(\Delta+l_-,p(\Delta)+l_-\bigr)\in\BP_{l_-}$ by
$\widetilde\Delta$.
\begin{utver} \label{msmv}
$\MS_{\mu}(\Delta_0,\ldots,\Delta_k)=\MV\bigl(\widetilde\Delta_0,\ldots,\widetilde\Delta_k\bigr)$
for every collection of convex bodies $\Delta_0,\ldots,\Delta_k$
in $\R\oplus\R^k$.
\end{utver}
\textsc{Proof.} This equality follows by definitions if
$\Delta_0=\ldots=\Delta_k$. The general statement follows from
this special case by uniqueness of the mixed shadow volume. $\Box$

\textsc{Billera-Sturmfels version of Minkowski integral.} The
original definition of the fiber integral is slightly different
from Definition \ref{def1}. Let $p:N\to K$ be a projection of an
$n$-dimensional real vector space to a $k$-dimensional one, and
let $\mu$ be a volume form on $K$.
\begin{defin}[\cite{bs}]\label{defbs} For a convex body $\Delta\subset N$, the set of all points of the form
$\int_{p(\Delta)}s\mu\in N$, where $s:p(\Delta)\to\Delta$ is a
continuous section of the projection $p$, is called the
\textit{Minkowski integral} of $\Delta$ and is denoted by
$\int^{BS}\Delta\, \mu$.
\end{defin}
Definitions \ref{def1} and \ref{defbs} are related as follows.
If, combining notation from these definitions, we assume that
$N=L\oplus M$ and $p$ is the projection $L\oplus M\to M$, then
the convex body $\int^{BS}\Delta\, \mu$ is contained in a fiber
of $p$, and $\int\Delta\, \mu$ is the image of $\int^{BS}\Delta\,
\mu$ under the projection $L\oplus M\to L$.

One can reduce Definition \ref{defbs} to Definition \ref{def1} as
well. This time, combining notation from these definitions,
suppose that $L=N$, $M=K$, and the body $\Delta^{diag}$ consists
of points $\bigl(a,p(a)\bigr)\in L\oplus M$, where $a$ runs over
all points of a convex body $\Delta\subset N$. Then
$\int^{BS}\Delta\, \mu=\int\Delta^{diag}\, \mu$. In particular,
one can denote
$\MF_{\mu}(\Delta_0^{diag},\ldots,\Delta_m^{diag})$ by
$\MF^{BS}_{\mu}(\Delta_0,\ldots,\Delta_m)$ and reformulate Theorem
\ref{thmixfib} for the Billera-Sturmfels version of mixed fiber
bodies.
\begin{theor} \label{thmixfibbs} There exists a unique symmetric
Minkowski-multilinear map
$\MF^{BS}_{\mu}:\underbrace{\mathcal{C}(N)\times\ldots\times
\mathcal{C}(N)}_{k+1}\to \mathcal{C}(N)$, such that
$\MF^{BS}_{\mu}(\Delta,\ldots,\Delta)=\int^{BS}\Delta\, \mu$ for
every convex body $\Delta\subset N$.
\end{theor}

\textsc{Monotonicity of mixed fiber bodies.} Proof of Theorem
\ref{thconv} gives the following fact as a byproduct.
\begin{theor}\label{thmonot1} In the notation of Theorem
\ref{thmixfib}, consider convex bodies
$\Delta_0,\ldots,\Delta_m,\Delta'_0,\ldots,\Delta'_m$ in the
space $L\oplus M$. If $\Delta_i\subset\Delta'_i$ and
$v(\Delta_i)=v(\Delta'_i)$ for every $i$, where $v$ is the
projection $L\oplus M\to M$, then
$\MF_{\mu}(\Delta_0,\ldots,\Delta_m)\subset\MF_{\mu}(\Delta'_0,\ldots,\Delta'_m)$.
\end{theor}
If $v(\Delta_i)\ne v(\Delta'_i)$, then the statement is not true
in general, but \linebreak
$\MF_{\mu}(\Delta_0,\ldots,\Delta_m)\subset\MF_{\mu}(\Delta'_0,\ldots,\Delta'_m)+a$
for a suitable $a\in L$ (see \cite{ekh}).


\begin{thebibliography}{9999}

\bibitem[B]{bernst}
D. N. Bernstein; The number of roots of a system of equations;
Functional Anal. Appl. 9 (1975), no. 3, 183--185.

\bibitem[BS]{bs}
L. J. Billera and B. Sturmfels, Fiber polytopes, Ann. of Math. (2)
135 (1992), no. 3, 527-549.

\bibitem[CDS]{cds}
E. Cattani, A. Dickenstein, B. Sturmfels, Rational hypergeometric
functions. Compositio Math. 128 (2001), 217-240.

\bibitem[CC]{cat}
R. Curran and E. Cattani, Restriction of A-Discriminants and Dual
Defect Toric Varieties, J. Symb. Comput. 42 (2007), 115-135.

\bibitem[D]{roc}
S. Di Rocco; Projective duality of toric manifolds and defect
polytopes; Proc. of the London Math. Soc. (3) 93 (2006), no. 1,
85-104.

\bibitem[DFS]{dfs}
A. Dickenstein, E. M. Feichtner and B. Sturmfels, Tropical
discriminants, J. Amer. Math. Soc., 20 (2007), 1111-1133.

\bibitem[E]{ernst}
L. Ernstr\"om, A Pl\"ucker formula for singular projective
varieties, Communications in algebra, 25 (1997), 2897-2901.

\bibitem[E05]{me05}
A. Esterov, Indices of 1-forms, resultants, and Newton polyhedra,
Russian Math. Surveys 60 (2005), no. 2, 352-353.

\bibitem[E06]{me06}
A. Esterov, Indices of 1-forms, intersection indices, and Newton
polyhedra, Sb. Math. 197 (2006), no. 7, 1085-1108.

\bibitem[E07]{detandnewt}
A. Esterov, Determinantal singularities and Newton polytopes,
Proc. of the Steklov inst. 259 (2007), 16-34.

\bibitem[E08]{mmj} A. Esterov, On the existence of mixed fiber bodies,
Moscow Mathematical Journal, 8 (2008), no. 3, 433-442.

\bibitem[E09]{mearx} A. Esterov, Determinantal singularities and Newton polyhedra,
arXiv:0906.5097.

\bibitem[EKh]{ekh}
A. Esterov and A. G. Khovanskii, Elimination theory and Newton
polytopes, Func. An. and Other Math., v. 2 (2008), no.1.


\bibitem[GKZ]{gkz}
I. M. Gelfand, M. M. Kapranov and A.V.Zelevinsky, Discriminants,
Resultants, and Miltidimentional Determinants, Birkh\"{a}user,
1994.

\bibitem[G]{gp}
P.~D. Gonzalez-Perez, Singularites quasi-ordinaires toriques et
polyedre de Newton du discriminant, Canad. J. Math. 52 (2000),
no. 2, 348-368.

\bibitem[H]{hall} M.Hall, Combinatorial Theory. Second Edition, John Wiley \& Sons, New
York etc., 1986.

\bibitem[Kar]{karp}
O. Karpenkov, Completely empty pyramids on integer lattices and
two-dimensional faces of multidimensional continued fractions.
Monatsh. Math. 152 (2007), no. 3, 217--249.

\bibitem[Kaz]{kaz}
B. Kazarnovskii, Truncations of systems of equations, ideals and
varieties, Izv. Math. 63 (1999), no. 3, 535--547.


\bibitem[Kh]{kh77}
A. G. Khovanskii, Newton polyhedra and the genus of complete
intersections, Func. Anal. Appl., 12 (1978), 38-46.

\bibitem[MT]{tak}
Y. Matsui, K. Takeuchi, A geometric degree formula for
$A$-discriminants and Euler obstructions of toric varieties,
arXiv:0807.3163.

\bibitem[MT2]{tak2}
Y. Matsui, K. Takeuchi, Milnor fibers over singular toric
varieties and nearby cycle sheaves, arXiv:0809.3148.

\bibitem[McD]{mcd}
J. McDonald, Fractional power series solutions for systems of
equations, Discrete Comput. Geom. 27 (2002), 501-529.

\bibitem[McM]{mcm}
P. McMullen, Mixed fibre polytopes, Discrete and Computational
Geometry, 32 (2004), 521-532.

\bibitem[O]{oka}
M. Oka; Principal zeta-function of non-degenerate complete
intersection singularity; J.~Fac.~Sci.~Univ.~Tokyo 37 (1990),
11--32.

\bibitem[S94]{sturmf}
B. Sturmfels, On the Newton polytope of the resultant, J.
Algebraic Combin. 3 (1994), no. 2, 207-236.

\bibitem[S96]{sturmf96}
B. Sturmfels, Gro\"bner bases and convex polytopes, University
Lecture Series, 8. AMS, Providence, RI, 1996.

\bibitem[ST]{st}
B. Sturmfels, J. Tevelev, Elimination theory for tropical
varieties, Math. Res. Lett. 15 (2008), 543-562.

\bibitem[SY]{sy} B. Sturmfels, J. Yu, Tropical Implicitization and Mixed Fiber
Polytopes, The IMA Volumes in Mathematics and its Applications,
Vol. 148 (Software for Algebraic Geometry), 111-131, Springer New
York, 2008.

\bibitem[T]{tev}
E. Tevelev, Compactifications of Subvarieties of Tori, Amer. J.
Math, 129, no. 4 (2007), 1087-1104.

\bibitem[V]{varch}
A. N. Varchenko; Zeta-Function of Monodromy and Newton's Diagram;
Inventiones math. 37 (1976), 253-262.

\end{thebibliography}
\end{document}